\documentclass{amsart}
\usepackage[margin=3cm]{geometry}
\usepackage{amsmath}
\usepackage{tikz-cd}
\usepackage{pinlabel}
\usepackage{bookmark}
\hypersetup{colorlinks=true,citecolor=black,urlcolor=blue,linkcolor=black}
\usepackage{amssymb}
\usepackage{extarrows}
\usepackage{mathabx}
\usepackage[all]{xy}
\usepackage{amsbsy}
\usepackage{graphicx}
\usepackage{appendix}
\usepackage{float}
\usepackage{subfigure}
\usepackage[numbers,sort&compress]{natbib}
\usepackage{amsthm}
\usepackage{geometry}
\usepackage{fancyhdr}
\usepackage{color} 
\usepackage{overpic}
\usepackage{enumerate} 
\usepackage{mathrsfs}
\usepackage{colonequals}
\usepackage{microtype}
\usepackage{hyperref}
\usepackage{diagbox}
\usepackage[thinlines]{easytable}

\newtheorem{thm}{Theorem}[section]
\newtheorem{cor}[thm]{Corollary}
\newtheorem{prop}[thm]{Proposition}
\newtheorem{lem}[thm]{Lemma}

\theoremstyle{definition}
\newtheorem{defn}[thm]{Definition}
\newtheorem{cons}[thm]{Construction}
\newtheorem{exmp}[thm]{Example}
\newtheorem{conj}[thm]{Conjecture}
\newtheorem*{fact}{Fact}

\newtheorem*{ack}{Acknowledgement}

\theoremstyle{remark}
\newtheorem{rem}[thm]{Remark}

\numberwithin{equation}{section}

\usepackage{xcolor}

\definecolor{lygreen}{HTML}{016646}

\newcommand{\beq}{\begin{equation*}\begin{aligned}}
\newcommand{\eeq}{\end{aligned}\end{equation*}}
\newcommand{\bpf}{\begin{proof}}
\newcommand{\epf}{\end{proof}}
\newcommand{\bthm}{\begin{thm}}
\newcommand{\ethm}{\end{thm}}
\newcommand{\bprop}{\begin{prop}}
\newcommand{\eprop}{\end{prop}}
\newcommand{\bcor}{\begin{cor}}
\newcommand{\ecor}{\end{cor}}
\newcommand{\blem}{\begin{lem}}
\newcommand{\elem}{\end{lem}}
\newcommand{\bdefn}{\begin{defn}}
\newcommand{\edefn}{\end{defn}}
\newcommand{\bcons}{\begin{cons}}
\newcommand{\econs}{\end{cons}}
\newcommand{\bexmp}{\begin{exmp}}
\newcommand{\eexmp}{\end{exmp}}
\newcommand{\brem}{\begin{rem}}
\newcommand{\erem}{\end{rem}}
\newcommand{\bfa}{\begin{fact}}
\newcommand{\efa}{\end{fact}}
\newcommand{\benu}{\begin{enumerate}[(1)]}
\newcommand{\eenu}{\end{enumerate}}

\newcommand{\al}{\alpha}

\newcommand{\ga}{\gamma}

\newcommand{\de}{\delta}

\newcommand{\ot}{\otimes}
\newcommand{\op}{\oplus}

\newcommand{\p}{\prime}

\newcommand{\bs}{\boldsymbol}

\newcommand{\aand}{~{\rm and}~}

\newcommand{\intg}{\mathbb{Z}}
\newcommand{\ft}{{\mathbb{F}_2}}

\newcommand{\real}{\mathbb{R}}

\newcommand{\comp}{\mathbb{C}}

\newcommand{\xra}{\xrightarrow}

\DeclareMathOperator{\im}{Im}
\DeclareMathOperator{\ke}{ker}

\DeclareMathOperator{\cone}{Cone}

\def\<#1>{\mathinner{\langle#1\rangle}}
\def\|#1|{\mathinner{\lVert #1 \rVert}}

\newcommand{\rk}{\operatorname{rank}}

\newcommand{\ptt}{\operatorname{pt}}

\newcommand{\dtop}[1]{D_{\rm top}(#1)}
\newcommand{\ti}[1]{t_2\left(#1\right)}
\newcommand{\dg}[2]{D^{\gamma}_{#1}(#2)}
\newcommand{\dsc}[2]{D^{\sharp,\comp}_{#1}(#2)}
\newcommand{\dst}[2]{D^{\sharp,\ft,0}_{#1}(#2)}
\newcommand{\dstm}[2]{D^{\sharp,\ft,\mu}_{#1}(#2)}
\newcommand{\nftp}[1]{\nu^{\sharp,\ft}_+(#1)}
\newcommand{\nftm}[1]{\nu^{\sharp,\ft}_-(#1)}

\newcommand{\suts}[2]{\mathbf{\Gamma}^{#2}_{#1}}
\newcommand{\sutg}[3]{(\mathbf{\Gamma}^{#2}_{#1},#3)}

\newcommand{\psp}[2]{\psi^{#1}_{+,#2}}
\newcommand{\psm}[2]{\psi^{#1}_{-,#2}}

\newcommand{\ish}[1]{I^\sharp  \left(#1\right)}
\newcommand{\ina}[1]{I^\natural \left(#1\right)}
\begin{document}

\title{\texorpdfstring{$2$}{2}-torsion in instanton Floer homology}

\author{Zhenkun Li}
\address{Department of Mathematics and Statistics, University of South Florida}
\curraddr{}
\email{zhenkun@usf.edu}
\thanks{}

\author{Fan Ye}
\address{Department of Mathematics, Harvard University}
\curraddr{}
\email{fanye@math.harvard.edu}
\thanks{}

\date{}
\dedicatory{}
\begin{abstract}
This paper studies the existence of $2$-torsion in instanton Floer homology with $\mathbb{Z}$ coefficients for closed $3$-manifolds and singular knots. First, we show that the non-existence of $2$-torsion in the framed instanton Floer homology $I^\sharp(S_n^3(K);\mathbb{Z})$ of any nonzero integral $n$-surgery along a knot $K$ in $S^3$ would imply that $K$ is fibered. Also, we show that $I^\sharp(S_{r}^3(K);\mathbb{Z})$ for any nontrivial $K$ with $r=1,1/2,1/4$ always has $2$-torsion. These two results indicate that the existence of $2$-torsion is expected to be a generic phenomenon for Dehn surgeries along knots. Second, we show that for genus-one knots with nontrivial Alexander polynomials and for unknotting-number-one knots, the unreduced singular instanton knot homology $I^\sharp(S^3,K;\mathbb{Z})$ always has $2$-torsion. Finally, some crucial lemmas that help us demonstrate the existence of $2$-torsion are motivated by analogous results in Heegaard Floer theory, which may be of independent interest. In particular, we show that, for a knot $K$ in $S^3$, if there is a nonzero rational number $r$ such that the dual knot $\widetilde{K}_r$ inside $S^3_r(K)$ is Floer simple, then $S^3_r(K)$ must be an L-space and $K$ must be an L-space knot.
\end{abstract}
\maketitle


\section{Introduction}

\subsection{Statement of results}\label{sec: statement of results}
In this paper, we study $2$-torsion in various versions of instanton Floer homology with $\intg$ coefficients. The first version is the framed instanton Floer homology $\ish{Y}$ for a $3$-manifold $Y$ introduced by Kronheimer and Mrowka \cite{kronheimer2011knot}. This version, over $\comp$, is conjectured to be isomorphic to $\widehat{HF}(Y)$, the hat version of the Heegaard Floer homology of $Y$, introduced by Ozsv\'ath and Szab\'o in \cite{ozsvath2004holomorphic}. Currently, with $\intg$ coefficients, there is no torsion observed for any known examples of $\widehat{HF}(Y;\intg)$; however, on the other hand, in this paper, we present abundant examples where $\ish{Y;\intg}$ has $2$-torsion, indicating that the two theories behave quite distinctly when working over $\intg$ or $\ft=\intg\slash(2\intg)$ instead of fields of characteristic $0$.

\bthm\label{thm: half surgery}
Suppose $K\subset S^3$ is a nontrivial knot ({\it i.e.}, not the unknot). Then $\ish{S^3_{r}(K);\intg}$ for $r=1,1/2,1/4$ all have $2$-torsion.
\ethm

\bthm\label{thm: non-fibered knot}
Suppose $K\subset S^3$ is a knot. If there exists $n\in\intg\backslash\{0\}$ such that $\ish{S^3_n(K);\intg}$ has no $2$-torsion, then $K$ must be a fibered knot.
\ethm

The fibered result in Theorem \ref{thm: non-fibered knot} is expected to be far from sharp. Indeed, we make the following conjecture, in which $KHI(S^3,K)=KHI(S^3,K;\comp)$ is the (sutured) instanton knot homology over $\comp$ constructed in \cite{kronheimer2010knots}, which is also conjectured to be isomorphic to the hat version of Heegaard knot Floer homology, $\widehat{HFK}(S^3,K;\comp)$.
\begin{conj}\label{conj: ribbon minimum}
	Suppose $K\subset S^3$ is a knot and $p/q\in\mathbb{Q}\backslash\{0\}$ such that $\ish{S^3_{p/q}(K);\intg}$ has no $2$-torsion. Then $S^3_{p/q}(K)$ must be an instanton L-space, {\it i.e.}, \[\dim \ish{S^3_{p/q}(K);\comp}=|H_1(S^3_{p/q}(K);\intg|=|p|,\]and $K$ is an instanton L-space knot. Consequently, we know $K$ is fibered and strongly quasi-positive and $|p/q|>2g(K)-1$ ({\it cf.} \cite[Theorem 1.15]{baldwin2019lspace} and \cite[Theorem 1.17]{LY2021large}).
\end{conj}

To the authors' knowledge, all known examples of $3$-manifolds $Y$ with no $2$-torsion in $\ish{Y;\intg}$ are instanton L-spaces ({\it cf.} Appendix \ref{sec: examples for no 2-torsion}). Moreover, if $\ish{S^3_n(K);\intg}$ has no $2$-torsion for a positive integer $n$, then we will show in Proposition \ref{prop: no 2-torsion implies Floer simple} that \[KHI(S^3_n(K),\widetilde{K}_n)\cong \ish{S^3_n(K);\comp},\]where $\widetilde{K}_n$ is the dual knot in the surgery manifold. Meanwhile, we establish the following theorem in Heegaard Floer theory via the immersed curve invariant \cite{HRW2024curve1,HRW2022curve2}, which is currently not available in the instanton theory and of independent interest.

\bthm\label{thm: Heegaard dual simple}
Suppose $K\subset S^3$. If there exists $p/q\in \mathbb{Q}\backslash\{0\}$ satisfying \begin{equation*}\label{eq: HF simple}
    \dim \widehat{HFK}(S^3_{p/q}(K),\widetilde{K}_{p/q};\ft)=\dim \widehat{HF}(S^3_{p/q}(K);\ft),
\end{equation*}then both dimensions must equal to $|H_1(S^3_{p/q}(K);\intg)|=|p|$ and hence $K$ is a (Heegaard Floer) L-space knot. Consequently, we know $K$ is fibered and strongly quasi-positive \cite{Ozsvath2005,ni2007knot,hedden2010positivity}, and $|p/q|>2g(K)-1$ \cite[Corollary 3.6]{Rasmussen2017}.
\ethm

In addition to the above theorems that confirm the existence of $2$-torsion for closed $3$-manifold case, we also prove the following results for the unreduced singular instanton knot homology of a knot $K$ inside a $3$-manifold $Y$, which is denoted by $\ish{Y,K}$ and was constructed by Kronheimer-Mrowka \cite{kronheimer2011khovanov}.

\bprop\label{prop: genus-one knot has 2-torsion}
If $K\subset S^3$ is a knot of genus one and with Alexander polynomial $\Delta_K(t)\neq 1$, then $\ish{S^3,K;\intg}$ has $2$-torsion.
\eprop

\bprop\label{prop: unknotting number one}
If $K\subset S^3$ is a knot with unknotting number one, then $\ish{S^3,K;\intg}$ has $2$-torsion.
\eprop

\bprop\label{prop: dual knot has 2-torsion}
If $K\subset S^3$ is a non-trivial knot, then $\ish{S^3_1(K),\widetilde{K}_1;\intg}$ has $2$-torsion.
\eprop

Theorems \ref{thm: half surgery} and \ref{thm: non-fibered knot}, together with Propositions \ref{prop: genus-one knot has 2-torsion}, \ref{prop: unknotting number one} and \ref{prop: dual knot has 2-torsion}, collectively suggest that the presence of 2-torsion in instanton Floer homology constitutes a prevalent phenomenon. Note that, in addition to the above new results about the $2$-torsion for unreduced singular instanton Floer homology, it has already been known to experts that $\ish{S^3,K;\intg}$ for quasi-alternating knots has $2$-torsion. See Appendix \ref{sec: alternating knot}.

While the homology groups associated to various coefficients are interconnected through the universal coefficient theorem, these groups can exhibit distinctive behaviors: homology groups over $\comp$ are evidently integrated into a broader axiomatic framework of Floer theory; however, homology groups over $\ft$ become more wild and may demonstrate unique characteristics, thereby offering potential novel applications.

One key property that characterizes this distinction between instanton Floer homology over $\comp$ and $\ft$ is the adjunction inequality, stated below specific to framed instanton homology over $\comp$ for simplicity.

\bthm[{\cite[Theorem 1.16]{baldwin2019lspace}}]\label{thm: adjunction}
Suppose $X:Y_1\to Y_2$ is an oriented cobordism and there exists an embedded surface $\Sigma\subset X$ of genus at least one such that
\[
\Sigma\cdot\Sigma\geq2g(\Sigma)-1.
\]
Then $\ish{X}=0: \ish{Y_1;\comp}\to \ish{Y_2;\comp}$.
\ethm

The adjunction inequality is a powerful tool to regulate the behavior of framed instanton homology. For example, Baldwin-Sivek's dimension formula \cite[Theorem 1.2]{baldwin2020concordance} for framed instanton homology of Dehn surgeries of knots in $S^3$ is almost simply built on this theorem together with the surgery exact triangle. Failure of Theorem \ref{thm: adjunction} over $\ft$ can be concluded from the following two facts.

\begin{lem}\label{lem: false statement for Poincare sphere}
	If Theorem \ref{thm: adjunction} were true over $\ft$ as well, then the Poincar\'{e} sphere $P=S^3_1(3_1)$ would satisfy the following condition: $\ish{P;\ft}\cong \ft$.\footnote{For the completeness of the paper, we state this result and present a proof, but we claim no originality of this fact. The authors first heard about this fact from John A. Baldwin.}
\end{lem}
\bpf
Note that the trefoil knot $K=3_1$ has a genus-one Seifert surface. When performing integral $n$-surgery along the trefoil, we obtain a trace cobordism
\[X_n:S^3\to S^3_n(K).\]
The Seifert surface gets closed up by the core of the $2$-handle in $X_n$, and we obtain a closed oriented genus-one surface $\Sigma_n\subset X_n$ with self-intersection number $n$. When $n\geq 1$, the adjunction inequality then implies that $X_n$ induces the zero map. The map induced by $X_n$ fits into an exact triangle ({\it cf.} Lemma \ref{lem: surgery exact triangle, BS}) in which the third term is $\ish{S^3_{n+1}(K);\ft}$. Since $\ish{S^3;\ft}\cong \ft$, we conclude that
\[
\dim\ish{S^3_1(K);\ft}=\dim\ish{S^3_2(K);\ft} -1 =\cdots =\dim\ish{S^3_5(K);\ft} - 4.
\]
Since $S^3_5(K)$ is a lens space, we know from Appendix \ref{sec: examples for no 2-torsion} that $\dim\ish{S^3_5(K);\ft} = 5$ and hence \[\dim\ish{P;\ft} = \dim \ish{S^3_1(K);\ft} = 1.\]
\epf
\begin{prop}[{Bhat \cite[Theorem 1.4]{bhat2023newtriangle}}]\label{prop: wrong for trefoil}
	Let $P$ denote the Poincar\'{e} sphere as above. We have
	\[
	\dim\ish{P;\ft}\geq 3.
	\]
\end{prop}

\begin{figure}[ht]
	\centering
	\begin{minipage}{0.49\linewidth}
		\begin{tikzpicture}[scale=0.75]
				\draw[very thick,->] (-3.5,0) -- (5.5,0) node[anchor=north] {$n$};
    				\draw[gray,very thin] (-3,-1) grid (5,5);
 					
    				\draw[blue, very thick] (-3,5) -- (1,1);
    				\draw[blue, very thick] (1,1) -- (5,5);
    				
    				\filldraw[blue] (1,1) circle (4pt);
    				
    				\draw[black, very thick, dashed] (1,1) -- (1,0);
    				
    				\node at (1,-0.5) {$\nftm{K}=\nftp{K}$};
		\end{tikzpicture}
	\end{minipage}
	\begin{minipage}{0.49\linewidth}
		\begin{tikzpicture}[scale=0.75]
				\draw[very thick,->] (-3.5,0) -- (5.5,0) node[anchor=north] {$n$};
    				\draw[gray,very thin] (-3,-1) grid (5,5);
 					
    				\draw[blue, very thick] (-3,4) -- (0,1);
    				\draw[blue, very thick] (2,1) -- (5,4);
    				\draw[blue, very thick] (0,1) -- (1,2);
    				\draw[blue, very thick] (1,2) -- (2,1);
    				
    				
    				\filldraw[blue] (0,1) circle (4pt);
    				\filldraw[blue] (2,1) circle (4pt);
    				\filldraw[blue] (1,2) circle (4pt);
    				
    				
    				\draw[black, very thick, dashed] (0,1) -- (0,0);
    				\draw[black, very thick, dashed] (2,1) -- (2,0);
    				
    				\node at (0,-0.5) {$\nftm{K}$};
    				\node at (2,-0.5) {$\nftp{K}$};
		\end{tikzpicture}
		\end{minipage}
	
	\caption{Left: V-shape. Right: W-shape.}\label{fig: V-shape and W-shape}
\end{figure}
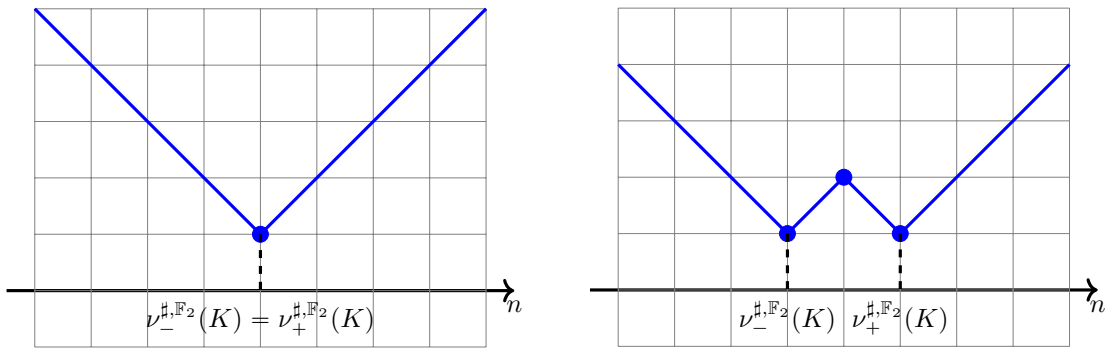

As mentioned above, the failure of the adjunction inequality implies that the dimension formula in \cite[Theorem 1.2]{baldwin2020concordance} no longer holds true for $\ft$. Using some partial vanishing result to replace the role played by the adjunction inequality, we can describe the behavior of the sequence $\{\dim\ish{S^3_n(K);\ft}\}_{n\in\intg}$ for a knot $K\subset S^3$ as in the following proposition.
\begin{prop}\label{prop: shape of dst, intro}
	Suppose $K\subset S^3$ is a knot, then the sequence $\{\dim\ish{S^3_n(K);\ft}\}_{n\in\intg}$ has one of the following three shapes.
	\begin{itemize}
		\item V-shape, see Figure \ref{fig: V-shape and W-shape}.
		\item W-shape, see Figure \ref{fig: V-shape and W-shape}.
		\item Generalized W-shape, see Figure \ref{fig: generalized W-shape, intro}.
	\end{itemize}
\end{prop}

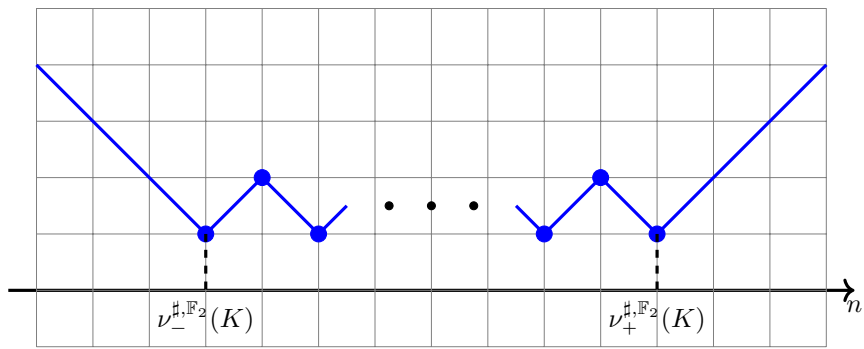
\begin{figure}[ht]
	\centering
		\begin{tikzpicture}[scale=0.75]
				\draw[very thick,->] (-7.5,0) -- (7.5,0) node[anchor=north] {$n$};
    				\draw[gray,very thin] (-7,-1) grid (7,5);
 					
    				\draw[blue, very thick] (-7,4) -- (-4,1);
    				\draw[blue, very thick] (7,4) -- (4,1);
    				\draw[blue, very thick] (-4,1) -- (-3,2);
    				\draw[blue, very thick] (4,1) -- (3,2);
    				\draw[blue, very thick] (-3,2) -- (-2,1);
    				\draw[blue, very thick] (3,2) -- (2,1);
    				\draw[blue, very thick] (-2,1) -- (-1.5,1.5);
    				\draw[blue, very thick] (2,1) -- (1.5,1.5);
    				

    				\filldraw[blue] (-4,1) circle (4pt);
    				\filldraw[blue] (-3,2) circle (4pt);
    				\filldraw[blue] (3,2) circle (4pt);
    				\filldraw[blue] (4,1) circle (4pt);
    				\filldraw[blue] (-2,1) circle (4pt);
    				\filldraw[blue] (2,1) circle (4pt);
    				
    				
    				\filldraw[black] (0,1.5) circle (2pt);
    				\filldraw[black] (-0.75,1.5) circle (2pt);
    				\filldraw[black] (0.75,1.5) circle (2pt);
    				
    				\draw[black, very thick, dashed] (-4,1) -- (-4,0);
    				\draw[black, very thick, dashed] (4,1) -- (4,0);
    				
    				\node at (-4,-0.5) {$\nftm{K}$};
    				\node at (4,-0.5) {$\nftp{K}$};
		\end{tikzpicture}
	
	\caption{Generalized W-shape.}\label{fig: generalized W-shape, intro}
\end{figure}

As shown in Figure \ref{fig: V-shape and W-shape} and Figure \ref{fig: generalized W-shape, intro}, we define two invariants:
\[
\nftp{K} = \min\{n~|~\forall k>n,~\dim\ish{S^3_{k}(K);\ft} = \dim\ish{S^3_{k-1}(K);\ft}+1\}
\]
\[
\aand \nftm{K} = \max\{n~|~\forall k<n,~\dim\ish{S^3_{k}(K);\ft} = \dim\ish{S^3_{k+1}(K);\ft}+1\}.
\]
We have the following.
\bprop\label{prop: nftp is concordance inv, intro}
The invariants $\nu^{\sharp,\ft}_{\pm}(K)$ are concordance invariants.
\eprop
\brem\label{rem: width}
Indeed, the three cases in Proposition \ref{prop: shape of dst, intro} can all be regarded as generalized W-shapes. If we introduce the \emph{width} of a generalized W-shape to be \[w(K)=\frac{\nftp{K}-\nftm{K}}{2},\] then the V-shape has width $0$ and W-shape has width $1$. By Proposition \ref{prop: nftp is concordance inv, intro}, the width of a knot is also a concordance invariant. A direct computation will show that the unknot has $\nu^{\sharp,\ft}_{\pm}=\pm 1$ and width $1$.
\erem
One interesting observation about these two invariants is that, since Theorem \ref{thm: adjunction} fails over $\ft$, the usual adjunction inequality that is expected for some gauge theoretic concordance invariants ({\it e.g.}, the cousin $\nu^{\sharp}$ invariant in \cite{baldwin2020concordance} and the $\tau$ invariant in \cite{ozsvath2003tau}) does not necessarily hold for $\nu^{\sharp,\ft}_{\pm}(K)$. This means that these invariants have the potential to lead to applications beyond the limitation of adjunction inequality.

\subsection{Strategy of the proof}\label{subsec: Strategy of the proof}
Our approach to proving the existence of $2$-torsion for instanton Floer homology over $\intg$ has a few ingredients involving the computation of the corresponding homology group over $\comp$ and $\ft$. To quantify the existence of $2$-torsion, for a (closed connected oriented) $3$-manifold $Y$, we define
\begin{equation}\label{eq: number of 2-torsion}
	\ti{Y} = \frac{1}{2}(\dim\ish{Y;\ft} - \dim\ish{Y;\comp})
\end{equation}
The universal coefficient theorem implies that $\ti{Y}$ is always a non-negative integer, and $\ish{Y,\intg}$ has $2$-torsion if and only if $\ti{Y}> 0$.

To estimate $\ti{Y}$, we need to estimate the dimensions of $\ish{Y}$ over $\ft$ and $\comp$, respectively. The estimation over $\ft$ is based on a new exact triangle recently developed by Bhat \cite{bhat2023newtriangle} in his thesis. The estimation over $\comp$ utilizes the dimension formula from Baldwin-Sivek \cite{baldwin2020concordance} and the surgery formulae \cite{LY2021large,LY2022integral1,LY2022integral2} proved by the authors of the current paper. The connection between different coefficients is from Kronheimer-Mrowka's results about the reduced version of singular instanton knot homology $\ina{Y,K}$. More precisely, from \cite[Lemma 7.7]{kronheimer2019web1} (see also \cite[Corollary 8.7]{daemi2019equivariant}) and \cite[Proposition 1.4]{kronheimer2011khovanov}, for any knot $K$ in a $3$-manifold $Y$, we have \begin{equation}\label{eq: 2dim}
    \dim \ish{Y,K;\ft}=2\dim \ina{Y,K;\ft} \aand \dim\ina{Y,K;\comp}=\dim KHI(Y,K)
\end{equation}
Hence the universal coefficient theorem implies
\begin{equation}\label{eq: KHI inequality}
    \dim \ish{Y,K;\ft}\ge 2\dim KHI(Y,K),
\end{equation}where the equality holds if and only if $\ina{Y,K}$ has no $2$-torsion. 

As mentioned above, we also need to use the (sutured) instanton knot homology $KHI(Y,K)$. It is a special case of sutured instanton homology $SHI(M,\ga)$ for a balanced sutured manifold $(M,\ga)$ \cite{kronheimer2010knots}, which is always defined over $\comp$, so we omit the coefficients for simplicity. If $K$ is null-homologous, there is a $\intg$-grading on $KHI(Y,K)$ \cite{kronheimer2010instanton}, usually called the \emph{Alexander grading}. When $K\subset S^3$, Kronheimer-Mrowka \cite[Proposition 7.16]{kronheimer2010knots} showed that the maximal nontrivial grading of $KHI(S^3,K)$ is the genus $g(K)$. Hence, we introduce the following notation.
\begin{equation}\label{eq: D_top}
	\dtop{K} = \dim KHI(S^3,K,g(K))\ge 1.
\end{equation}

For a knot $K\subset S^3$, we use $\widetilde{K}_1\subset S^3_1(K)$ to denote the dual knot of the $1$-surgery. We will prove the following crucial technical results in \S \ref{sec: Floer simple and fiberedness}.

\blem\label{lem: dual knot of 1-surgery}
Suppose $K\subset S^3$ is a nontrivial knot, then 
\[
\dim KHI(S^3_1(K),\widetilde{K}_1) \geq \dim\ish{S^3_1(K)} + 2  \dtop{K}.
\]
\elem

\blem\label{lem: fiberness detection}
Suppose $K\subset S^3$ is a nontrivial knot. Suppose $n > 0$ is an integer such that
\begin{equation}\label{eq: floer simple knots}
	\dim KHI(S^3_n(K),\widetilde{K}_n) = \dim\ish{S^3_n(K)}.
\end{equation}
Then $K$ must be fibered, $V$-shaped, and $0<\nu^{\sharp}(K)= 2  \tau^{\sharp}(K) - 1<n$, where the notions are from \cite{baldwin2020concordance}. 
\elem

It is worth mentioning that the proofs of Lemma \ref{lem: dual knot of 1-surgery} and Lemma \ref{lem: fiberness detection} (as crucial technical inputs to approach the main results of the paper) and our proposal of Conjecture \ref{conj: ribbon minimum} are all motivated by the analogous questions and results from Heegaard Floer theory, and to establish these results in Heegaard Floer theory, we also make essential use of the immersed curve invariant techniques developed in \cite{HRW2024curve1,HRW2022curve2}. We include discussion on Heegaard Floer theory in \S \ref{sec: HF theory}.

A knot $\widetilde{K}_n$ satisfying (\ref{eq: floer simple knots}) is usually called \emph{Floer-simple} \cite[\S 1.2]{Hedden2011simple}. Hence, we call $K\subset S^3$ \emph{dually Floer simple} ({\it cf.} Definition \ref{defn: dually Floer simple}) if one of its dual knots $\widetilde{K}_n$ is Floer simple. The fiberedness conclusion is to the satisfaction of the applications in this paper, but we expect more restrictions to the dually Floer simple knot ({\it cf.} Theorem \ref{thm: Heegaard dual simple}). Currently known examples of dually Floer simple knots are only instanton $L$-space knots whose properties were studied in \cite{baldwin2019lspace} and \cite{LY2021large}. 


With the help of Lemmas \ref{lem: dual knot of 1-surgery} and \ref{lem: fiberness detection}, we can bound $\ti{S^3_r(K)}$ for various surgery slopes $r=p/q$ from below. Simply put, the gap between $\dim KHI(S^3_n(K),\widetilde{K}_n)$ and $\dim\ish{S^3_n(K)}$ characterizes the non-vanishing of $\ti{S^3_r(K)}$.

Finally, we mention the following dimension inequality, which directly implies Proposition \ref{prop: unknotting number one}. It is of independent interest and provides a new obstruction to unknotting-number-one knots. The proof is a nice application of the large surgery formula developed in \cite{LY2021large}.

\bthm\label{thm: unknotting number one knot}
If $J\subset S^3$ is a knot with unknotting number one, then
\[\dim\ish{S^3,J;\comp} \leq \dim \ina{S^3,J;\comp} + 3=\dim KHI(S^3,J)+3.\]
\ethm
 \bpf[Proof of Proposition \ref{prop: unknotting number one}]
 From Theorem \ref{thm: unknotting number one knot} and the inequality (\ref{eq: KHI inequality}), if $\dim KHI(S^3,K) > 3$ then \[\dim\ish{S^3,K;\ft}-\dim\ish{S^3,K;\comp}\ge 2\dim KHI(S^3,K)-(\dim KHI(S^3,K)+3)>0\]and $\ish{S^3,K;\intg}$ has $2$-torsion. From \cite[Theorem 1.6]{baldwin2018khovanov}, we know knots with $\dim \ina{S^3,K;\comp} \leq 3$ must be either the unknot or the trefoil $3_1$. The unknot has unknotting number zero and it is already known from Theorem \ref{thm: thin knot has 2 torsion} that $\ish{S^3,K;\intg}$ has $2$-torsion for $K=3_1$.
 \epf
\begin{ack}
The authors thank Deeparaj Bhat for the valuable discussion in their joint work (under preparation). A preliminary version of Theorem \ref{thm: half surgery} was established under the discussion with Bhat. The second author thanks Yi Ni for helpful communication, which motivates Lemma \ref{lem: dual knot of 1-surgery} and Theorem \ref{thm: Heegaard dual simple}. The authors thank Steven Sivek to pointing out the case of $1/4$-surgery in Theorem \ref{thm: half surgery}. The authors also thank Peter Kronheimer and Jake Rasmussen for helpful conversations. 
\end{ack}

\section{Preliminaries}
\subsection{Basic setups}\label{subsec: basic setups}
Unless otherwise mentioned, we will always let $Y$ be a closed connected oriented $3$-manifold and $K$ to be a framed knot in $Y$. If $K$ is null-homologous, then the framing by default comes from any Seifert surface of $K$. If $Y=S^3$, we may omit it in the notion for various Floer homologies. We will use $\mathbb{K}$ to denote either a commutative ring with a unit or a field. We write $\mathbb{N}=\{0,1,\dots\}$ for natural numbers and $\ft=\intg\slash(2\intg)$.

For any smooth unoriented $1$-submanifold (called a \emph{1-cycle} or the \emph{bundle data}) $\lambda\subset Y$ and any basepoint $\ptt$ in $Y$, Kronheimer-Mrowka constructed the framed instanton homology $\ish{Y,\lambda}$ in the following three ways.
\benu
  \item From \cite[\S 4.3]{kronheimer2011khovanov}, they considered the singular instanton homology $I(Y,H,\lambda\sqcup {\lambda_0})$, where $H$ is a Hopf link in the neighborhood of $\ptt$, and $\lambda_0$ is a standard arc connecting two components of $H$.
  \item From \cite[\S 1.2]{kronheimer2011knot} (see also \cite[\S 5.4]{kronheimer2011khovanov}), they considered the ordinary instanton Floer homology $I(Y\# T^3,\lambda\sqcup S^1)^{\psi}$, where the connected sum is made at $\ptt$, $S^1$ is a standard circle in $T^3$, and $\psi\subset H^1(Y\# T^3;\ft)$ is the 2-element subgroup generated by the dual torus $R$ of $S^1$ in $T^3$, and we take the quotient of the usual configuration space by this extra $2$-element group when constructing instanton Floer homology.
  \item From \cite[\S 7.6]{kronheimer2010knots}, they considered the sutured instanton knot homology $KHI(Y,U,\lambda)$ of an unknot $U$ in the neighborhood of $\ptt$, which is ultimately defined to be the sutured instanton homology $SHI(Y\backslash  N(U),\ga_U=\mu_U\cup(-\mu_U),\lambda)$ for an oriented meridian $\mu_U$ of $U$ (see more discussion in \S \ref{sec: sutured theory}). This is only defined over $\comp$ and on the homology level since the construction uses the generalized eigenspace decomposition of some operators on the ordinary instanton Floer homology (with extra bundle data from $\lambda$). As an isomorphism class in \cite[\S 7]{kronheimer2010knots} or as a projectively transitive system by Baldwin-Sivek's naturality result \cite[\S 9]{baldwin2015naturality}, the space $KHI(Y,U,\lambda)$ is identical to $SHI(Y\backslash N(\ptt),S^1,\lambda)$, which is another definition of $\ish{Y,\lambda;\comp}$.
\eenu 
Note that the first two definitions are over $\intg$ and can be generalized to a $3$-manifold $Y$ with a singular knot $K$, which we denote by $\ish{Y,K,\lambda}$. The third definition can also be generalized to such a $(Y,K)$ by the work of Xie-Zhang \cite{xie2019tangle}, again over $\comp$. The first two definitions have $\intg/4$ homological grading, while the third only has a canonical $\intg/2$ homological grading from \cite[\S 2.6]{kronheimer2010instanton}.

Note that the isomorphism class of the Floer homology only depends on the homology class $[\lambda]\in H_1(Y;\ft)$ since it represents the bundle data for an $SO(3)$ bundle. When $\lambda=\emptyset$, we simply write $\ish{Y}$. 

From \cite[\S 5.4]{kronheimer2011khovanov}, the first two definitions are isomorphic via Floer's excision cobordism, hence the isomorphism intertwines with any cobordism maps supported in $Y\backslash N(\ptt)$ naturally. Also from \cite[\S 5.4]{kronheimer2011khovanov}, the second definition is isomorphic to the third definition via a \emph{special} choice of closure in the construction of $KHI(Y,U,\lambda)$, which is exactly $(Y\# T^3,R,\lambda\sqcup S^1)$. Since $g(R)=1$, by \cite[Theorem 2.5]{baldwin2019lspace} (see also \cite[Lemma 4]{froyshov02equi}), we know $\mu(\ptt)^2=4$ for the $\mu$ map in the construction of $KHI(Y,U,\lambda)$, and then the isomorphism can be made to intertwine with the cobordism map. Note that Baldwin-Sivek's naturality result \cite{baldwin2015naturality} only works for closures of genus larger than one, while the relation between closures of genus one and larger than one remains to be isomorphism rather than \emph{canonical} isomorphism. Hence, the best result we can state about the cobordism map is that, for any \emph{fixed} closure and any \emph{fixed} isomorphism between the special closure above and the fixed closure, the isomorphism intertwines the cobordism map supported in $Y\backslash N(U)$. Similar isomorphism results hold when there is a singular knot $K$ by Floer's excision theorem (for sutured version, see \cite[Remark 7.8]{xie2019tangle}). All isomorphisms above respect the homological gradings.
\subsection{The exact triangles}
We will use various exact triangles to establish the main results of this paper. We start with the surgery exact triangle. 

\blem[\cite{scaduto2015instanton}]\label{lem: surgery exact triangle, original}
Suppose $K\subset Y$ is a framed knot in a $3$-manifold $Y$ and suppose $\lambda\subset Y$ is a $1$-cycle. Let $\mu\subset Y$ be a push-off of the meridian of $K$. Then there is an exact triangle
\[
\xymatrix{
\ish{Y_0(K),\lambda\cup\mu;\mathbb{K}}\ar[rr]&&\ish{Y_1(K),\lambda;\mathbb{K}}\ar[dl]\\
&\ish{Y,\lambda;\mathbb{K}}\ar[ul]&
}
\]
\elem
Note that the three $3$-manifolds in the triangle in Lemma \ref{lem: surgery exact triangle, original} are cyclic, which means that we can choose any of the three manifolds as the one we start with, and obtain different bundle data in \cite[Figure 1]{scaduto2015instanton}. Alternatively, one could add the extra bundle data $\mu$ to all of the three $3$-manifolds (not just restricted to $0$-surgery). Also, since the (isomorphism class of) Floer homology only depends on the homology class of $\lambda$ mod $2$, one could start with a possibly different $1$-cycle $\lambda^{\prime}$ representing the same homology class as $\lambda$ on $Y$ and $Y_1(K)$ but cancelling $\mu$ on $Y_0(K)$. Hence, we obtain a variation of Lemma \ref{lem: surgery exact triangle, original}.
\blem[{\cite[\S 2.2]{baldwin2020concordance}}]\label{lem: surgery exact triangle, BS}
Suppose $K\subset Y$ is a framed knot in a $3$-manifold $Y$ and suppose $\lambda\subset Y$ is a $1$-cycle. Then there is an exact triangle
\[
\xymatrix{
\ish{Y_0(K),\lambda;\mathbb{K}}\ar[rr]&&\ish{Y_1(K),\lambda;\mathbb{K}}\ar[dl]\\
&\ish{Y,\lambda;\mathbb{K}}\ar[ul]&
}
\]
\elem

In addition to the above triangles, there is also an exact triangle relating the framed instanton homology of closed $3$-manifolds to singular instanton knot homology of knots.
\blem[{Bhat \cite[Theorem 1.1]{bhat2023newtriangle}}]\label{lem: mixed triangle}
Suppose $K\subset Y$ is a framed knot in a $3$-manifold $Y$ and suppose $\lambda\subset Y$ is a $1$-cycle. Then there is an exact triangle
\[
\xymatrix{
\ish{Y_0(K),\lambda;\mathbb{K}}\ar[rr]&&\ish{Y_2(K),\lambda;\mathbb{K}}\ar[dl]\\
&\ish{Y,K,\lambda;\mathbb{K}}\ar[ul]&
}
\]
\elem
\brem
For the proofs of Theorems \ref{thm: half surgery} and \ref{thm: non-fibered knot}, we use the triangle in (\ref{lem: mixed triangle}) over $\mathbb{K}=\ft$, while the proofs of Propositions \ref{prop: genus-one knot has 2-torsion} and \ref{prop: dual knot has 2-torsion} use the triangle over $\mathbb{K}=\comp$. Additionally, the triangle over $\comp$ will be used in the joint work of Bhat and the authors (in preparation). 
\erem



Here $\ish{Y,K,\lambda;\mathbb{K}}$ is the unreduced singular instanton knot homology of the triple $(Y,K,\lambda)$ introduced by Kronheimer-Mrowka \cite[\S 4.3]{kronheimer2011khovanov} and mentioned in \S \ref{subsec: basic setups}. They also introduced a reduced version, which is denoted by $\ina{Y,K,\lambda}$. More precisely, we write $K^\natural$ for the connected sum of $K$ and the Hopf link $H\subset S^3$ with a standard arc $\lambda_0$ connecting two components of $H$ as in the first definition of $\ish{Y,\lambda}$ in \S \ref{subsec: basic setups}. Then, we define \begin{equation}\label{eq: defn of ina}
    \ina{Y,K,\lambda}=I(Y,K^{\natural},\lambda\sqcup \lambda_0).
\end{equation} 


\subsection{The sutured theory}\label{sec: sutured theory}
For a balanced sutured manifold $(M,\gamma)$, Kronheimer-Mrowka \cite{kronheimer2011knot} introduced the sutured instanton homology $SHI(M,\gamma)$. Note that this homology is only defined over $\comp$ and not $\ft$ nor $\intg$. The balanced sutured manifolds we will use in this paper all come from knot complements. For a knot $K\subset S^3$, we equip the knot with the Seifert framing. Let $\Gamma_n$ be the suture on $\partial (S^3\backslash N(K))\cong T^2$ consisting of a pair of oppositely oriented non-separating simple closed curves of slope $-n$ (adding the minus sign is to be consistent with the notations in the authors' previous work). Let $\Gamma_{\mu}=\mu\cup(-\mu)$ consist of two meridians with opposite orientations. As defined in \cite[\S 7.6]{kronheimer2011knot}, we have
\[
KHI(S^3,K) = SHI(S^3\backslash N(K),\Gamma_{\mu})\aand KHI(S^3_n(K),\widetilde{K}_n) = SHI(S^3\backslash N(K),\Gamma_{-n}).
\]
We will need results from \cite[\S 5]{li2019tau}. Note that the symmetry \cite[Theorem 1.2 (3)]{li2018gluing} between a knot $K$ and its mirror $\widebar{K}$ implies that
\begin{equation}\label{eq: symmetry}
    SHI(-(S^3\backslash N(\widebar{K})),\Gamma_{-n})\cong \operatorname{Hom}(SHI(S^3\backslash N(K),\Gamma_{n}),\comp),
\end{equation}
where the minus sign denotes the opposite orientation of the manifold. This allows us to translate the results from \cite{li2019tau} into our current setup. In particular, for $n\in\intg$, we define
\begin{equation}\label{eq: D^gamma_n}
	\dg{n}{K} = \dim KHI(S^3_n(K),\widetilde{K}_n) = \dim SHI(S^3\backslash N(K),\Gamma_{-n})
\end{equation}
\blem[{\cite[Lemmas 5.2 and 5.5]{li2019tau}}]\label{lem: d^gamma is unimodular}
For any knot $K\subset S^3$, the sequence $\{\dg{n}{K}\}_{n\in\intg}$ is unimodal with unique minimum $n=2\tau_I(K)$, where $\tau_I(K)$ is the concordance tau invariant from instanton theory. That is, we have
\[
\dg{n+1}{K}=
\begin{cases}
	\dg{n}{K}-1& n < 2  \tau_I(K);\\

	\dg{n}{K}+1& n \geq 2  \tau_I(K).
\end{cases}
\]
\elem

\brem
In \cite[Lemmas 5.2 and 5.5]{li2019tau}, the authors only stated that the unique minimum of the sequence $\{\dg{n}{K}\}_{n\in\intg}$ has an index $n_0(K)\leq 2\tau_I(K)$. The inequality is indeed an equality for the following reason. The inequality holds for any knot, in particular, for the mirror knot. Since $\tau_I(\widebar{K})=-\tau_I(K)$ by \cite[Corollary 1.3]{li2019tau}, the inequality in the other direction also holds, and we conclude the equality. Also, Baldwin-Sivek \cite{baldwin2020concordance} introduced another version of tau invariant denoted by $\tau^{\sharp}$ and in \cite{li2019tau} it is shown that $\tau_I=\tau^{\sharp}$ for any knot $K\subset S^3$. 
\erem

\subsection{Dimension formula for Dehn surgeries}
In \cite{baldwin2020concordance}, Baldwin-Sivek proved a dimension formula for the framed instanton homology of Dehn surgeries along a knot $K\subset S^3$. For a knot $K\subset S^3$ and $n\in\intg$, we write
\begin{equation}\label{eq: D^sharp_n}
	\dsc{n}{K} = \dim\ish{S^3_n(K);\comp}.
\end{equation}

\blem[{\cite[Theorem 1.2]{baldwin2020concordance}}]\label{lem: dimension formula}
Suppose $K\subset S^3$ is a knot. If $K$ is V-shaped, then
\[
\dsc{n}{K}= \dsc{\nu^{\sharp}(K)}{K}+|n-\nu^{\sharp}(K)|
\]where $\nu^{\sharp}(K)$ denotes the concordance invariant introduced by Baldwin-Sivek. If instead $K$ is W-shaped, then $\nu^{\sharp}(K) = 0$ and
\[
\dsc{n}{K}=\begin{cases}
	\dim\ish{S^3_0(K),\mu;\comp} +2 & n = 0\\

	\dim\ish{S^3_0(K),\mu;\comp} + |n| &n\in \intg\backslash\{0\}.
\end{cases}\]
\elem

\brem
In \cite{baldwin2020concordance}, Baldwin and Sivek introduced the notions of V-shaped and W-shaped knots based on the sequence $\{\dim I^{\sharp}(S^3_n(K);\mathbb{C})\}_{n\in\mathbb{Z}}$, specifically using complex coefficients. We adhere to these definitions in the current paper. In particular, in Proposition \ref{prop: shape of dst, intro}, we only characterize the corresponding sequence over $\mathbb{F}_2$ as having a V-shape or W-shape, without explicitly classifying the knot itself as V-shaped or W-shaped.
\erem

\blem[{\cite[Theorem 1.2]{baldwin2020concordance}}]\label{lem: dimension formula, half integral surgery}
Suppose $K\subset S^3$ is a knot and $n\in\mathbb{Z}\backslash \{0\}$. Then
\[
\dim \ish{S^3_{\frac{2n-1}{2}}(K);\comp} = \begin{cases}
	2  \dsc{n}{K} - 1 & \nu^{\sharp}(K) < n;\\

	2  \dsc{n}{K} + 1 & \nu^{\sharp}(K) \geq n.
\end{cases}
\]If $\nu^{\sharp}(K)\neq 0$, then the result also holds for $n=0$.
\elem
Results from \cite{baldwin2020concordance}, \cite{baldwin2022concordanceII}, and \cite{li2019tau} imply the following.
\blem\label{lem: relation between nu and tau}
Suppose $K\subset S^3$ is a knot. Then either $\nu^{\sharp}(K)\in\{2\tau_I(K)+1,2\tau_I(K)-1\}$ or $\nu^{\sharp}(K)=\tau_I(K)=0$.
\elem

Another useful result as a corollary of \cite[Proposition 3.14]{LY2020} and \cite[Theorem 1.2]{LY2021} is the following.
\blem\label{lem: relation between D^gamma and D^sharp}
Suppose $K\subset S^3$ is a knot and $n\in\intg\backslash\{0\}$. Then there exists $k_n\in\mathbb{N}$ such that
\[
\dg{n}{K}=\dsc{n}{K}+2  k_n.
\]
\elem
\bdefn\label{defn: dually Floer simple}
We say that a knot $K\subset S^3$ is \emph{dually Floer simple for $n\in\intg\backslash\{0\}$} if 
\[
\dg{n}{K}=\dsc{n}{K}, \text{ or equivalently }k_n=0.
\]
If such an $n$ exists but its value is not specified, we simply say that $K$ is \emph{dually Floer simple}.
\edefn

\section{\texorpdfstring{$2$}{2}-torsion in framed instanton homology}\label{sec: Floer simple and fiberedness}
In this section, we first prove Lemmas \ref{lem: dual knot of 1-surgery} and \ref{lem: fiberness detection}, and then use them to derive Theorems \ref{thm: half surgery} and \ref{thm: non-fibered knot}, and Propositions \ref{prop: dual knot has 2-torsion} and \ref{prop: nftp is concordance inv, intro}.

\subsection{Surjectivity on the top summand}
In this subsection, we use the integral surgery formula to prove the following technical proposition about the differentials $\widetilde{d}_{1,\pm}$ for $\widetilde{K}_1\subset S_1^3(K)$.

\bprop\label{prop: surjective map}
Suppose $K \subset S^3$ is a nontrivial knot and let $\widetilde{K}_1$ be the dual knot in $S_1^3(K)$. Write \[KHI(-S^3_1(K),\widetilde{K}_1)=\bigoplus_{i\in[-g,g]}KHI(-S^3_1(K),\widetilde{K}_1),i)\]for the Alexander grading decomposition for any Seifert surface of $\widetilde{K}_1$, where the minus sign denotes the opposite orientation of the manifold, and $g=g(\widetilde{K}_1)=g(K)$. Suppose \begin{equation}\label{eq: dual d1}
    \widetilde{d}_{1,\pm}: KHI(-S^3_1(K),\widetilde{K}_1,i)\to KHI(-S^3_1(K),\widetilde{K}_1,i\pm 1)
\end{equation}are the first differentials corresponding to $\widetilde{K}_1$ (see \cite[Theorem 3.20]{LY2021large} for the definition). Then $\widetilde{d}_{1,\pm}$ is surjective onto $KHI(-S^3_1(K),\widetilde{K}_1,\pm g)$, respectively.
\eprop
\bpf
Note that the dual knot $\widetilde{K}_1$ is also null-homologous, and hence canonically framed by Seifert longitude. We know $S^3$ is obtained from $S_1^3(K)$ 
 by $(-1)$-surgery on $\widetilde{K}_1$. Hence $-S^3$ is obtained from $-S_{1}^3(K)$ by $(+1)$-surgery on $\widetilde{K}_1$. We apply two cases in the truncated version of the integral surgery formula \cite[Proposition 2.22]{LY2022integral2} (note that $(p,q,m)=(1,0,-1)$ in our case, see also \cite{LY2022integral1} for the full proof).

If $g=1$, then the second case in \cite[Proposition 2.22]{LY2022integral2} applies and we have 
\[\ish{S^3;\comp}\cong H_*(A(0)),\]
where $A(0)$ is the bent complex \[(KHI(-S_{1}^3(K),\widetilde{K}_1),(\widetilde{d}_{1,+}+\widetilde{d}_{1,-})\big|_{KHI(-S_{1}^3(K),\widetilde{K}_1,0)}).\] Note that $\widetilde{d}_{1,\pm}$ shifts the Alexander grading by $\pm 1$. If $\widetilde{d}_{1,+}$ is not surjective onto the summand $KHI(-S^3_1(K),\widetilde{K}_1,1)$, by the duality between $\widetilde{d}_{1,\pm}$ \cite[Corollary 3.35]{LY2021large}, we know $\widetilde{d}_{1,-}$ is also not surjective onto the summand $KHI(-S^3_1(K),\widetilde{K}_1,-1)$. Then $H_*(A(0))$ is at least $2$-dimensional, which contradicts the fact that $\ish{S^3;\comp}$ is $1$-dimensional \cite[Proposition 4.2]{kronheimer2011khovanov}.

If $g>1$, then the first case in \cite[Proposition 2.22]{LY2022integral2} applies, and we have \[\ish{S^3;\comp}\cong H_*(\cone\big(\bigoplus_{|s|<g}H_*(A(s))\xra{\Pi_-+\Xi_{-1}\circ \Pi_+} \bigoplus_{s\in [2-g,g-1]}H_*(B^-(s))\big)).\]Since $\Xi_{-1}:H_*(B^+(s))\to H_*(B^-(s+1))$ is an isomorphism, and $\Pi_\pm$ sends $H_*(A(s))$ to $H_*(B^\pm(s))$, respectively for all $s$, we know \[\ke(\Pi_-:H_*(A(g-1))\to H_*(B^-(g-1)))\aand \ke (\Pi_+:H_*(A(1-g))\to H_*(B^+(1-g)))\]contribute to the homology of the mapping cone. 

If $\widetilde{d}_{1,+}$ is not surjective onto the summand $KHI(-S^3_1(K),\widetilde{K}_1,g)$, the elements in \[KHI(-S^3_1(K),\widetilde{K}_1,g)\backslash \im (\widetilde{d}_{1,+})\] contributes to the first kernel. By duality, we know $\widetilde{d}_{1,-}$ is also not surjective onto the summand $KHI(-S^3_1(K),\widetilde{K}_1,-g)$, and the elements in \[KHI(-S^3_1(K),\widetilde{K}_1,-g)\backslash \im (\widetilde{d}_{1,-})\]contributes to the second kernel. Hence, the homology of the mapping cone is at least $2$-dimensional, which again induces a contradiction.
\epf

A similar proof provides a general result.
\bcor\label{cor: generalization}
Suppose $K\subset Y$ is a nontrivial null-homologous knot with the Seifert framing and let $\widetilde{K}_1$ be the dual knot in $Y_1^3(K)$. Write $g=g(\widetilde{K}_1)=g(K)$. Suppose \begin{equation*}\label{eq: dual d1 2}
    \widetilde{d}_{1,\pm}: KHI(-Y_1(K),\widetilde{K}_1,i)\to KHI(-Y_1(K),\widetilde{K}_1,i\pm 1)
\end{equation*}are the first differentials corresponding to $\widetilde{K}_1$ (see \cite[Theorem 3.20]{LY2021large} for the definition). Then \[2\dim\operatorname{coker}\widetilde{d}_{1,\pm}|_{KHI(-Y_1(K),\widetilde{K}_1,\pm (g-1)} \le  \dim I^\sharp(Y;\comp).\]
\ecor

We also have the following corollary for the manifold without orientation reversal:
\bcor\label{cor: dual injective}
Suppose $K\subset S^3$ is a nontrivial knot and let $\widetilde{K}_1$ be the dual knot in $S_1^3(K)$. Then the maps $\widetilde{d}_{1,\pm}$ on $KHI(S^3_1(K),\widetilde{K}_1)$ are injective on $KHI(S^3_1(K),\widetilde{K}_1,\mp g(K))$, respectively, where the signs of the differentials are chosen such that $\widetilde{d}_{1,\pm}$ shift the Alexander gradings by $\pm 1$.
\ecor
\bpf
The proof of Proposition \ref{prop: surjective map} uses manifolds with opposite orientations simply to be consistent with the previous work of the authors. To obtain the corollary, we can either start with the manifolds with opposite orientations such that the orientation reversal in the statement cancels, or take the dual spaces and dual maps in the original statement via the symmetry in \eqref{eq: symmetry}. Note that the symmetry only reverses the orientation of the manifold but not that of the suture, so we need to further reverse the orientation of the suture by the symmetry in the proof of \cite[Lemma 2.5]{baldwin2018khovanov}. Unpacking the definitions of maps and isomorphisms, one can show these two constructions are indeed equivalent.
\epf
Then we prove the first crucial lemma in \S \ref{subsec: Strategy of the proof} and Proposition \ref{prop: dual knot has 2-torsion} as an immediate application. For readers' convenience, we restate these results below.

\quad

\noindent
{\bf Lemma \ref{lem: dual knot of 1-surgery}.} {\it Suppose $K\subset S^3$ is a nontrivial knot, then 
\[
\dim KHI(S^3_1(K),\widetilde{K}_1) \geq \dim\ish{S^3_1(K)} + 2  \dtop{K}.
\]
}
\bpf
From \cite[Theorem 3.20]{LY2021large}, if $K$ is a rationally null-homologous knot in a $3$-manifold $Y$, then the maps $d_{1,\pm}$ on $KHI(Y,K)$ (up to mirror) are the differentials on the first pages of the two spectral sequences from $KHI(Y,K)$ to $\ish{Y;\comp}$. From Corollary \ref{cor: dual injective} and (\ref{eq: D_top}), the ranks of $d_{1,\pm}$ are at least $D_{\operatorname{top}}(K)$, so the inequality in the lemma follows from the spectral sequences.
\epf

\noindent
{\bf Proposition \ref{prop: dual knot has 2-torsion}.} {\it If $K\subset S^3$ is a nontrivial knot, then $\ish{S^3_1(K),\widetilde{K}_1;\intg}$ has $2$-torsion.}

\bpf
Lemma \ref{lem: dimension formula, half integral surgery} and Lemma \ref{lem: mixed triangle} with $\mathbb{K}=\comp$ and suitable framings implies that
\[
\begin{aligned}
	\dim \ish{S^3_1(K),\widetilde{K}_1;\comp} &\leq  \dim  \ish{S^3_{1/2}(K);\comp} + \dim \ish{S^3;\comp}\\
 &=  \dim  \ish{S^3_{1/2}(K);\comp} + 1\\
&\leq 2  \dim \ish{S^3_1(K);\comp} +2
\end{aligned}
\]
On the other hand, the equations (\ref{eq: KHI inequality}) and Lemma \ref{lem: dual knot of 1-surgery} implies that
\[
\begin{aligned}
	\dim\ish{S^3_1(K),\widetilde{K}_1;\ft}
	&\geq 2  \dim KHI(S^3_1(K),\widetilde{K}_1)\\
	&\geq 2  \dim \ish{S^3_1(K);\comp} + 4 \dtop{K}\\
	&\geq 2  \dim \ish{S^3_1(K);\comp} + 4.
\end{aligned}
\]
Hence $\ish{S^3_1(K),\widetilde{K}_1;\intg}$ has $2$-torsion.
\epf

\subsection{Dually Floer simple knots}
In this subsection, we study dually Floer simple knots $K\subset S^3$. One initial observation is the following.

\blem\label{lem: half integer surgery has 2-torsion}
	Suppose $K\subset S^3$ is a knot and $n\in\intg\backslash\{0\}$. If $\dg{n}{K}\neq \dsc{n}{K}$, or equivalently $k_n>0$ in Lemma \ref{lem: relation between D^gamma and D^sharp}, then we have
	\[
	\ti{S^3_{\frac{2n-1}{2}}(K)}>0,
	\]where the notions are from (\ref{eq: D^gamma_n}), (\ref{eq: D^sharp_n}), and (\ref{eq: number of 2-torsion}). Under such an assumption, we conclude that $\ish{S^3_{\frac{2n-1}{2}}(K);\intg}$ must have $2$-torsion.
\elem

\bpf
	For a knot $K\subset S^3$, we take the dual knot $\widetilde{K}_n\subset S^3_n(K)$. We can take the meridian of the original knot $K\subset S^3$ to be the longitude of $\widetilde{K}_n$ to induce a framing. Applying Lemma \ref{lem: mixed triangle} to $\widetilde{K}_n$ and $\lambda = \emptyset$ over $\mathbb{K}=\ft$, we obtain the following exact triangle
\[
\xymatrix{
\ish{S^3;\ft}\ar[rr]&&\ish{S^3_{\frac{2n-1}{2}}(K);\ft}\ar[dl]\\
&\ish{S^3_n(K),\widetilde{K}_n;\ft}\ar[ul]&
}
\]
From the above triangle, we have
\begin{equation*}
	\begin{aligned}
	\dim\ish{S^3_{\frac{2n-1}{2}}(K);\ft}\geq& \dim\ish{S^3_{n}(K), \widetilde{K}_n;\ft} - 1\\
	[{\rm By ~ Equation~}(\ref{eq: 2dim})]~~~=& 2\cdot \dim\ina{S^3_{n}(K), \widetilde{K}_n;\ft} - 1\\
	[{\rm By~universal~coefficient~theorem}]~~~\geq& 2\cdot \dim\ina{S^3_{n}(K), \widetilde{K}_n;\comp} - 1\\
	[{\rm By ~ Equation~}(\ref{eq: 2dim})]~~~=& 2\cdot \dim KHI(S^3_{n}(K), \widetilde{K}_n) - 1\\
	=& 2  \dg{n}{K} - 1\\
	=& 2  \dsc{n}{K} + 4  k_n - 1.
\end{aligned}
\end{equation*}
On the other hand, Lemma \ref{lem: dimension formula, half integral surgery} implies that:
\[
\dim\ish{S^3_{\frac{2n-1}{2}}(K);\comp}\leq 2  \dsc{n}{K} +\dim \ish{S^3;\comp}=2  \dsc{n}{K} +1.
\] Hence, we conclude that
\[
\ti{S^3_{\frac{2n-1}{2}}(K)}=\frac{1}{2}\left(\dim\ish{S^3_{\frac{2n-1}{2}}(K);\ft}-\dim\ish{S^3_{\frac{2n-1}{2}}(K);\comp}\right)\geq 2  k_n-1>0.
\]
Note that the last inequality follows from the hypothesis that $k_n\neq 0$.
\epf

\bcor\label{cor: 1/2 surgery has 2-torsion}
Suppose $K\subset S^3$ is a nontrivial knot. Then $\ish{S^3_{1/2}(K);\intg}$ has $2$-torsion.
\ecor
\bpf
The corollary follows directly from Lemmas \ref{lem: half integer surgery has 2-torsion} and \ref{lem: dual knot of 1-surgery} by setting $n=1$.
\epf

Lemma \ref{lem: half integer surgery has 2-torsion} motivates us to study the situation that $\dg{n}{K}=\dsc{n}{K}$. Recall that we defined the dually Floer simple knot in Definition \ref{defn: dually Floer simple}. Now we prove the second crucial lemma in \S \ref{subsec: Strategy of the proof}.

\quad

\noindent
{\bf Lemma \ref{lem: fiberness detection}.} {\it Suppose $K\subset S^3$ is a knot. Suppose $n > 0$ is an integer such that
\begin{equation*}\label{eq: floer simple knots 2}
	\dim KHI(S^3_n(K),\widetilde{K}_n) = \dim\ish{S^3_n(K)}.
\end{equation*}
Then $K$ must be fibered, $V$-shaped, and $0<\nu^{\sharp}(K)= 2  \tau^{\sharp}(K) - 1<n$, where the notions are from \cite{baldwin2020concordance}.
}
\bpf
Suppose $K\subset S^3$ is dually Floer simple such that the dual knot $\widetilde{K}_n$ in the $n$-surgery is Floer simple. Write $\tau_I=\tau_I(K)$ and $\nu^\sharp=\nu^\sharp(K)$. We discuss several cases.

{\bf Case 1}. Suppose $K$ is $W$-shaped. Then $\nu^\sharp = \tau_I = 0$. Then Lemmas \ref{lem: d^gamma is unimodular} and \ref{lem: dimension formula} imply that
\[
\dg{n}{K}=\dsc{n}{K}
\] for all $n\in\intg\backslash\{0\}$. The equality for $n=1$ contradicts Lemma \ref{lem: dual knot of 1-surgery}; see Figure \ref{fig: case 1 and 2}.

{\bf Case 2}. Suppose $K$ is $V$-shaped and $\nu^\sharp>n$. Then $2\tau_I\geq \nu^\sharp-1\geq n$. As a result, we could still conclude that 
\[
\dg{1}{K}=\dsc{1}{K}
\] 
which again contradicts Lemma \ref{lem: dual knot of 1-surgery}; see Figure \ref{fig: case 1 and 2}.

\begin{figure}[ht]
	\centering
	\begin{minipage}{0.49\linewidth}
		\begin{tikzpicture}[scale=0.75]
				\draw[very thick,->] (-2.5,0) -- (5.5,0) node[anchor=north] {$x$};
    				\draw[very thick,->] (0,0) -- (0,6.5) node[anchor=east] {$y$};
    				\draw[gray,very thin] (-2,-1) grid (5,6);
 					
    				\draw[blue, very thick] (0,3) -- (1,2);
    				\draw[blue, very thick] (0,3) -- (-1,2);
    				\draw[lygreen, very thick] (1,2) -- (5,6);
    				\draw[lygreen, very thick] (-1,2) -- (-2,3);
    				\draw[red, very thick] (0,1) -- (1,2);
    				\draw[red, very thick] (0,1) -- (-1,2);
    				
    				\filldraw[black] (1,2) circle (4pt);
    				\filldraw[black] (3,4) circle (4pt);
    				
    				\draw[black, very thick, dashed] (1,2) -- (1,0);
    				\draw[black, very thick, dashed] (3,4) -- (3,0);
    				
    				\node at (1,-0.3) {$1$};
    				\node at (0,-0.3) {$0$};
    				\node at (2,-0.3) {$\cdots$};
    				\node at (3,-0.3) {$n$};
		\end{tikzpicture}
	\end{minipage}
	\begin{minipage}{0.49\linewidth}
		\begin{tikzpicture}[scale=0.75]
				\draw[very thick,->] (-2.5,0) -- (6.5,0) node[anchor=north] {$x$};
    				\draw[very thick,->] (0,-1.5) -- (0,6.5) node[anchor=east] {$y$};
    				\draw[gray,very thin] (-2,-1) grid (5,6);
 					
    				\draw[lygreen, very thick] (0,6) -- (5,1);
    				
    				\filldraw[black] (1,5) circle (4pt);
    				\filldraw[black] (3,3) circle (4pt);
    				\filldraw[black] (5,1) circle (4pt);
    				
    				\draw[black, very thick, dashed] (1,5) -- (1,0);
    				\draw[black, very thick, dashed] (3,3) -- (3,0);
    				\draw[black, very thick, dashed] (5,1) -- (5,0);
    				
    				\node at (1,-0.3) {$1$};
    				\node at (2,-0.3) {$\cdots$};
    				\node at (3,-0.3) {$n$};
    				\node at (3.75,-0.3) {$\cdots$};
    				\node at (5,-0.3) {$\nu^{\sharp}-1$};
	\end{tikzpicture}
	\end{minipage}
	
	\caption{Left: Case 1. Right: Case 2. The blue curve represents the sequence $\{\dsc{n}{K}\}_{n\in\intg}$ and the red curve represents the sequence $\{\dg{n}{K}\}_{n\in\intg}$. If the two curves coincide, they are drawn in green.}	\label{fig: case 1 and 2}	
\end{figure}

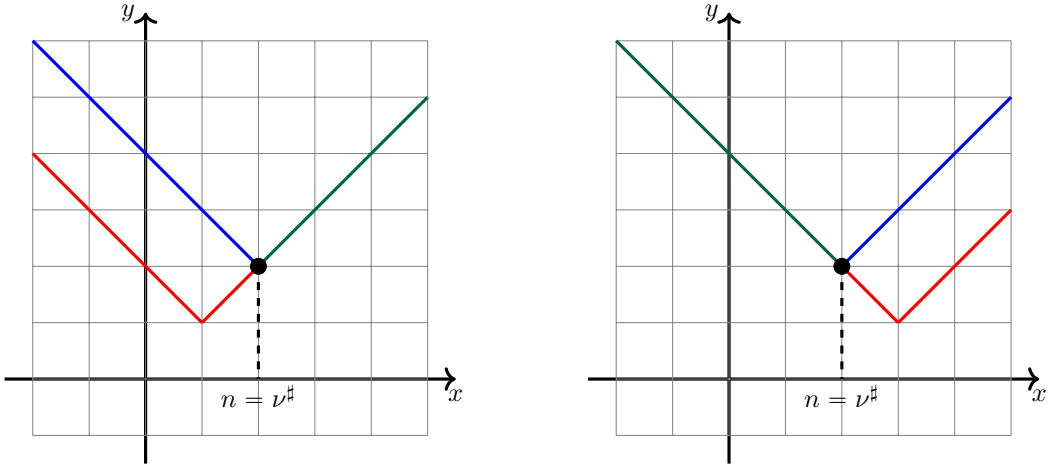
\begin{figure}[ht]
	\centering
	\begin{minipage}{0.49\linewidth}
		\begin{tikzpicture}[scale=0.75]
				\draw[very thick,->] (-2.5,0) -- (5.5,0) node[anchor=north] {$x$};
    				\draw[very thick,->] (0,-1.5) -- (0,6.5) node[anchor=east] {$y$};
    				\draw[gray,very thin] (-2,-1) grid (5,6);
 					
    				\draw[blue, very thick] (-2,6) -- (2,2);
    				\draw[lygreen, very thick] (2,2) -- (5,5);
    				\draw[red, very thick] (1,1) -- (2,2);
    				\draw[red, very thick] (-2,4) -- (1,1);
    				
    				\filldraw[black] (2,2) circle (4pt);
    				
    				\draw[black, very thick, dashed] (2,2) -- (2,0);
    				
    				\node at (2,-0.3) {$n=\nu^{\sharp}$};
		\end{tikzpicture}
	\end{minipage}
	\begin{minipage}{0.49\linewidth}
		\begin{tikzpicture}[scale=0.75]
				\draw[very thick,->] (-2.5,0) -- (5.5,0) node[anchor=north] {$x$};
    				\draw[very thick,->] (0,-1.5) -- (0,6.5) node[anchor=east] {$y$};
    				\draw[gray,very thin] (-2,-1) grid (5,6);
 					
    				\draw[lygreen, very thick] (-2,6) -- (2,2);
    				\draw[blue, very thick] (2,2) -- (5,5);
    				\draw[red, very thick] (3,1) -- (2,2);
    				\draw[red, very thick] (3,1) -- (5,3);
    				
    				\filldraw[black] (2,2) circle (4pt);
    				
    				\draw[black, very thick, dashed] (2,2) -- (2,0);
    				
    				\node at (2,-0.3) {$n=\nu^{\sharp}$};
			\end{tikzpicture}
		\end{minipage}
	
	\caption{Left: Case 3, $2\tau_I = \nu^\sharp - 1$. Right: Case 3, $2\tau_I = \nu^\sharp + 1$. The blue curve represents the sequence $\{\dsc{n}{K}\}_{n\in\intg}$ and the red curve represents the sequence $\{\dg{n}{K}\}_{n\in\intg}$. If the two curves coincide, they are drawn in green.}	\label{fig: case 3}	
\end{figure}

{\bf Case 3}. Suppose $K$ is $V$-shaped and $\nu^\sharp=n$. In this case either $2\tau_I = \nu^\sharp+1 = n+1$ which implies that
\[
\dg{n+1}{K}=\dsc{n+1}{K}-2
\]
or $2\tau_I = \nu^\sharp-1 = n-1$ which implies that
\[
\dg{n-1}{K}=\dsc{n-1}{K}-2
\]
and both situations violate Lemma \ref{lem: relation between D^gamma and D^sharp}; see Figure \ref{fig: case 3}.

{\bf Case 4}. Suppose $K$ is $V$-shaped, $\nu^\sharp<n$, and $2\tau_I = \nu^\sharp-1$. In this situation, we conclude that
\[
\dg{\nu^\sharp-1}{K}=\dsc{\nu^\sharp-1}{K}-2
\]
which is impossible by Lemma \ref{lem: relation between D^gamma and D^sharp}; see Figure \ref{fig: case 4 and case 5}.

{\bf Case 5}. Suppose $K$ is $V$-shaped, $\nu^\sharp<n$, and $2\tau_I = \nu^\sharp+1$. In this situation, we conclude that
\[
\dg{m}{K}=\dsc{m}{K}.
\]
for any $m\geq 2\tau_I$ and
\[
\dg{m}{K}=\dsc{m}{K}+2.
\]
for any $m\leq \nu^\sharp = 2\tau_I - 1$.
In particular, we must have $\nu^\sharp>0$ to make sure that
\[
\dg{1}{K}=\dsc{1}{K}+2\text{ rather than }\dg{1}{K}=\dsc{1}{K}.
\]
Then Lemma \ref{lem: dual knot of 1-surgery} implies that
\[
\dtop{K} = 1
\]
which further implies that $K$ is a fibered knot by \cite[Proposition 4.1]{kronheimer2010instanton}; see Figure \ref{fig: case 4 and case 5}.
\begin{figure}[ht]
	\centering
	\begin{minipage}{0.49\linewidth}
		\begin{tikzpicture}[scale=0.75]
				\draw[very thick,->] (-2.5,0) -- (5.5,0) node[anchor=north] {$x$};
    				\draw[very thick,->] (0,-1.5) -- (0,6.5) node[anchor=east] {$y$};
    				\draw[gray,very thin] (-2,-1) grid (5,6);
 					
    				\draw[blue, very thick] (-2,6) -- (2,2);
    				\draw[lygreen, very thick] (2,2) -- (5,5);
    				\draw[red, very thick] (1,1) -- (2,2);
    				\draw[red, very thick] (-2,4) -- (1,1);
    				
    				\filldraw[black] (4,4) circle (4pt);
    				\filldraw[black] (2,2) circle (4pt);
    				
    				\draw[black, very thick, dashed] (4,4) -- (4,0);
    				\draw[black, very thick, dashed] (2,2) -- (2,0);
    				
    				\node at (4,-0.3) {$n$};
    				\node at (3,-0.3) {$\cdots$};
    				\node at (2,-0.3) {$\nu^{\sharp}$};
		\end{tikzpicture}
	\end{minipage}
	\begin{minipage}{0.49\linewidth}
		\begin{tikzpicture}[scale=0.75]
				\draw[very thick,->] (-2.5,0) -- (5.5,0) node[anchor=north] {$x$};
    				\draw[very thick,->] (-2,-1.5) -- (-2,6.5) node[anchor=east] {$y$};
    				\draw[gray,very thin] (-2,-1) grid (5,6);
 					
    				\draw[lygreen, very thick] (1,1) -- (5,5);
    				\draw[blue, very thick] (1,1) -- (-2,4);
    				\draw[red, very thick] (2,2) -- (-2,6);
    				
    				\filldraw[black] (4,4) circle (4pt);
    				\filldraw[black] (1,1) circle (4pt);
    				\filldraw[black] (2,2) circle (4pt);
    				\filldraw[black] (-1,3) circle (4pt);
    				\filldraw[black] (-1,5) circle (4pt);
    				
    				\draw[black, very thick, dashed] (4,4) -- (4,0);
    				\draw[black, very thick, dashed] (2,2) -- (2,0);
    				\draw[black, very thick, dashed] (1,1) -- (1,0);
    				\draw[black, very thick, dashed] (-1,5) -- (-1,0);
    				
    				\node at (4,-0.3) {$n$};
    				\node at (3.25,-0.3) {$\cdots$};
    				\node at (1,-0.3) {$\nu^{\sharp}$};
    				\node at (2,-0.3) {$\nu^{\sharp}+1$};
    				\node at (0,-0.3) {$\cdots$};
    				\node at (-1,-0.3) {$1$};
			\end{tikzpicture}
		\end{minipage}
	
	\caption{Left: Case 4. Right: Case 5. The blue curve represents the sequence $\{\dsc{n}{K}\}_{n\in\intg}$ and the red curve represents the sequence $\{\dg{n}{K}\}_{n\in\intg}$. If the two curves coincide, they are drawn in green.}	\label{fig: case 4 and case 5}	
\end{figure}
In summary, under the assumption, only Case 5 will happen. Hence, we conclude the result in the lemma.
\epf

\subsection{\texorpdfstring{$2$}{2}-torsion in integral surgeries}\label{sec: 2-torsion in integral surgeries}
In this subsection, we study the framed instanton homology of Dehn surgeries over $\ft$, aiming to show the existence of $2$-torsion. For a knot $K\subset S^3$, define a new sequence similar to (\ref{eq: D^sharp_n}) for any $n\in \intg$.
\begin{equation}\label{eq: D^sharp_n, ft}
	\dst{n}{K} = \dim \ish{S^3_n(K);\ft},
\end{equation}where the superscript $0$ denotes the trivial bundle data, while we use the superscript $\mu$ to denote the nontrivial bundle data in \S \ref{sec: Shapes of dimension sequences}. Note that with this definition we have
\begin{equation}\label{eq: defn torsion}
    \ti{S^3_n(K)} = \frac{1}{2}(\dst{n}{K} - \dsc{n}{K}).
\end{equation}

A simple observation is the following
\blem\label{lem: t_2 eventually decreases}
Suppose $K\subset S^3$ is a knot. If either $K$ is V-shaped and $n\geq \nu^{\sharp}(K)$ or $K$ is W-shaped and $n\geq 1$, then
\[
\ti{S^3_{n+1}(K)} \leq \ti{S^3_{n}(K)}.
\]
\elem
\bpf
From Lemma \ref{lem: dimension formula}, we know that when $n\ge \nu^{\sharp}(K)$, we have
\[
\dsc{n+1}{K} = \dsc{n}{K} + 1.
\]
Lemma \ref{lem: surgery exact triangle, BS} with $\mathbb{K}=\ft$ implies that
\[
\dst{n+1}{K} \leq \dst{n}{K} + \dim \ish{S^3}\leq \dst{n}{K} +1.
\]Hence, we obtain the result by (\ref{eq: defn torsion}). The case of W-shaped is similar.
\epf

\bcor\label{cor: t_2 eventually constant}
For any fixed knot $K\subset S^3$, there exists $N\in\mathbb{Z}$ such that for any integer $n\geq N$, we have
\[
\ti{S^3_{n+1}(K)} = \ti{S^3_{n}(K)}.
\] 
\ecor
\bpf
Lemma \ref{lem: t_2 eventually decreases} implies that the sequence $\{ \ti{S^3_n(K)}\}_{n\in\intg_+}$ when $n$ large is non-increasing, non-negative, and always takes integral values. So it will eventually be a constant.
\epf

\bcor\label{cor: D^ft_n eventually increasing}
For any fixed knot $K\subset S^3$, there exists $N_+\in\intg$ such that for any integer $n\geq N_+$, we have
\[
\dst{n+1}{K} = \dst{n}{K} + 1.
\]
Similarly, there exists an $N_-\in\intg$ such that for any integer $n<N_-$, we have
\[
\dst{n-1}{K} = \dst{n}{K} + 1.
\]
\ecor
\bpf
The existence of $N_+$ follows immediately from Corollary \ref{cor: t_2 eventually constant} and Lemma \ref{lem: dimension formula}. The existence of $N_-$ is obtained by passing to the mirror of $K$.
\epf

\bprop\label{prop: no 2-torsion implies Floer simple}
Suppose $K\subset S^3$ is a nontrivial knot and $n\in\mathbb{N}$ such that 
\[
\ti{S^3_n(K)} = 0.
\]
Then one of the following statements holds.
\begin{enumerate}
	\item Either $n=0$ and $K$ is W-shaped; or
	\item $K$ is dually Floer simple and $n>\nu^{\sharp}(K)$.
\end{enumerate}
\eprop

\bpf
We first deal with the case that $n=0$. If $K$ is $V$-shaped, then by passing to the mirror of $K$, we can assume that $\nu^{\sharp}(K)\leq 0$. Then we conclude from Lemma \ref{lem: t_2 eventually decreases} that
\[
\ti{S^3_1(K)} = 0
\]
as well. An application of Lemma \ref{lem: surgery exact triangle, BS} over $\ft$ implies \begin{equation}\label{eq: reverse inequ}
    \begin{aligned}
	\dim\ish{S^3_{1/2}(K);\ft} &\leq \dst{0}{K} + \dst{1}{K}\\
	&=  \dsc{0}{K} + \dsc{1}{K}\\
	&= \dim\ish{S^3_{1/2}(K);\comp}
\end{aligned}
\end{equation}where the second equation follows from the assumption $\ti{S^3_0(K)}=\ti{S^3_1(K)}=0$ and the third equation follows from Lemma \ref{lem: dimension formula} to $n=0,1$, Lemma \ref{lem: dimension formula, half integral surgery} to $n=1$, and the assumption $\nu^{\sharp}(K)\leq 0$. From the universal coefficient theorem, the inequality in (\ref{eq: reverse inequ}) must be equality and we have
\[
\ti{S^3_{1/2}(K)} = 0,
\]
which violates Corollary \ref{cor: 1/2 surgery has 2-torsion}. From now on we assume that $n> 0$.

{\bf Case 1}. $\nu^{\sharp}(K) < n$. In this case, we conclude from Lemma \ref{lem: t_2 eventually decreases} that
\[
\ti{S^3_{n+1}(K)} = 0.
\]
Then a similar application of Lemmas \ref{lem: surgery exact triangle, BS}, \ref{lem: dimension formula}, and \ref{lem: dimension formula, half integral surgery} implies that
\[
\begin{aligned}
	\dim\ish{S^3_{\frac{2n+1}{2}}(K);\ft} &\leq \dst{n}{K} + \dst{n+1}{K}\\
	&=  \dsc{n}{K} + \dsc{n+1}{K}\\
	&= \dim\ish{S^3_{\frac{2n+1}{2}}(K);\comp}
\end{aligned}
\]
Hence the universal coefficient theorem implies that 
\[
\ti{S^3_{\frac{2n+1}{2}}(K)} = 0,
\]
which by Lemma \ref{lem: half integer surgery has 2-torsion} implies that $K$ must be dually Floer simple.

{\bf Case 2}. $\nu^{\sharp}(K) \geq n$. In this case one could go backward (or switch to mirror knot) and run the argument in Case 1 to conclude that
\[
\ti{S^3_{1}(K)} = \ti{S^3_{0}(K)} = 0
\]
and this further implies that
\[
\ti{S^3_{1/2}(K)} = 0,
\]
which contradicts Corollary \ref{cor: 1/2 surgery has 2-torsion}.
\epf

Now we are able to prove the main theorems in the introduction.

\quad

\noindent
{\bf  Theorem \ref{thm: half surgery}.} {\it Suppose $K\subset S^3$ is a nontrivial knot ({\it i.e.}, not the unknot). Then $\ish{S^3_{r}(K);\intg}$ for $r=1,1/2,1/4$ all have $2$-torsion.
}
\bpf
The case of $1/2$-surgery is done in Corollary \ref{cor: 1/2 surgery has 2-torsion}. The case of $1$-surgery follows from Proposition \ref{prop: no 2-torsion implies Floer simple} and Lemma \ref{lem: fiberness detection}: if $\ish{S^3_{n}(K);\intg}$ has no $2$-torsion for some positive integer $n$, then we must have
\[n>\nu^{\sharp}(K)>0,\]which implies $n\ge 2$. The case of $1/4$-surgery follows from the fact that $1$-surgery on $(2,1)$-cable of any knot $K$ produces a $3$-manifold diffeomorphic to the one obtained from the $1/4$-surgery on $K$ \cite[Corollary 7.3]{gorden1983satellite} (with $m=n=p=1$ and $q=2$).
\epf

\quad

\noindent
{\bf Theorem \ref{thm: non-fibered knot}.} {\it Suppose $K\subset S^3$ is a knot. If there exists $n\in\intg\backslash\{0\}$ such that $\ish{S^3_n(K);\intg}$ has no $2$-torsion, then $K$ must be a fibered knot.
}
\bpf
Taking the mirror knot if necessary, we can assume there exists a positive integer $n$ such that $\ti{S^3_n(K)} = 0$. Then Proposition \ref{prop: no 2-torsion implies Floer simple} implies that $K$ is dually Floer simple. Then the theorem follows from Lemma \ref{lem: fiberness detection}.
\epf

\subsection{Shapes of dimension sequences}\label{sec: Shapes of dimension sequences}Recall that in \S \ref{sec: 2-torsion in integral surgeries}, we defined the sequence
\[
\dst{n}{K} = \dim \ish{S^3_n(K);\ft}.
\]
We further define the following sequence
\begin{equation}\label{eq: dstm}
	\dstm{n}{K} = \dim \ish{S^3_n(K),\mu;\ft}
\end{equation}
to take the extra bundle data $\mu$ into consideration. In this subsection, we study those sequences. First, we review the following technical lemma. 
\blem\label{lem: vanishing result for doubling.}
Suppose $W:Y_0\to Y_1$ is a cobordism between closed $3$-manifolds. Let $\nu\subset W$ be a properly embedded $2$-manifold representing the bundle data. We write $W:Y_0\to Y_1$ and $\lambda_{i} = Y_{i}\cap \nu$ for $i=0,1$. If $\nu\cdot S$ is odd for some embedded sphere $S\subset W$ with zero self-intersection number, then the map
\[
\ish{W,\nu;\mathbb{K}}: \ish{Y_0,\lambda_0;\mathbb{K}} \to \ish{Y_1,\lambda_1;\mathbb{K}}
\]
vanishes for any commutative ring or field $\mathbb{K}$.
\elem
\bpf
It follows from the argument in the proof of \cite[Proposition 3.3]{baldwin2020concordance}: the boundary of the neighborhood of $S$ is $S^1\times S^2$ by the self-intersection number. If the bundle data is non-trivial, then there are no flat connections, and hence a neck-stretching argument shows that the map must be zero.
\epf

We have a few basic properties of the sequences.
\blem\label{lem: dst and dstm, 1}
Suppose $K\subset S^3$ is a knot. Then we have the following.
\begin{enumerate}
	\item If $n\in\intg$ is odd, then
	\[
	\dstm{n}{K} = \dst{n}{K}.
	\]
	\item If $n\in\intg$ is even, then
	\[
	 \dstm{n}{K}-\dst{n}{K}\in\{-2,0,2\}.
	\]
	\item If $n\in\intg$ is even such that $\dstm{n}{K} \neq \dst{n}{K}$, then
	\[
	\dst{n-1}{K} = \dst{n+1}{K}.
	\]
\end{enumerate}
\elem
\bpf
Part (1) follows from the facts that $H_1(S^3_n(K);\ft)$ is trivial for $n$ odd and the isomorphism class of $\ish{Y,\lambda}$ only depends on $[\lambda]\in H_1(Y;\ft)$.

For Part (2), we take $\lambda = \mu$ and $\lambda = \emptyset$ in Lemma \ref{lem: surgery exact triangle, BS} when $n$ is even. This yields two exact triangles
\begin{equation}\label{eq: dstm and dst}
	\xymatrix{
\ish{S^3_{n}(K),\mu;\ft}\ar[rr]&&\ish{S^3_{n+1}(K);\ft}\ar@<-2pt>[d]\ar@<2pt>[d]&&\ish{S^3_{n}(K);\ft}\ar[ll]\\
&&\ish{S^3;\ft}\ar[llu]\ar[rru]&&
}
\end{equation}
This implies Part (2) since $\dim\ish{S^3;\ft}=1$.

Part (3) follows from the triangles in (\ref{eq: dstm and dst}) as well. Since $\dim\ish{S^3;\ft}=\dim \ish{S^3,\mu;\ft}=1$, we have
\[
|\dst{n\pm1}{K}-\dst{n}{K}|\leq 1\text{ and }|\dst{n\pm1}{K}-\dstm{n}{K}|\leq 1.
\]So when $\dstm{n}{K} \neq \dst{n}{K}$, we must have $\dstm{n}{K} = \dst{n}{K}\pm 2$ from Part (2) and hence
\[
\dst{n-1}{K} = \dst{n+1}{K}.
\]
\epf


\blem\label{lem: dst and dstm, 3}
Suppose $n$ is odd and 
\[
\dst{n}{K} = \dstm{n+1}{K} + 1,
\]
Then we have
\[
\dst{n-1}{K} = \dst{n}{K} +1.
\]
\elem

\bpf
The proof follows the idea in the proof of \cite[Lemma 3.1]{baldwin2020concordance}. First, we describe the maps in the triangles from Lemma \ref{lem: surgery exact triangle, BS} explicitly.
Taking $\lambda = \emptyset$ in Lemma \ref{lem: surgery exact triangle, BS}, we have
\begin{equation}\label{eq: infty, n, n+1, original}
	\xymatrix{
	\ish{S^3_n(K),\mu;\ft}\ar[rr]&&\ish{S^3_{n+1}(K);\ft}\ar[dl]^{G_{n+1}}\\
	&\ish{S^3;\ft}\ar[ul]^{F_n^{\mu}}&
	}
\end{equation}
Here, the map $F^{\mu}_n$ is induced by the cobordism $X_n$ that is obtained from $[0,1]\times S^3$ by attaching an $n$-framed $2$-handle along $K$. The bundle data $\nu(F_n^{\mu})$ on $X_n$ inducing $F^{\mu}_n$ is the co-core disk of the $2$-handle. The map $G^{\mu}_{n+1}$ is induced by the cobordism $Z_{n+1}$ that is obtained from $[0,1]\times S^3_{n+1}(K)$ by attaching a $0$-framed $2$-handle along a meridian of the knot $K$ (that survives from the $(n+1)$-surgery along $K$ and hence becoming a curve in $S^3_{n+1}(K)$). The bundle data $\nu(G_{n+1})$ is trivial. Since we have not specify the homology orientation, the maps are determined by the cobordisms and the bundle data up to signs.

Taking $\lambda = \mu$ in Lemma \ref{lem: surgery exact triangle, BS}, we have
\begin{equation}\label{eq: infty, n, n+1, mu}
	\xymatrix{
	\ish{S^3_n(K);\ft}\ar[rr]&&\ish{S^3_{n+1}(K),\mu;\ft}\ar[dl]^{G_{n+1}^{\mu}}\\
	&\ish{S^3;\ft}\ar[ul]^{F_n}&
	}
\end{equation}
Note that we have
\[
\ish{S^3;\ft} \cong \ish{S^3,\mu;\ft}
\]
since $\mu$ bounds a disk inside $S^3$ and hence represent a trivial homology class. The cobordisms inducing $F_n$ and $G_{n+1}^{\mu}$ are the same as those for $F_n^{\mu}$ and $G_{n+1}$, but the bundle data becomes different. For $F_n$, the bundle data $\nu(F_n)$ is the union of $\nu(F_n^{\mu})$ with an annulus $[0,1]\times\mu$ capped off inside $S^3$ by a disk bounded by $\mu$. This actually yields a second copy of $\nu(F_n^{\mu})$ so the bundle data $\nu(F_n)$ can be chosen to be trivial. Similarly, we conclude that $\nu(G^{\mu}_{n+1})$ can be taken as a core disk of the $2$-handle forming $Z_{n+1}$.

Now let $n$ be the odd number as in the hypothesis of the lemma. We have the following diagram.
\[
\xymatrix{
\ish{S^3_{n-1}(K);\ft}\ar[rr]&&\ish{S^3_{n}(K),\mu;\ft}\cong \ish{S^3_{n}(K);\ft}\ar@<-2pt>[d]_{G^{\mu}_n}\ar[rr]&&\ish{S^3_{n+1}(K),\mu;\ft}\ar[dll]\\
&&\ish{S^3;\ft}\ar[llu]\ar@<-2pt>[u]_{F_n}&&
}
\]
Note that the union of the co-core disk in $Z_n$ and the core disk in $X_n$ form an embedded sphere $S$ with zero self-intersection number. The union of the bundle data for $G_n^\mu$ and $F_n$ is just the core disk in $Z_n$, which intersects $S$ once. We know from Lemma \ref{lem: vanishing result for doubling.} that
\[F_n\circ G^{\mu}_{n} = 0 : \ish{S^3_{n}(K),\mu;\ft}\to \ish{S^3_{n}(K);\ft}.\]
By hypothesis, we know that $F_n$ is injective, so $G_{n}^{\mu} = 0$, which implies that
\[
\dst{n-1}{K} = \dst{n}{K} + 1.
\]
\epf

\blem\label{lem: dst and dstm, 4}
Suppose $n$ is odd such that
\[
\dst{n}{K} = \dst{n+1}{K} + 1,
\]
Then we know that
\[
\dstm{n-1}{K} = \dst{n}{K} +1.
\]
\elem
\bpf
Following the proof of Lemma \ref{lem: dst and dstm, 1}, we have a diagram
\[
\xymatrix{
\ish{S^3_{n-1}(K),\mu;\ft}\ar[rr]&&\ish{S^3_{n}(K);\ft}\cong\ish{S^3_{n}(K),\mu;\ft}\ar@<-2pt>[d]_{G_n}\ar[rr]&&\ish{S^3_{n+1}(K);\ft}\ar[dll]\\
&&\ish{S^3;\ft}\ar[llu]\ar@<-2pt>[u]_{F^{\mu}_n}&&
}
\]
which implies (following the proof of Lemma \ref{lem: dst and dstm, 3}) that
\[
\dstm{n-1}{K} = \dst{n}{K} +1.
\]
\epf

\blem\label{lem: dst and dstm, 5}
Suppose $n$ is even such that
\[
\dstm{n}{K} = \dst{n+1}{K} + 1,
\]
Then we know that
\[
\dst{n-1}{K} = \dst{n}{K} +1.
\]
\elem
\bpf
In this case, we have the following diagram.
\[
\xymatrix{
&&\ish{S^3_{n}(K),\mu;\ft}\ar[rr]&&\ish{S^3_{n+1}(K);\ft}\ar[dll]\\
&&\ish{S^3;\ft}\ar[lld]\ar[u]_{F^{\mu}_n}&&\\
\ish{S^3_{n-1}(K),\mu;\ft}\ar[rr]&&\ish{S^3_{n}(K);\ft}\ar[u]_{G_n}&&
}
\]
The hypothesis then implies that $F_{n}^{\mu}$ is injective and by Lemma \ref{lem: vanishing result for doubling.} we have $G_n$ is injective, which implies that
\[
\dst{n-1}{K} = \dst{n}{K} +1.
\]
\epf

\blem\label{lem: dst and dstm, 6}
Suppose $n$ is even such that
\[
\dst{n}{K} = \dst{n+1}{K} + 1,
\]
Then we know that
\[
\dst{n-1}{K} = \dstm{n}{K} +1.
\]
\elem
\bpf
In this case, we have the following diagram.
\[
\xymatrix{
&&\ish{S^3_{n}(K);\ft}\ar[rr]&&\ish{S^3_{n+1}(K),\mu;\ft}\ar[dll]\\
&&\ish{S^3;\ft}\ar[lld]\ar[u]_{F_n}&&\\
\ish{S^3_{n-1}(K);\ft}\ar[rr]&&\ish{S^3_{n}(K),\mu;\ft}\ar[u]_{G^{\mu}_n}&&
}
\]
The hypothesis then implies that $F_{n}$ is injective and by Lemma \ref{lem: vanishing result for doubling.} we have $G_n^{\mu}$ is injective, which implies that
\[
\dst{n-1}{K} = \dstm{n}{K} +1.
\]
\epf

Now we are ready to prove Propositions \ref{prop: shape of dst, intro} and \ref{prop: nftp is concordance inv, intro}. We restate Proposition \ref{prop: shape of dst, intro} with more details. We start with the definitions of two invariants $\nu^{\sharp,\ft}_{\pm}(K)$.
\bdefn
Suppose $K\subset S^3$ is a knot. Define
\[
\nftp{K} = \min\{n~|~\forall k\ge n,~\dim\ish{S^3_{k+1}(K);\ft} = \dim\ish{S^3_{k}(K);\ft}+1\}
\]
\[
\aand \nftm{K} = \max\{n~|~\forall k\le n,~\dim\ish{S^3_{k-1}(K);\ft} = \dim\ish{S^3_{k}(K);\ft}+1\}.
\]
\edefn
\noindent
{\bf Proposition \ref{prop: nftp is concordance inv, intro}.} {\it The invariants $\nu^{\sharp,\ft}_{\pm}(K)$ are concordance invariants.}
\bpf
We first observe that from the triangles (\ref{eq: infty, n, n+1, original}) and (\ref{eq: infty, n, n+1, mu}), we know the following.
\[
\nftp{K} = \min \{n~|~\forall k\geq n,~F^{\mu}_k = 0\text{ if }k\text{ odd or }F_k = 0\text{ if }k\text{ even}.\}
\]

Now, if $K_0$ and $K_1$ are two knots that are concordant to each other, then the proof of \cite[Theorem 3.7]{baldwin2020concordance} (about the concordance invariant $N(K)$) applies here as well, and we conclude that
\[F_n(K_0)\neq 0\Leftrightarrow F_n(K_1)\neq 0\text{ and }F^{\mu}_n(K_0)\neq 0\Leftrightarrow F^{\mu}_n(K_1)\neq 0\]
As a result, we conclude that $\nftp{K}$ is a concordance invariant and $\nftm{K}$ is as well due to the symmetry between $K$ and $\widebar{K}$.
\epf

\bprop[Proposition \ref{prop: shape of dst, intro}]\label{prop: shape of dst}
For a knot $K\subset S^3$, we have the following two results regarding $\nu^{\sharp,\ft}_{\pm}$.
\begin{itemize}
		\item For any $n\ge\nftp{K}$, we have
		\begin{equation}\label{eq: dst increasing}
			\dst{n+1}{K} = \dstm{n+1}{K} = \dst{n}{K}+1 = \dstm{n}{K}+1
		\end{equation}
		
		\item For any $n\le \nftm{K}$, we have
		\begin{equation}\label{eq: dst decreasing}
			\dst{n-1}{K} = \dstm{n-1}{K} = \dst{n}{K}+1 = \dstm{n}{K}+1
		\end{equation}
\end{itemize}
Furthermore, the sequence $\{\dst{n}{K}\}_{n\in\intg}$ has one of the following three shapes.
\begin{enumerate}
	\item V-shaped (see Figure \ref{fig: V-shape}): we have that $\{\dst{n}{K}\}_{n\in\intg}$ is unimodal, {\it i.e.}, it has a unique minimum at integer $m=\nu^{\sharp,\ft}_+(K) = \nu^{\sharp,\ft}_-(K)$. Furthermore, we have
	\[\dstm{n}{K} = \dst{n}{K}\text{ for }n\neq m\aand \]\[\begin{cases}
	    \dstm{m}{K} - \dst{m}{K}\in\{0,2\} &\text{ if }m\text{ is even};\\
     \dstm{m}{K} = \dst{m}{K}&\text{ if }m\text{ is odd}.
	\end{cases}\] 
	\begin{figure}[ht]
	\centering
	\begin{minipage}{0.49\linewidth}
		\begin{tikzpicture}[scale=0.75]
				\draw[very thick,->] (-3.5,0) -- (5.5,0) node[anchor=north] {$n$};
    				\draw[gray,very thin] (-3,-1) grid (5,6);
 					
    				\draw[blue, very thick] (-3,5) -- (1,1);
    				\draw[blue, very thick] (1,1) -- (5,5);
    				
    				\filldraw[blue] (1,1) circle (4pt);
    				
    				\draw[black, very thick, dashed] (1,1) -- (1,0);
    				
    				\node at (1,-0.5) {$\nftm{K}=\nftp{K}$};
		\end{tikzpicture}
	\end{minipage}
	\begin{minipage}{0.49\linewidth}
		\begin{tikzpicture}[scale=0.75]
				\draw[very thick,->] (-3.5,0) -- (5.5,0) node[anchor=north] {$n$};
    				\draw[gray,very thin] (-3,-1) grid (5,6);
 					
    				\draw[blue, very thick] (-3,5) -- (1,1);
    				\draw[blue, very thick] (1,1) -- (5,5);
    				\draw[red, very thick, dashed] (0,2) -- (1,3);
    				\draw[red, very thick, dashed] (1,3) -- (2,2);
    				
    				\filldraw[blue] (1,1) circle (4pt);
    				\filldraw[red] (1,3) circle (4pt);
    				
    				\draw[black, very thick, dashed] (1,1) -- (1,0);
    				
    				\node at (1,-0.5) {$\nftm{K}=\nftp{K}$};
		\end{tikzpicture}
		\end{minipage}
	
	\caption{The case of a V-shaped knot. The blue curve represents the sequence $\{\dst{n}{K}\}_{n\in\intg}$ and the red (dashed) curve and dots represent the sequence $\{\dstm{n}{K}\}_{n\in\intg}$. Left: General situation. Right: Note that this could only happen when $m=\nftm{K}=\nftp{K}$ is even.}\label{fig: V-shape}
\end{figure}
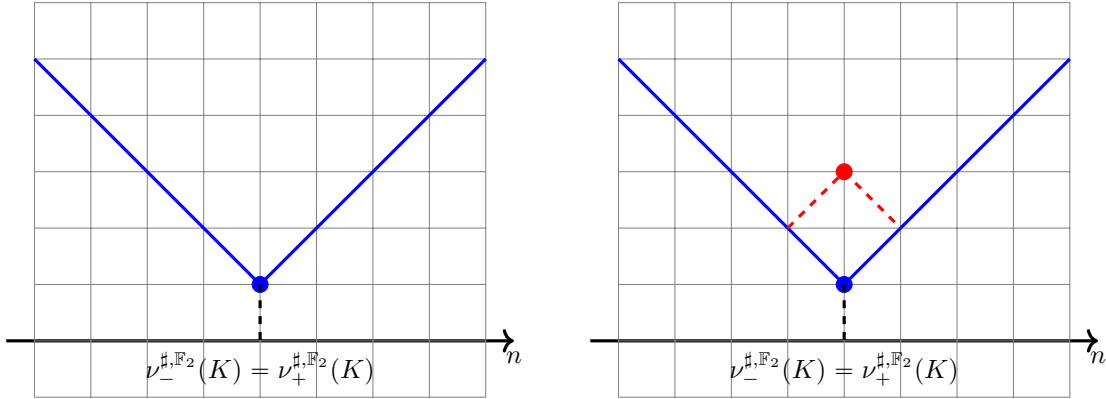

	\item W-shaped (see Figure \ref{fig: W-shape}): we have $\nftp{K} = \nftm{K}+2$ and for $m = \nftp{K} - 1 = \nftm{K} + 1$, we have
		\[\begin{cases}
		    \dst{m}{K} = \dstm{m}{K} + 2 \aand \dstm{n}{K} = \dst{n}{K} &\text{ if } m \text{ is even}\aand n\neq m;\\
                \dst{m\pm1}{K} = \dstm{m\pm 1}{K} - 2\aand \dstm{n}{K} = \dst{n}{K}&\text{ if } m \text{ is odd}\aand n\neq m\pm 1.
		\end{cases}\]
	
	\begin{figure}[ht]
	\centering
	\begin{minipage}{0.49\linewidth}
		\begin{tikzpicture}[scale=0.75]
				\draw[very thick,->] (-3.5,0) -- (5.5,0) node[anchor=north] {$n$};
    				\draw[gray,very thin] (-3,-1) grid (5,6);
 					
    				\draw[blue, very thick] (-3,5) -- (0,2);
    				\draw[blue, very thick] (2,2) -- (5,5);
    				\draw[blue, very thick] (0,2) -- (1,3);
    				\draw[blue, very thick] (1,3) -- (2,2);
    				\draw[red, very thick, dashed] (1,1) -- (0,2);
    				\draw[red, very thick, dashed] (1,1) -- (2,2);
    				
    				\filldraw[blue] (0,2) circle (4pt);
    				\filldraw[blue] (2,2) circle (4pt);
    				\filldraw[blue] (1,3) circle (4pt);
    				\filldraw[red] (1,1) circle (4pt);
    				
    				\draw[black, very thick, dashed] (0,2) -- (0,0);
    				\draw[black, very thick, dashed] (2,2) -- (2,0);
    				
    				\node at (0,-0.5) {$\nftm{K}$};
    				\node at (2,-0.5) {$\nftp{K}$};
		\end{tikzpicture}
	\end{minipage}
	\begin{minipage}{0.49\linewidth}
		\begin{tikzpicture}[scale=0.75]
				\draw[very thick,->] (-3.5,0) -- (5.5,0) node[anchor=north] {$n$};
    				\draw[gray,very thin] (-3,-1) grid (5,6);
 					
    				\draw[blue, very thick] (-3,4) -- (0,1);
    				\draw[blue, very thick] (2,1) -- (5,4);
    				\draw[blue, very thick] (0,1) -- (1,2);
    				\draw[blue, very thick] (1,2) -- (2,1);
    				
    				\draw[red, very thick, dashed] (1,2) -- (0,3);
    				\draw[red, very thick, dashed] (1,2) -- (2,3);
    				\draw[red, very thick, dashed] (-1,2) -- (0,3);
    				\draw[red, very thick, dashed] (3,2) -- (2,3);
    				
    				\filldraw[blue] (0,1) circle (4pt);
    				\filldraw[blue] (2,1) circle (4pt);
    				\filldraw[blue] (1,2) circle (4pt);
    				
    				\filldraw[red] (0,3) circle (4pt);
    				\filldraw[red] (2,3) circle (4pt);
    				
    				\draw[black, very thick, dashed] (0,1) -- (0,0);
    				\draw[black, very thick, dashed] (2,1) -- (2,0);
    				
    				\node at (0,-0.5) {$\nftm{K}$};
    				\node at (2,-0.5) {$\nftp{K}$};
		\end{tikzpicture}
	\end{minipage}
	
	\caption{The case of a W-shaped knot. The blue curve represents the sequence $\{\dst{n}{K}\}_{n\in\intg}$ and the red (dashed) curve and dots represent the sequence $\{\dstm{n}{K}\}_{n\in\intg}$. Left: $\nu^{\sharp,\ft}_{\pm}(K)$ are odd and $m$ is even. Right: $\nu^{\sharp,\ft}_{\pm}(K)$ are even and $m$ is odd.}\label{fig: W-shape}
\end{figure}
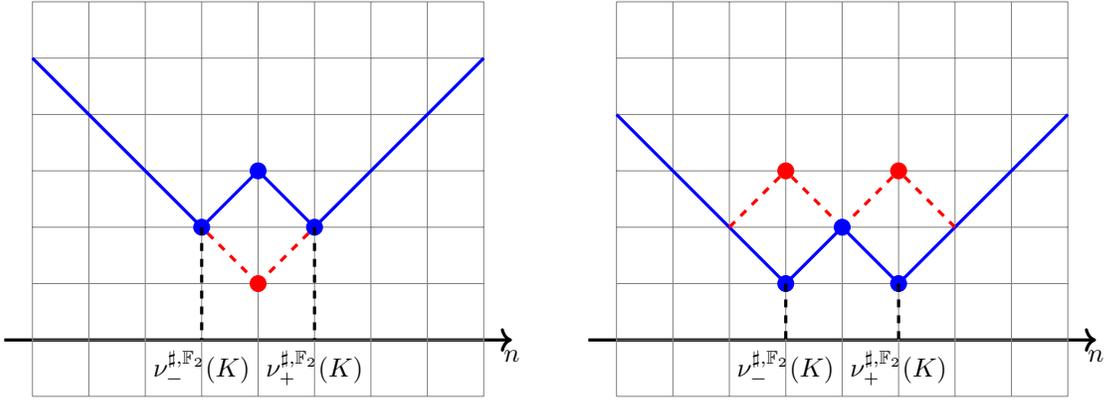

	\item Generalized W-shaped (see Figure \ref{fig: generalized W-shape}): we have $m=\nftp{K} > \nftm{K}+2$ and the following holds. For an even integer $n$ such that $n\in[\nftm{K},\nftp{K}]$, we have
		\[
		\begin{cases}
	    \dst{n}{K} = \dstm{n}{K} + 2 &\text{ if } m \text{ is odd};\\	
            \dst{n}{K} = \dstm{n}{K} - 2 &\text{ if } m \text{ is even}
		\end{cases}
		\]and for other integer $n$, we have \[\dst{n}{K} = \dstm{n}{K}.\]

\end{enumerate}
\eprop

\begin{figure}[ht]
	\centering
		\begin{minipage}{0.49\linewidth}
			\begin{tikzpicture}[scale=0.6]
				\draw[very thick,->] (-5.5,0) -- (5.5,0) node[anchor=north] {$n$};
    				\draw[gray,very thin] (-5,-1) grid (5,5);
 					
    				\draw[blue, very thick] (-5,4) -- (-3,2);
    				\draw[blue, very thick] (5,4) -- (3,2);
    				\draw[blue, very thick] (-3,2) -- (-2,3);
    				\draw[blue, very thick] (3,2) -- (2,3);
    				\draw[blue, very thick] (-2,3) -- (-1,2);
    				\draw[blue, very thick] (2,3) -- (1,2);
    				
    				\draw[red, very thick, dashed] (-3,2) -- (-2,1);
    				\draw[red, very thick, dashed] (3,2) -- (2,1);
    				\draw[red, very thick, dashed] (-2,1) -- (-1,2);
    				\draw[red, very thick, dashed] (2,1) -- (1,2);

    				\filldraw[blue] (-2,3) circle (4pt);
    				\filldraw[blue] (-1,2) circle (4pt);
    				\filldraw[blue] (1,2) circle (4pt);
    				\filldraw[blue] (2,3) circle (4pt);
    				\filldraw[blue] (-3,2) circle (4pt);
    				\filldraw[blue] (3,2) circle (4pt);
    				
    				\filldraw[red] (-2,1) circle (4pt);
    				\filldraw[red] (2,1) circle (4pt);
    				
    				\filldraw[black] (0,2) circle (1pt);
    				\filldraw[black] (-0.5,2) circle (1pt);
    				\filldraw[black] (0.5,2) circle (1pt);
    				
    				\draw[black, very thick, dashed] (-3,2) -- (-3,0);
    				\draw[black, very thick, dashed] (3,2) -- (3,0);
    				
    				\node at (-3,-0.5) {$\nftm{K}$};
    				\node at (3,-0.5) {$\nftp{K}$};
		\end{tikzpicture}
		\end{minipage}\begin{minipage}{0.49\linewidth}
			\begin{tikzpicture}[scale=0.6]
				\draw[very thick,->] (-5.5,0) -- (5.5,0) node[anchor=north] {$n$};
    				\draw[gray,very thin] (-5,-1) grid (5,5);
 					
    				\draw[blue, very thick] (-5,4) -- (-2,1);
    				\draw[blue, very thick] (5,4) -- (2,1);
    				\draw[blue, very thick] (-2,1) -- (-1,2);
    				\draw[blue, very thick] (2,1) -- (1,2);
    				
    				\draw[red, very thick, dashed] (-3,2) -- (-2,3);
    				\draw[red, very thick, dashed] (3,2) -- (2,3);
    				\draw[red, very thick, dashed] (-2,3) -- (-1,2);
    				\draw[red, very thick, dashed] (2,3) -- (1,2);

    				\filldraw[blue] (-2,1) circle (4pt);
    				\filldraw[blue] (-1,2) circle (4pt);
    				\filldraw[blue] (1,2) circle (4pt);
    				\filldraw[blue] (2,1) circle (4pt);
    				
    				\filldraw[red] (-2,3) circle (4pt);
    				\filldraw[red] (2,3) circle (4pt);
    				
    				\filldraw[black] (0,2) circle (1pt);
    				\filldraw[black] (-0.5,2) circle (1pt);
    				\filldraw[black] (0.5,2) circle (1pt);
    				
    				\draw[black, very thick, dashed] (-2,1) -- (-2,0);
    				\draw[black, very thick, dashed] (2,1) -- (2,0);
    				
    				\node at (-2,-0.5) {$\nftm{K}$};
    				\node at (2,-0.5) {$\nftp{K}$};
		\end{tikzpicture}
		\end{minipage}
	
	\caption{The case of a generalized W-shaped knot. The blue curve represents the sequence $\{\dst{n}{K}\}_{n\in\intg}$ and the red (dashed) curve and dots represent the sequence $\{\dstm{n}{K}\}_{n\in\intg}$. Left: $\nftp{K}$ is odd. Left: $\nftp{K}$ is even}\label{fig: generalized W-shape}
\end{figure}
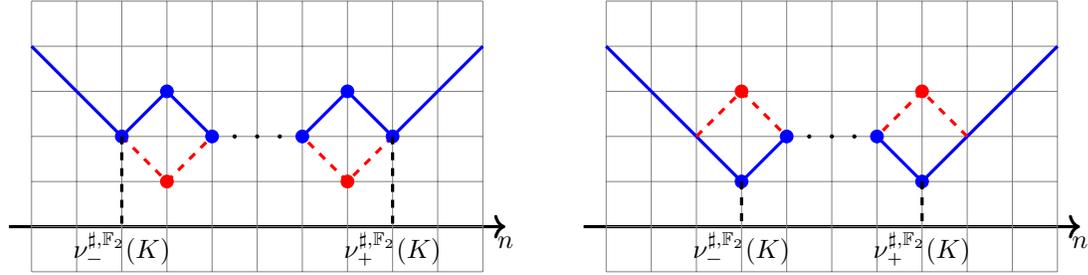
\brem
Note that Proposition \ref{prop: shape of dst} also describes the behavior of the sequence $\{\dstm{n}{K}\}_{n\in\intg}$, which is indeed also in either V-, W-, or generalized W-shapes. Following Remark \ref{rem: width}, we can also consider the widths of the generalized W-shapes for both $\{\dst{n}{K}\}_{n\in\intg}$ and $\{\dstm{n}{K}\}_{n\in\intg}$, which we write as $w_0(K)$ and $w_\mu(K)$, respectively. To define the latter width, we first define numbers similar to $\nftp{K}$ and $\nftm{K}$ by considering the dimensions with $\mu$, and then consider half of the difference. Note that Proposition \ref{prop: nftp is concordance inv, intro} and its variation involving $\mu$ show that both widths are concordance invariants. Then Proposition \ref{prop: shape of dst} implies \[|w_0(K)-w_{\mu}(K)|=1\]unless both widths are equal to zero (the first case in Figure \ref{fig: V-shape}). By a direct computation, we know $w_0(U)=1$ and $w_{\mu}(U)=0$ for the unknot $U$.
\erem
\bpf[Proof of Proposition \ref{prop: shape of dst}]
We observe that if $n>\nftp{K}$ and $\dst{n}{K} \neq \dstm{n}{K}$, then Lemma \ref{lem: dst and dstm, 1} Part (3) implies that 
\[
\dst{n-1}{K} = \dst{n+1}{K}
\]
which contradicts the definition of $\nftp{K}$. Therefore, Equation (\ref{eq: dst increasing}) follows. The proof of Equation (\ref{eq: dst decreasing}) follows similarly. 

Next, we deal with the shape. Note that, by definition, we have
\begin{equation}\label{eq: dst at nftp-1}
	\dst{\nftp{K}-1}{K} = \dst{\nftp{K}}{K} + 1 \aand \dst{\nftm{K}+1}{K} = \dst{\nftm{K}}{K} + 1.
\end{equation}
This implies that $\nftp{K}\neq \nftm{K}+1$. We discuss several cases.

{\bf Case 1.} We have $\nftp{K} = \nftm{K}$. Then $K$ is $V$-shaped by (\ref{eq: dst increasing}), (\ref{eq: dst decreasing}), and Parts (1) and (2) of Lemma \ref{lem: dst and dstm, 1}.

{\bf Case 2.} We have $\nftp{K} = \nftm{K} + 2$ and they are both odd (and $m$ is even). Equation (\ref{eq: dst at nftp-1}) implies that
\[
\dst{n}{K} = \dst{\nftp{K}}{K} + 1 = \dst{\nftm{K}}{K} + 1
\]
And hence $K$ is W-shaped. Lemma \ref{lem: dst and dstm, 6} then implies that 
\[
\dst{n}{K} = \dstm{n}{K} + 2.
\]

{\bf Case 3.} We have $\nftp{K} = \nftm{K} + 2$ and they are both even (and $m$ is odd). Equation (\ref{eq: dst at nftp-1}) directly implies that $K$ is W-shaped. If 
\[
\dstm{\nftp{K}}{K} = \dst{\nftp{K}}{K},
\]
then Equation (\ref{eq: dst at nftp-1}) together with Lemmas \ref{lem: dst and dstm, 3}, \ref{lem: dst and dstm, 4}, \ref{lem: dst and dstm, 5}, and \ref{lem: dst and dstm, 6} implies that $K$ is V-shaped and $\nftp{K} = \nftm{K}$, which is a contradiction. Hence, we have 
\[
\dstm{\nftp{K}}{K} = \dst{\nftp{K}}{K}+2.
\]
The possibility of $-2$ in the above equation can also be ruled out by Equation (\ref{eq: dst at nftp-1}). A similar argument shows that 
\[
\dstm{\nftm{K}}{K} = \dst{\nftm{K}}{K}+2.
\]
{\bf Case 4.} We have $\nftp{K} > \nftm{K} + 2$, and $\nftp{K}$ is even. Equation (\ref{eq: dst at nftp-1}) and Lemma \ref{lem: dst and dstm, 5} together imply that
\[
\dstm{\nftp{K}}{K} = \dst{\nftp{K}}{K}.
\]
Now suppose $n<\nftp{K}$ is an odd integer such that 
\[
\dstm{n+1}{K} = \dst{n+1}{K} +2,
\]
Then we know from Lemma \ref{lem: dst and dstm, 4} that
\[
\dstm{n-1}{K} = \dst{n}{K} + 1.
\]
If $\dst{n-1}{K} = \dst{n}{K} + 1$ as well, then Lemmas \ref{lem: dst and dstm, 3}, \ref{lem: dst and dstm, 4}, \ref{lem: dst and dstm, 5}, and \ref{lem: dst and dstm, 6} together imply that
\[
\nftm{K} = n+1.
\]
If $\dst{n-1}{K} = \dst{n}{K} -1$, then Lemma \ref{lem: dst and dstm, 1} Part (3) implies that
\[
\dst{n-2}{K} = \dst{n}{K}
\]
and hence the conditions that hold for $n$ also hold for $n-2$. We can perform an induction starting from $n=\nftp{K}$, which terminates when we arrive at the situation when $n+1 = \nftm{K}$. In this case, $K$ is generalized W-shaped.

{\bf Case 5.} We have $\nftp{K} > \nftm{K} + 2$, and $\nftp{K}$ is odd. If
\[
\dstm{\nftp{K}-1}{K} = \dstm{\nftp{K}}{K} + 1
\]
Then Equation (\ref{eq: dst at nftp-1}) together with Lemmas \ref{lem: dst and dstm, 3}, \ref{lem: dst and dstm, 4}, \ref{lem: dst and dstm, 5}, and \ref{lem: dst and dstm, 6}, imply that $\nftm{K}=\nftp{K}$ which contradicts the assumption.
If \[
\dstm{\nftp{K}-1}{K} = \dstm{\nftp{K}}{K} - 1
\]
then we have
\[
\dst{\nftp{K}-2}{K} = \dst{\nftp{K}}{K}
\]
and we can perform an induction similarly as in Case 4 to conclude that $K$ must be generalized W-shaped.
\epf

\section{\texorpdfstring{$2$}{2}-torsion in singular instanton knot homology}
In this section, we prove Proposition \ref{prop: genus-one knot has 2-torsion} and Theorem \ref{thm: unknotting number one knot}.
\subsection{Genus-one knots}
In this subsection, we focus on knots $K\subset S^3$ with $g(K)=1$. Since the genus is small, we can compute most of the relevant information over $\comp$ from our previous work. Recall the definition of $\dtop{K}$ in (\ref{eq: D_top}) when $g(K)=1$
\[\dtop{K}=\text{dim}_{\comp}KHI(S^3,K,1).\]

\blem\label{lem: KHI of genus one knot}
For a knot $K\subset S^3$ with $g(K)=1$, and let $\tau_I$ be the instanton tau invariant, we have
\[
\dim KHI(S^3,K,i)=
\begin{cases}
	\dtop{K} & i=\pm1;\\
	\equiv  1\pmod 2 & i=0;\\
        0 & \text{otherwise.}
\end{cases}\]
\elem

\bpf
The case $i=1$ is by definition. The cases $|i|>1$ and $i=-1$ follow from the adjunction inequality and the symmetry from \cite[\S 7]{kronheimer2010knots} respectively. The case $i=0$ follows from \cite[Theorem 1.1]{kronheimer2010instanton} and the fact that the sum of the coefficients the Alexander polynomial $\Delta_K(t)$ of a knot is always odd.
\epf

\blem\label{lem: 1-surgery}
Suppose $K\subset S^3$ is a knot with $g(K)=1$. Then we have
\[
\dim\ish{S^3_1(K);\comp} = \begin{cases}
	2\dtop{K} - 1 & \tau_{I}(K) = 1\\

	2\dtop{K} + 1 & \tau_{I}(K) < 1\\\end{cases}\]
\elem
\bpf
This is simply a restatement of \cite[Corollary 8.4]{LY2022integral2}.
\epf

We are ready to prove the main result of this subsection.

\quad

\noindent
{\bf Proposition \ref{prop: genus-one knot has 2-torsion}.} {\it If $K\subset S^3$ is a knot with genus one and the Alexander polynomial $\Delta_K(t)\neq 1$, then $\ish{S^3,K;\intg}$ has $2$-torsion.
}
\bpf
Applying Lemma \ref{lem: 1-surgery} to the mirror knot $\widebar{K}$, we conclude that
\[
\dim\ish{S^3_{-1}(K);\comp} = \begin{cases}
	2\dtop{K} + 1 & \tau_{I}(K) > -1\\

	2\dtop{K} - 1 & \tau_{I}(K) = -1\\
\end{cases}\]
Lemma \ref{lem: mixed triangle} with $(-1)$-framing then implies that
\begin{equation}\label{eq: I-sharp of K, genus one}
	\dim\ish{S^3,K;\comp}\leq
\begin{cases}
	4\dtop{K} & \tau_I(K) = \pm 1\\

	4\dtop{K} + 2 & \tau_I(K) = 0\\
\end{cases}\end{equation}

When $\tau_{I}(K) = \pm 1$, from the inequalities (\ref{eq: KHI inequality}), (\ref{eq: I-sharp of K, genus one}), and Lemma \ref{lem: KHI of genus one knot}, we have \[\begin{aligned}
    \dim \ish{S^3,K;\ft}-\dim\ish{S^3,K;\comp}\ge&2\dim KHI(S^3,K)-4\dtop{K}\\=&2  (2\dtop{K}+ 1)-4\dtop{K}\\\ge& 2
\end{aligned}\]
When $\tau_I(K) = 0$ and $\Delta_{K}(t)\neq 1$, we then claim that 
\[
\dim KHI(S^3,K,0) \geq 3
\]
 Recall from \cite[Theorem 3.20]{LY2021large} that we introduced differentials
\[
d_{1,-}: KHI(S^3,K,i) \to KHI(S^3,K,i-1)\aand d_{2,-}: KHI(S^3,K,1) \to KHI(S^3,K,-1)
\]
such that the homology \[H_*(KHI(S^3,K),d_{1,-}+d_{2,-})\cong \ish{S^3;\comp}\cong\comp,\]and $\tau_I(K)$ is the grading that supports a generator of $\ish{S^3;\comp}$ by a reinterpretation. Furthermore, we know that $d_{2,-}$ shifts the $\intg/2$ homological grading on $KHI(S^3,K)$ by $1$ from \cite[\S 3.7]{LY2021large}. Assume that
\[\Delta_K(t) = at+b+at^{-1}.\]
Since we have assumed that $\Delta_{K}(t)\neq 1$, we know $a\neq 0$. From \cite[Theorem 1.1]{kronheimer2010instanton}, we know that
\[
\chi(KHI(S^3,K,\pm 1)) = -a\neq 0.
\]

We claim that $d_{1,-}\neq 0$. If $d_{1,-}=0$, we know that
\[
\dim\ish{S^3;\comp} = \dim H_*(KHI(S^3,K),d_{2,-}) \geq \dim{\rm ker}(d_{2,-}^\p)\geq a+1,
\]
where $d_{2,-}^\p$ is the restriction of $d_{2,-}$ on the Alexander gradings $1$ and $0$, the contribution $a$ comes from $KHI(S^3,K,1)$ since $d_{2,-}$ shifts the homological grading, and the contribution $1$ comes from $KHI(S^3,K,0)$ by Lemma \ref{lem: KHI of genus one knot}. Since $a\neq 0$, we conclude that $\dim\ish{S^3;\comp}\geq 2$, which is impossible.

Now, we have $d_{1,-}\neq 0$. If $\dim KHI(S^3,K,0) = 1$, then we know that
\[
\text{either } d_{1,-}\left(KHI(S^3,K,0)\right)\neq 0,\text{ or } KHI(S^3,K,0)\subset \im (d_{1,-}),
\] 
while either case contradicts the fact that $\tau_I(K) = 0$. Hence, we conclude the claim $\dim KHI(S^3,K,0) \geq 3$ by Lemma \ref{lem: KHI of genus one knot}.

Since $\dim KHI(S^3,K,0) \geq 3$, the inequalities (\ref{eq: KHI inequality}), (\ref{eq: I-sharp of K, genus one}), and Lemma \ref{lem: KHI of genus one knot} again imply that \[\begin{aligned}
    \dim \ish{S^3,K;\ft}-\dim\ish{S^3,K;\comp}\ge&2\dim KHI(S^3,K)-(4\dtop{K}+2)\\=&2  (2\dtop{K}+ 3)-4\dtop{K}-2\\\ge& 4.
\end{aligned}\]
This concludes the proof of the proposition.
\epf



\subsection{Unknotting-number-one knots}
In this subsection, we focus on knots with unknotting number one. Suppose $L^\prime =K\cup L\subset S^3$ is a link, where $K$ is a knot and has zero linking number with any component of $L$. For any rational number $r=p/q\in\mathbb{Q}$, let $\Gamma_r$ be the suture on $S^3\backslash N(K)$ consisting of two curves of slope $-r$ with opposite orientations, with respect to the Seifert framing of $K$. Let $\Gamma_{\mu}$ again be the suture consisting of two meridians of $K$. We write
\[
\suts{r}{L} = SHI(-(S^3\backslash N(K)),-\Gamma_r,L),
\]
where $SHI$ denotes the instanton Floer homology of a balanced sutured manifold with tangle constructed by Xie-Zhang \cite{xie2019tangle}, and we regard $L\subset S^3\backslash N(K)$ as the tangle (indeed the link is the singular locus). As in \cite[\S 2]{LY2022integral1}, a Seifert surface of $K$ disjoint from $L$ induces a grading on $\suts{r}{L}$, which takes value in $\intg$ when $p$ odd and $\intg+1/2$ when $p$ even. We write the grading as
\[
\suts{r}{L} = \bigoplus_i \sutg{r}{L}{i}
\]
Since $L$ is inside the interior of $S^3\backslash N(K)$, all arguments and constructions we did without the singular locus apply verbatim. Suppose $g$ is the minimal genus of the Seifert surface of $K$ disjoint from $L$. Then we have the following lemmas.
\blem[{\cite[Theorem 2.21]{LY2020}}]\label{lem: adjunction and sut decomp}
Suppose $r=p/q\in\mathbb{Q}$, where $p$ and $q$ are co-prime and $p>0$. Then we have the following
\begin{enumerate}
	\item For $|i|>g+\frac{p-1}{2}$, we have
	\[
	\sutg{r}{L}{i} = 0.
	\]
	\item For $r\neq 0$ and $|i|=g+\frac{p-1}{2}$, we have
	\[
	\sutg{r}{L}{i} \cong SHI (-M,-\gamma,L),
	\]
	where $(M,\gamma,L)$ is obtained from $(S^3\backslash N(K),\Gamma_n,L)$ by decomposing along a minimal genus Seifert surface of $K$ disjoint from $L$.
\end{enumerate}
\elem
\blem[{\cite[Proposition 4.26]{LY2020}, or the restatement in \cite[Lemma 2.19]{LY2022integral1}}]\label{lem: str of Gamma (2n-1)/2}
Suppose $n$ is an odd integer with $n\geq 2g+1$. Then we have the following. (We have similar results when $n$ is even if we consider $i\in\intg+1/2$).
\begin{enumerate}
	\item For integer $|i|\le \frac{n-1}{2}-g$, we have
	\[
	\sutg{\frac{2n-1}{2}}{L}{i} \cong \ish{-S^3,L}.
	\]
	\item We have
	\[
	\ish{-S^3_{-n}(K),L} \cong \bigoplus_{i=\frac{1-n}{2}}^{\frac{n-1}{2}} \sutg{\frac{2n-1}{2}}{L}{i}.
	\]
\end{enumerate}
\elem

\blem[{\cite[Proposition 4.14]{LY2020}, or the restatement in \cite[Lemma 2.13]{LY2022integral1}}]\label{lem: bypass triangle 2n-1}
For any $n\in\intg$, we have two (graded) exact triangles
\[
\xymatrix{
\sutg{n}{L}{i-\frac{n-1}{2}}\ar[rr]^{\psp{n}{n-1}}&&\sutg{n-1}{L}{i+\frac{n}{2}}\ar[dl]^{\psp{n-1}{\frac{2n-1}{2}}}\\
&\sutg{\frac{2n-1}{2}}{L}{i}\ar[ul]^{\psp{\frac{2n-1}{2}}{n}}&
}
\xymatrix{
\sutg{n}{L}{i+\frac{n-1}{2}}\ar[rr]^{\psm{n}{n-1}}&&\sutg{n-1}{L}{i-\frac{n}{2}}\ar[dl]^{\psm{n-1}{\frac{2n-1}{2}}}\\
&\sutg{\frac{2n-1}{2}}{L}{i}\ar[ul]^{\psm{\frac{2n-1}{2}}{n}}&
}
\]
\elem

\blem[{\cite[Proposition 5.5]{LY2021large}}]\label{lem: composition zero}
Suppose $n\ge 2g+1$ and $|i|\le g$. Then we have
\[
\psp{n+1}{n}\circ\psm{n+2}{n+1} = 0 : \sutg{n+2}{L}{i+\frac{n+1}{2}} \to \sutg{n}{L}{i+\frac{n-1}{2}};
\]\[\psm{n+1}{n}\circ\psp{n+2}{n+1}=0: \sutg{n+2}{L}{i-\frac{n+1}{2}} \to \sutg{n}{L}{i-\frac{n-1}{2}}. \]
\elem
Now we are ready to prove Theorem \ref{thm: unknotting number one knot}.

\quad

\noindent
{\bf Theorem \ref{thm: unknotting number one knot}.} {\it If $J\subset S^3$ is a knot with unknotting number one, then
\[\dim\ish{S^3,J;\comp} \leq \dim \ina{S^3,J;\comp} + 3=\dim KHI(S^3,J)+3.\]
}

\begin{figure}[ht]
	\centering
	\begin{overpic}[width = 4.0 in]{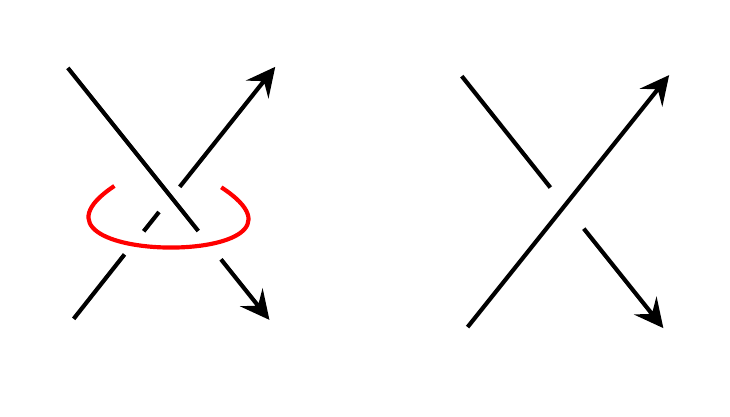}
		\put(18,5){$U\subset S^3$}
		\put(61,5){$J=U\subset S_{-1}^3(K)\cong S^3$}
		\put(35,26){$K$}
	\end{overpic}
	\vspace{-0.2in}
	\caption{The two component link $K\cup U\subset S^3$.}\label{fig: K cup L}
\end{figure}

\bpf
For simplicity, we assume all Floer homologies in the proof are over $\comp$ and omit the notation. Since $J$ has unknotting number one, we can find a two-component link $K\cup U$ such that the following hold. See Figure \ref{fig: K cup L}.
\begin{itemize}
	\item $K$ and $U$ are both unknots inside $S^3$.
	\item $K$ bounds a disk that intersects $U$ twice with opposite signs. 
	\item $U\subset S^3_{-1}(K)\cong S^3$ coincides with either $J$ or its mirror.
\end{itemize}

\begin{figure}[ht]
	\centering
	\begin{overpic}[width = 4.0 in]{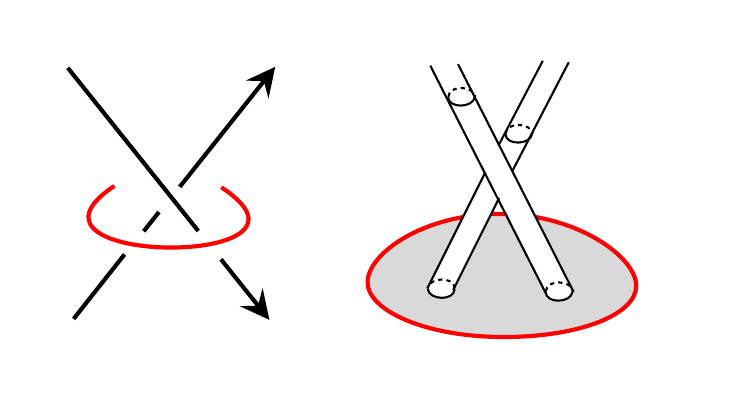}
		\put(18,5){$U\subset S^3$}
		\put(68,5){$K$}
		\put(76,35){$S$}
		\put(35,26){$K$}
	\end{overpic}
	\vspace{-0.2in}
	\caption{The Seifert surface $S$ of $K$, which is the union of the shadow twice-punctured disk and the tube along $U$.}\label{fig: Seifert surface of K}
\end{figure}
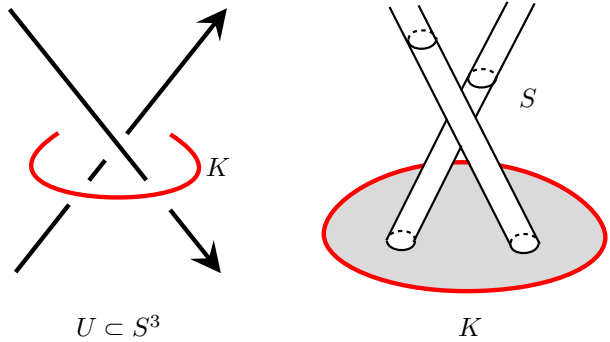

Let $L$ denote either $U$ or $U^{\natural}$ (see the definition of $U^\natural$ from (\ref{eq: defn of ina})). We first pass to the mirror of $K$, or equivalently, $K\subset -S^3$. For any $n\in \intg$, the surgery triangle from Lemma \ref{lem: surgery exact triangle, BS} along $K$ yields an exact triangle
\[
\xymatrix{
\ish{(-S^3)_n(K),L}\ar[rr]&&\ish{(-S^3)_{n+1}(K),L}\ar[dl]\\
&\ish{-S^3,L}\ar[ul]^{G_n}&
}
\]
Note that $G_n$ is the map induced by a cobordism $W_n$ obtained from $[0,1]\times (-S^3,L)$ by attaching a $4$-dimensional $2$-handle with framing $n$. Since $K$ has a genus-one Seifert surface $S$ disjoint from $L$, as shown in Figure \ref{fig: Seifert surface of K}, we can cap it off by the core of the $2$-handle. This yields a closed oriented surface $\widehat{\Sigma}_n\subset W$ with $g(\widehat{\Sigma}_n)=1$ and $\widehat{\Sigma}_n\cdot \widehat{\Sigma}_n = n$. Hence the adjunction inequality originated from \cite[Theorem 1.1]{kronheimer1995gaugeII} implies that $G_n = 0$ when $n\geq 1$. As a result, we conclude that, for any integer $n\geq 1$,
\[
\dim \ish{(-S^3)_n(K),L} = \dim \ish{(-S^3)_1(K),L} + (n-1)\cdot \dim \ish{-S^3,L}
\]
Taking $n\geq 3$, and changing back to the original knot $K\subset S^3$, Lemma \ref{lem: str of Gamma (2n-1)/2} implies that 
\begin{equation}\label{eq: 0-grading = 1-surgery}
	\sutg{\frac{2n-1}{2}}{L}{0}\cong \ish{-S^3_{-1}(K),L}.
\end{equation}
Indeed, the isomorphism (\ref{eq: 0-grading = 1-surgery}) can also be obtained from the integral surgery formula and the discussion before \cite[Proposition 2.23]{LY2022integral2}.

\begin{figure}[ht]
	\centering
	\begin{overpic}[width = 1.5 in]{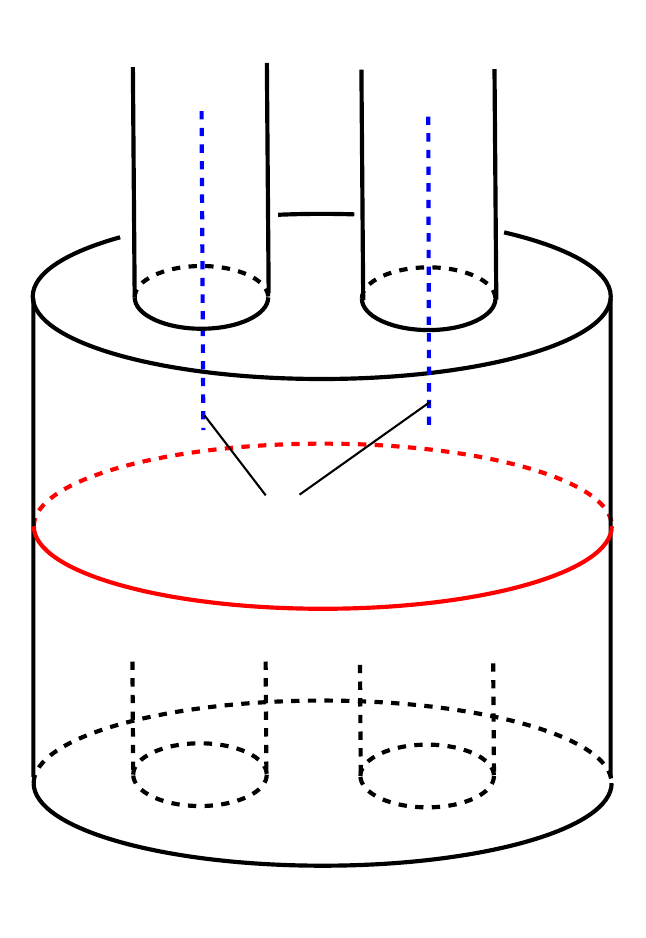}
		\put(33,2){$M$}
		\put(67,43){$\gamma$}
		\put(29,42){$L$}
	\end{overpic}
	\caption{The sutured manifold $(M,\gamma)$.}\label{fig: decompose along S}
\end{figure}

Taking $n=3$ and $4$, Lemma \ref{lem: bypass triangle 2n-1} then implies the following two exact triangles.
\[
\xymatrix{
\sutg{\frac{7}{2}}{L}{0}\ar[rr]&&\sutg{4}{L}{\frac{3}{2}}\ar[rr]\ar[dl]^{\psm{4}{3}}&&\sutg{\frac{9}{2}}{L}{-1}\ar[dl]\\
&\sutg{3}{L}{-2}\ar[ul]&&\sutg{5}{L}{-3}\ar[ul]^{\psp{5}{4}}&
}
\]
Applying the homological algebra result in \cite[Lemma A.3.10]{ozsvathbookgrid}, and using the fact that $\psm{4}{3}\circ\psp{5}{4} = 0$ in the above grading from Lemma \ref{lem: composition zero} for $g=1,n=3,i=-1$, we obtain a triangle
\[
\xymatrix{
\sutg{\frac{7}{2}}{L}{0}\ar[rr]&&\sutg{\frac{9}{2}}{L}{-1}\ar[dl]\\
&\sutg{3}{L}{-2}\oplus \sutg{5}{L}{-3}\ar[ul]& 
}
\]
From Lemma \ref{lem: adjunction and sut decomp}, Lemma \ref{lem: str of Gamma (2n-1)/2}, the equation (\ref{eq: 0-grading = 1-surgery}), we obtain the following exact triangle
\begin{equation}\label{eq: key triangle}
	\xymatrix{
\ish{-S^3_{-1}(K),L}\ar[rr]&&\ish{-S^3,L}\ar[dl]\\
&SHI(-M,-\gamma,L)\oplus SHI(-M,-\gamma,L)\ar[ul]& 
}
\end{equation}
where  $(M,\gamma)$ is obtained from $(S^3\backslash N(K),\Gamma_n)$ by decomposing along the surface $S$, and $L$ survives in $M$. See Figure \ref{fig: decompose along S}.

Taking $L=U$ in (\ref{eq: key triangle}), recalling that $J = U\subset S^3_{-1}(K)\cong S^3$, we have 
\begin{equation}\label{eq: I^sharp J}
	\begin{aligned}
	\dim \ish{-S^3,J} &\leq 2 \dim SHI(-M,-\gamma,U) + \dim \ish{-S^3, U}\\
	&=2 \dim SHI(-M,-\gamma,U)+2,
\end{aligned}
\end{equation}where the last equation follows from \cite[Lemma 8.3]{kronheimer2011khovanov}.

Taking $L=U^{\natural}$ in (\ref{eq: key triangle}), we have 
\begin{equation*}
	\begin{aligned}
	\dim \ish{-S^3,J^{\natural}} &\geq 2 \dim SHI(-M,-\gamma,U^{\natural}) - \dim \ish{-S^3, U^{\natural}}\\
	&=2 \dim SHI(-M,-\gamma,U^{\natural}) - 2,
\end{aligned}
\end{equation*}
where the last equation follows from the fact in \cite[Lemma 5.3]{xie2021earring} that for any knot $J\subset S^3$, we have
\[
\dim\ish{-S^3,J^{\natural}}=2\dim I(-S^3,J^{\natural})=2\dim \ina{-S^3,J}.
\]
Hence, we conclude
\begin{equation}\label{eq: I^sharp J^natural}
	\dim \ina{-S^3,J}\geq \dim SHI(-M,-\gamma,U^{\natural}) -1.
\end{equation}
Comparing (\ref{eq: I^sharp J}) and (\ref{eq: I^sharp J^natural}), we can conclude the theorem as long as
\begin{equation}\label{eq: SHI J vs SHI J^natural}
	\dim SHI(-M,-\gamma,U^{\natural}) =2 \dim SHI(-M,-\gamma,U).
\end{equation}

\begin{figure}[ht]
	\centering
	\begin{overpic}[width = 4.0 in]{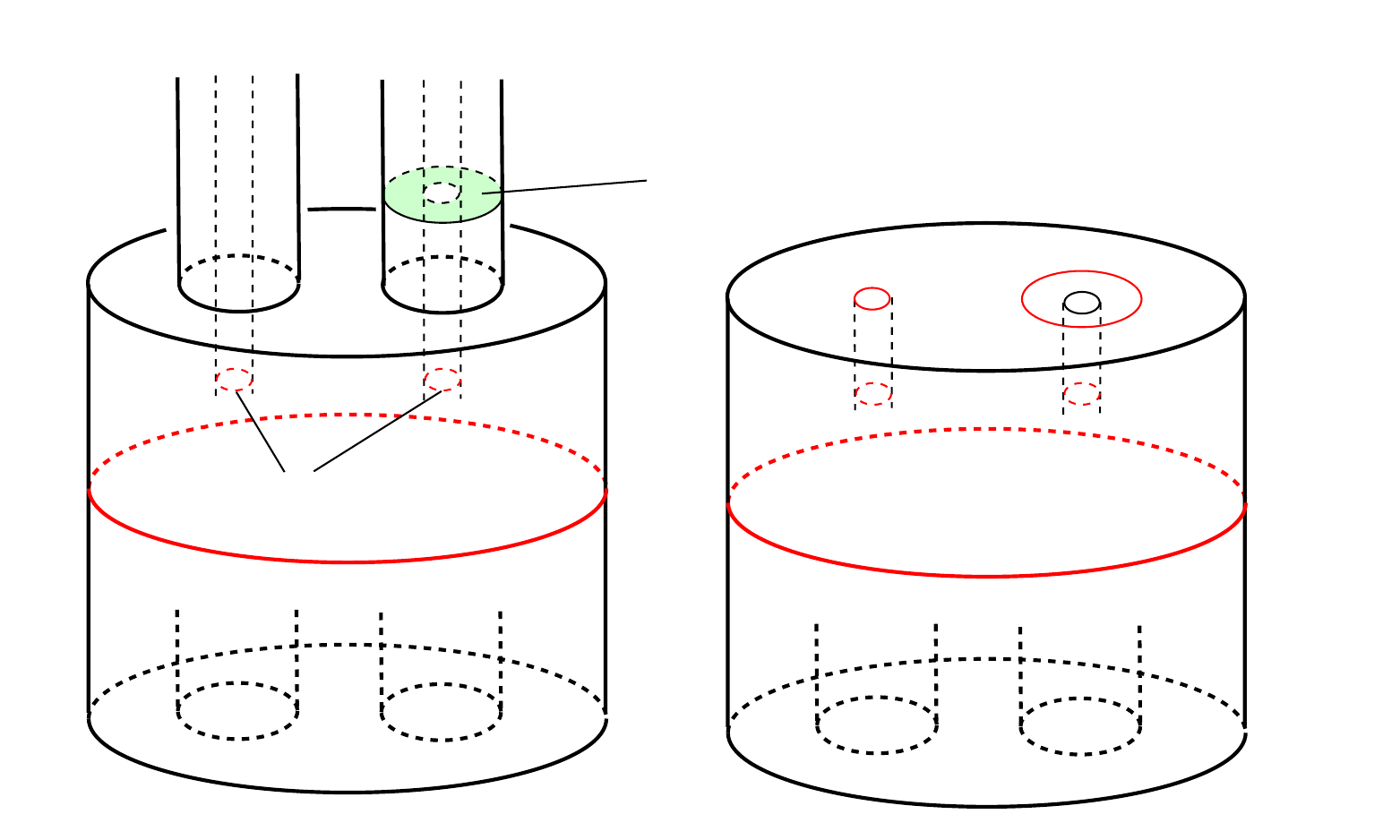}
		\put(22,-1){$M_1$}
		\put(68,-1){$M_2$}
		\put(20,23){$\mu_{U}$}
		\put(45,25){$\gamma$}
		\put(47,46){$A$}
	\end{overpic}
	\caption{Left: The sutured manifold $(M_1,\gamma_1)$. Right: The sutured manifold $(M_2,\gamma_2)$.}\label{fig: decompose along A}
\end{figure}

We first deal with the left-hand side of (\ref{eq: SHI J vs SHI J^natural}). Take $M_1 = M\backslash N(U)$ and $\gamma_1 = \gamma\cup \mu_U$, where $\mu_U$ is the union of two meridians; see Figure \ref{fig: decompose along A}. The proof of \cite[Proposition 1.4]{kronheimer2011khovanov} applies here, and we conclude that
\[
SHI(-M,-\gamma,U^{\natural})\cong SHI(-M_1,-\gamma_1).
\]
Furthermore, we observe that $(M_1,\gamma_1)$ admits a product annulus $A$. Let $(M_2,\gamma_2)$ be the result of the sutured manifold decomposition of $(M_1,\gamma_1)$ along $A$, we conclude from the instanton version of \cite[Proposition 6.7]{kronheimer2010knots} that
\begin{equation}\label{eq: U^natural = gamma_2}
	SHI(-M,-\gamma,U^{\natural})\cong SHI(-M_1,-\gamma_1)\cong SHI(-M_2,-\gamma_2).
\end{equation}
See Figure \ref{fig: decompose along A}.

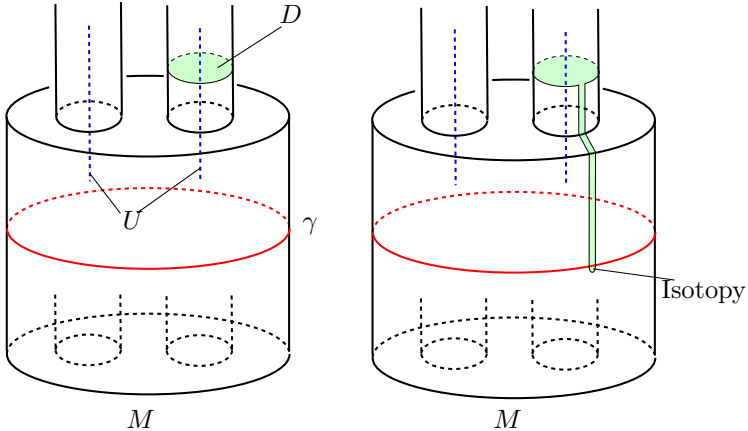
\begin{figure}[ht]
	\centering
	\begin{overpic}[width = 4.0 in]{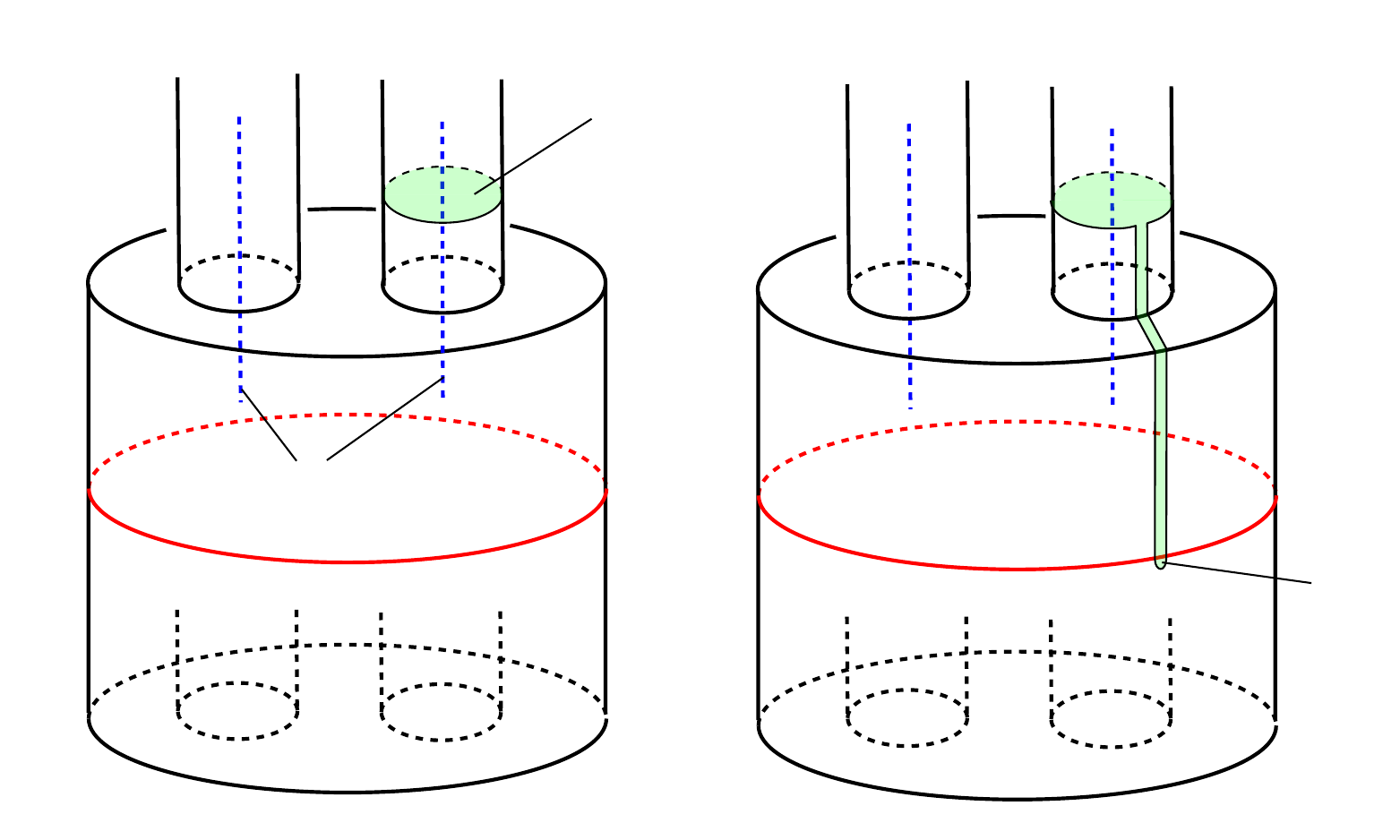}
		\put(22,-1){$M$}
		\put(70,-1){$M$}
		\put(21.5,25){${U}$}
		\put(45,25){$\gamma$}
		\put(42,52){$D$}
		\put(92,16){Isotopy}
	\end{overpic}
	\caption{Left: The disk $D\subset M$. Right: The isotopy of $\partial D$.}\label{fig: disk D}
\end{figure}

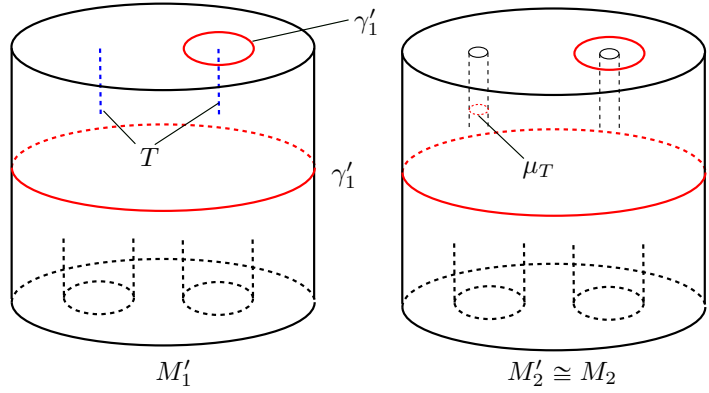
\begin{figure}[ht]
	\centering
	\begin{overpic}[width = 4.0 in]{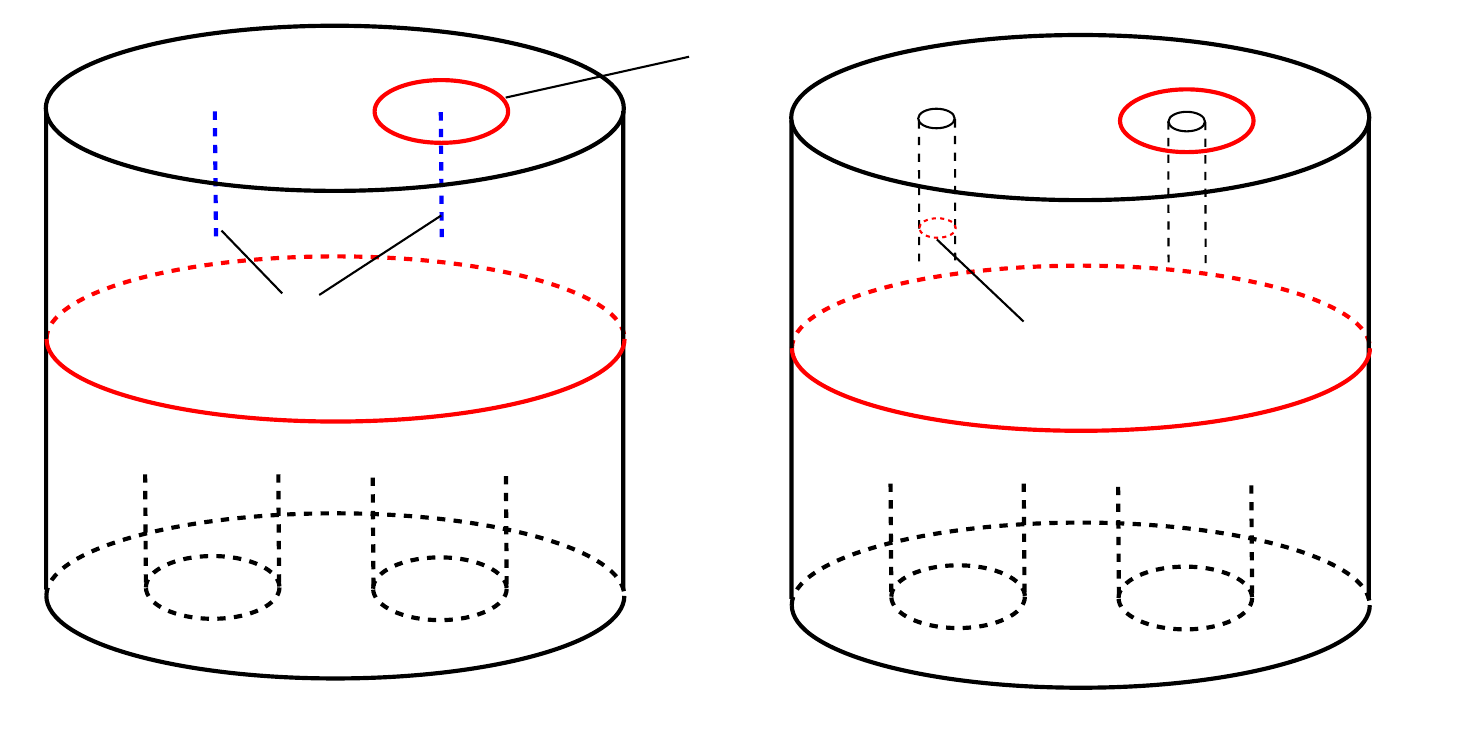}
		\put(22,-1){$M_1^{\prime}$}
		\put(68,-1){$M_2^{\prime}\cong M_2$}
		\put(20,27){${T}$}
		\put(45,25){$\gamma_1^{\prime}$}
		\put(48,45){$\gamma_1^{\prime}$}
		\put(70,26){$\mu_T$}
	\end{overpic}
	\caption{Left: The sutured manifold $(M_1^{\prime},\gamma_1^{\prime})$. Right: The sutured manifold $(M_2^{\prime}\cong M_2,\gamma_2^{\prime})$.}\label{fig: decompose along D}
\end{figure}

Next, we deal with $(M,\gamma,U)$. Since $U$ no longer has an earring, we do not perform an excision as in the proof of \cite[Proposition 1.4]{kronheimer2011khovanov}. The annulus $A\subset M_1 = M\backslash N(U)$ now becomes a disk $D$ intersecting $U$ uniquely at one point. Isotope $\partial D$ to make it intersect $\ga$ at two points, and still write the perturbed disk as $D$; see Figure \ref{fig: disk D}. Choose the orientation of $D$ such that the sutured manifold decomposition
\[
(-M,-\ga, U)\stackrel{D}{\leadsto}(-M_{1,+}^\prime,-\ga_{1,+}^\prime,T)
\]
is taut; see Figure \ref{fig: decompose along D}.

Note that we have another sutured manifold decomposition
\[
(-M,-\ga, U)\stackrel{-D}{\leadsto}(-M_{1,-}^{\prime},-\ga_{1,-}^{\prime},T),
\]
where $(-M_{1,-}^{\prime},-\ga_{1,-}^{\prime},T)$ is not taut since the isotopy creates a component of $-\gamma_{1,-}$ bounding a disk and hence $R(-\gamma_{1,-})$ is compressible. Note that, here $T\subset S^3$ is a tangle, and is vertical in the sense of \cite[Definition 1.1]{xie2019tangle}. The construction of \cite[\S 3]{li2019direct} applies to $D\subset (-M,-\gamma)$, and we can find a closure $Y$ for $(-M,-\gamma)$ such that $D$ extends to a torus $\Sigma_1\subset Y$ that intersects the singular locus $U$ transversely at one point. As a result, \cite[Proposition 6.1]{xie2019tangle} applies and hence this torus $\Sigma_1$ leads to an operator $\mu^{\operatorname{orb}}(\Sigma_1)$ on $SHI(-M,-\gamma, U)$ that has only two possible eigenvalues $1$ and $-1$. Hence, we have a generalized eigenspace decomposition
\[
SHI(-M,-\gamma, U) = SHI(-M,-\gamma, U,1)\oplus SHI(-M,-\gamma, U,-1).
\] 
Proof of \cite[Proposition 7.11]{kronheimer2010knots} then implies that
\[
SHI(-M,-\gamma, U,\pm 1)\cong SHI (-M^{\prime}_{1,\pm},-\gamma_{1,\pm},T).
\]
Since $R(-\gamma_{1,-})$ is compressible, we conclude that $SHI(-M_{1,-}^{\prime},-\gamma_{1,-}, T) = 0$ and hence
\[
SHI(-M,-\gamma, U) \cong SHI(-M_{1,+}^{\prime},-\gamma_{1,+}, T).
\]
Applying \cite[Lemma 7.10]{xie2019tangle}, we conclude that
\[
SHI(-M,-\gamma, U) \cong SHI(-M_{1,+}^{\prime},-\gamma_{1,+}, T) \cong SHI (-M_2^{\prime},-\gamma_2^{\prime}),
\]
where $M_2^{\prime} = M_{1,+}^{\prime} \backslash N(T)$ and $\gamma_2=\gamma_{1,+}\cup\mu_T$ for a meridian $\mu_T$ of $T$; see Figure \ref{fig: decompose along D}. It is straightforward to check that $M_2^{\prime}\cong M_2$ and $\gamma_2$ is the union of $\gamma_2^{\prime}$ with two more parallel copies of $\mu_T$. Then (\ref{eq: U^natural = gamma_2}) and the proof of \cite[Theorem 3.1]{kronheimer2010instanton} implies that
\[
SHI(-M,-\gamma, U^{\natural}) \cong SHI (-M_2,-\gamma_2)\cong SHI (-M_2,-\gamma_2^{\prime})\otimes \comp^2\cong SHI(-M,-\gamma, U)\otimes \comp^2.
\] 
This completes the proof of the equation in (\ref{eq: SHI J vs SHI J^natural}) and hence the inequalities (\ref{eq: I^sharp J}) and (\ref{eq: I^sharp J^natural}) conclude the proof of the theorem.
\epf





\section{Results in Heegaard Floer theory}\label{sec: HF theory}
In this section, we prove Theorem \ref{thm: Heegaard dual simple} using Heegaard Floer theory. Indeed, we obtain much stronger results, which may be of independent interest.

\subsection{Setups of immersed curve invariants}\label{subsec: curve inv}
Proposition \ref{prop: surjective map} and Corollary \ref{cor: dual injective} are motivated by observation of the immersed multi-curve invariant for $3$-manifold with a torus boundary constructed by Hanselman-Rasmussen-Watson \cite{HRW2024curve1,HRW2022curve2}, usually called \emph{curve invariant} for short. The curve invariant is based on bordered Heegaard Floer homology constructed by Lipshitz-Ozsv\'{a}th-Thurston \cite{LOT2018bordered}, whose analogue in instanton theory is still unknown. Using the curve invariant, we can prove much stronger results than Proposition \ref{prop: surjective map}. We expect analogous results to hold in instanton theory, but proofs based on current techniques such as integral surgery formulae are still beyond the authors' knowledge.

In this subsection, we first review basic facts about curve invariants. We prove results about the behavior of curves and corollaries in Heegaard Floer theory in the next subsection. We always assume that the base field is $\ft$ since the bordered theory and hence the curve invariant only work over $\ft$. We will omit the $3$-manifold and the coefficient in the Floer homology for simplicity.

Suppose $M$ is a compact, connected, oriented $3$-manifold with torus boundary. We write \begin{equation}\label{eq: multicurve inv}
    \widehat{HF}(M)\aand \bs{\ga}(M)=\{\ga_0,\dots,\ga_n\}
\end{equation}for the multi-curve invariant in \cite{HRW2024curve1,HRW2022curve2} and the underlying set of immersed curves in $\partial M\backslash z$, where $z$ is a basepoint. Note that $\widehat{HF}(M)$ consist of compact immersed curves in $\bs{\ga}(M)$ with local systems, where a curve with a local system of dimension $k$ can be roughly regarded as consisting of $k$ copies of the curve with twists between copies. All known examples from bordered Floer homology have no local systems, but a general classification theorem about an extendable type D structure over torus algebra should include ones with nontrivial local systems.

Suppose  $M_0$ and $M_1$ are two such manifolds and $h:\partial M_1\to \partial M_0$ is an orientation reversing homeomorphism preserving the basepoints $z_i$ for $i=0,1$. Then the gluing theorem states the hat version of Heegaard Floer homology (over $\ft$) of the glued closed $3$-manifold $\widehat{HF}(M_0\cup_h M_1)$ is isomorphic to the (immersed) Lagrangian intersection Floer homology of $\widehat{HF}(M_0)$ and $h(\widehat{HF}(M_1))$ in $\partial M_0\backslash z_0$. Since we work on the punctured surface, the dimension of the Lagrangian intersection Floer homology is just the geometric intersection number of the curves, {\it i.e.}, the intersection number when all bigons that do not cover the basepoint are canceled by regular homotopy of the curves. Note that a local system of dimension $k$ contributes to $k$ times the intersection number since we can always move the twists disjoint from the intersection point. In \cite[\S 7.1]{HRW2024curve1}, a standard way to cancel all bigons was introduced, and the resulting curve is called \emph{peg-board diagram}. Roughly speaking, the curve is pulled tight as a geodesic under some fixed metric on the punctured torus. 

Other than $\widehat{HF}$ of a closed $3$-manifold from gluing, the curve invariant can also be used to compute the hat and minus versions of the knot Floer homology $\widehat{HFK}$ and $HFK^-$ of all dual knots in the Dehn fillings (note that the information is not enough to compute the full knot Floer filtered complex $CFK^-$; see \cite{Hanselman2023surgery,Hanselman2023immersed} for more construction). More precisely, following \cite[\S 4.3]{HRW2022curve2}, let $\mu$ be a simple closed curve on $\partial M$ and let $\widetilde{K}_\mu$ be the dual knot in the Dehn filling along $\mu$. We consider a neighborhood disk $D$ disjoint from the curves in $\bs{\ga}(M)$, and replace the basepoint $z$ with two basepoints $z$ and $w$. We suppose $\al$ intersects $D$ and separates the two basepoints. We always put $z$ on the left of the curve as in \cite[Figure 41]{HRW2022curve2}. 

We define $C^-(\widehat{HF}(M),\mu)$ to be the chain complex generated by intersection points of $\widehat{HF}(M)$ and $\mu$ over $\ft[U]$, with the differential\begin{equation}\label{eq: diff HFK^-}
    \partial x=\sum_{i=0}^\infty\sum_yU^iN_i^{\omega}(x,y)\cdot y,
\end{equation}where $N_i^{w}(x,y)$ is the number of bigons from $x$ to $y$, counted modulo $2$, covering the $w$ basepoint $i$ times. In general, we do not expect $\partial^2=0$; again see \cite{Hanselman2023surgery,Hanselman2023immersed} for more discussion. The $z$ basepoint induces an Alexander filtration $A$ on both complexes by $A(x)-A(y)=n_z(B)-n_w(B)$ and $A(U\cdot x)=A(x)-1$, where $n_z(B)$ and $n_w(B)$ are multiplicities of basepoints in the bigon $B$. Let $gC^-(\widehat{HF}(M),\mu)$ be the associated graded chain complex, where we do not count bigons covering $z$ and do obtain a chain complex. Let $g\widehat{C}(\widehat{HF}(M),\mu)$ be obtained from $gC^-(\widehat{HF}(M),\mu)$ by setting $U=0$. Then we have filtered chain homotopy equivalences\begin{equation}\label{eq: chain equ}
    gC^-(\widehat{HF}(M),\mu)\simeq gCFK^-(\widetilde{K}_\al)\aand g\widehat{C}(\widehat{HF}(M),\mu)\simeq g\widehat{CFK}(\widetilde{K}_\mu),
\end{equation}where the homologies are $HFK^-(\widetilde{K}_\mu)$ and $\widehat{HFK}(\widetilde{K}_\mu)$, respectively, both inherited with the Alexander grading.

To draw the curves more clearly, we can lift them to some covers of $\partial M\backslash z$ and count intersection points and bigons with different underlying images. There are usually two choices of lifts. The first is the universal abelian cover $\real^2/\intg^2$ with respect to some boundary framing. The second is more standard, which corresponds to the kernel of the composite homomorphism \[\pi_1(\partial M,z)\xra{\operatorname{abelianization}} H_1(\partial M;\intg)\xra{i_*} H_1(M;\intg),\]and denoted by $\widebar{T}_{M,\mathfrak{s}}$ for each given spin$^c$ structure $\mathfrak{s}$ on $M$. The group of the deck transformations is \[H_M=\operatorname{Im}(i_*:H_1(\partial M;\intg)\to H_1(M;\intg)).\]

The second lift is useful to describe the Alexander grading of knot Floer homology. Following \cite[\S 4.4]{HRW2022curve2}, let $\lambda$ be the homological longitude which generates the kernel of $i_*$, and fix a class $[\Sigma]\in H_2(M,\partial M;\intg)$ with $\partial [\Sigma]=\lambda$. Let $\mu$ be a simple closed curve on $\partial M$ as above. We assume that $\mu\neq \lambda$ such that the dual knot $\widetilde{K}_\mu$ is rationally null-homologous. Then the set $\operatorname{Spin}^c(M,\ga_\mu)$ of spin$^c$ structures on the balanced sutured manifold $(M,\ga_\mu=\mu\cup -\mu)$ is an affine space ({\it i.e.}, torsor) over $H^2(M,\partial M;\intg)\cong H_1(M;\intg)$. Let $\operatorname{Spin}^c(M)$ be the set of spin$^c$ structures on $M$, which is an affine space over $H^2(M;\intg)\cong H_1(M,\partial M;\intg)$. Note that the fiber of the restriction map $\operatorname{Spin}^c(M,\ga_\mu)\to \operatorname{Spin}^c(M)$ is a $H_M$-affine space. Since $\partial M$ is a torus, each $\widebar{\mathfrak{s}}\in \operatorname{Spin}^c(M,\ga_\mu)$ has a well-defined first Chern class (or equivalently, the relative Euler class) in $H^2(M,\partial M;\intg)$. The half of the pairing with $[\Sigma]$ provides a function \[A: \operatorname{Spin}^c(M,\ga_\mu)\to \frac{1}{2}\intg,\]which provides an Alexander grading on $\widehat{HFK}(\widetilde{K}_\mu)\cong SFH(M,\ga_\mu)$. We write $\widehat{HFK}(\widetilde{K}_\mu,\widebar{\mathfrak{s}})$ for the summand corresponding to $\widebar{\mathfrak{s}}\in \operatorname{Spin}^c(M,\ga_\mu)$.

On the other hand, we choose $\mathfrak{s}$ to be the image of $\widebar{\mathfrak{s}}$ in $\operatorname{Spin}^c(M)$. The cover $\widebar{T}_{M,\mathfrak{s}}$ can be identified with $H_1(\partial M;\real)/\langle \lambda\rangle$. There is a natural height function \[h: \widebar{T}_{M,\mathfrak{s}}\to \real\]given by $h(v)=v\cdot \lambda$. In the computation of $g\widehat{C}(\widehat{HF}(M),\mu)$, we can replace the disk $D$ and the two basepoints $z,w$ by a single basepoint $z$ and the curve $\mu$ by a noncompact arc $L_\mu$ (as a Lagrangian submanifold) connecting the basepoint to itself with the same slope as $\mu$. Since the group of deck transformation is $H_M$, the set of lifts of $L_\mu$ is also an affine space over $H_M$, which can be identified with $\operatorname{Spin}^c(M,\ga_\mu)$ up to an overall shift. We write $L_{\mu,\widebar{\mathfrak{s}}}$ for the lift of $L_\mu$ whose midpoint has height equals to $A(\widebar{\mathfrak{s}})$. Then we have a refinement of (\ref{eq: chain equ})\begin{equation}\label{eq: alexander}
    g\widehat{HFK}(\widetilde{K}_\mu,\widebar{\mathfrak{s}})\simeq g\widehat{C}(\widehat{HF}(M),L_{\mu,\widebar{\mathfrak{s}}}),
\end{equation}where we only count intersection points with different underlying images. Though not mentioned in the reference, it is clear that the Alexander grading summand of $gCFK^-(\widetilde{K}_\mu,\widebar{\mathfrak{s}})$ corresponds to $gC^-(\widehat{HF}(M),L_{\mu,\widebar{\mathfrak{s}}})$. Note that the differentials in $gC^-$ will shift the Alexander grading, and we should use the lifts of the previous model with $z$ and $w$ to count the differentials.

Note that the filtered chain homotopy equivalence in (\ref{eq: chain equ}) implies an isomorphism between induced spectral sequences, in particular a relation between the differentials on the first pages. 

On the side of $gCFK^-$, the first differential was studied by Sarkar \cite[\S 4]{sarkar15moving} and denoted by $\Psi$. There is also another map $\Phi$ which can be regarded as the first differential on the spectral sequence corresponding to the knot with opposite orientation. Roughly speaking, the maps $\Psi$ and $\Phi$ are obtained by counting holomorphic disks passing through the basepoints $w$ and $z$ once, respectively. In Zemke's reconstruction \cite{zemke17moving} (see also \cite[\S 4.2]{Zemke2019link}), these two maps are related to some dividing sets on $K\times I\subset Y\times I$ for the given knot $K$ in the $3$-manifold $Y$, which are exactly the dividing sets of the contact structures we use to define $d_{1,\pm}$. Hence $d_{1,\pm}$ are indeed analogs of $\Psi$ and $\Phi$ (up to mirror manifolds and hence taking dual spaces and dual maps).

On the side of $gC^-$, we can first use the peg-board diagram of the curve invariant such that there is no differential on $g\widehat{C}$, then bigons with $n_z=0,n_w=1$ correspond to the differential on the first page, which coincides with $\Psi$ by the equivalence in (\ref{eq: chain equ}). Considering the Alexander grading in (\ref{eq: alexander}), we can use any two arcs $L_{\mu,\widebar{\mathfrak{s}}}$ and $L_{\mu,\widebar{\mathfrak{s}}^\p}$ that share an endpoint in the cover $T_{M,\mathfrak{s}}$ to compute $\Psi$ on the corresponding Alexander grading summand, where the basepoints $z$ and $w$ are regarded near the common endpoint and on different sides of the arcs.

For the other map $\Phi$, since we reverse the orientation of the knot rather than the orientation of $M$ (or the ambient $3$-manifold $Y$), the symmetry in \cite[\S 3.1]{HRW2022curve2} does not apply. The only difference is to replace $\lambda$ by $-\lambda$ and hence the height functions have opposite signs. From the elliptic involution in \cite[\S 3.2]{HRW2022curve2}, the map $\Phi$ can be obtained by counting bigons with $n_w=0,n_z=1$.

Finally, let us focus on the special case where $M$ is the complement of a knot $K$ in $S^3$. Then $\lambda$ is just the Seifert longitude and $[\Sigma]$ is the class of any Seifert surface which generates $H_2(M,\partial M;\intg)\cong \intg$. We have $H_1(M;\intg)\cong \intg$ and $H_1(M,\partial M;\intg)\cong \{e\}$. Then there is a unique choice of $\mathfrak{s}\in \operatorname{Spin}^c(M)$ and the cover $\widebar{T}_{M,\mathfrak{s}}$ can be identified with the infinite cylinder $[-1/2,1/2]\times \real$ with basepoints $\{0\}\times (\intg+1/2)$, where $\{\pm 1/2\}\times \real$ are glued together. To be clear, we will draw the (lift of) curve in \begin{equation}\label{eq: plane with points}
    (\real\times \real)\backslash (\intg\times (\intg+\frac{1}{2})),
\end{equation}which is periodic in the first coordinate. We call the first coordinate \emph{horizontal} and the second coordinate \emph{vertical}, parameterized by the meridian of the knot $K$ and the longitude $\lambda$, respectively. 

Hence the height function is just the projection onto the second coordinate. The \emph{height} of an arc $L_{\mu,\widebar{\mathfrak{s}}}$ is used to denote the height of its midpoint for short. For $p/q\in \mathbb{Q}\cup\{1/0\}$ and $h\in \intg+(p-1)/2$, we write \begin{equation}\label{eq: arc}
    L_{p/q,h}
\end{equation} for the arc of slope $p/q$ and height $h$. Then the intersection points of the peg-board curve $\widehat{HF}(M)$ and $L_{p/q,h}$ contribute to $\widehat{HFK}(\widetilde{K}_{p/q},h)$, the knot Floer homology of the dual knot $\widetilde{K}_{p/q}\subset S^3_{p/q}(K)$ in Alexander grading $h$. Since $\widehat{HFK}(K)$ with the Alexander grading detects the genus of $K$, the curve intersects the vertical arc of height $g$ ({\it i.e.}, $L_{1/0,g}$) nontrivially and any vertical arc of larger height trivially. Moreover, the symmetry of the knot Floer chain complex for knots in $S^3$ implies that the whole curve invariant $\widehat{HF}(M)$ is symmetric under the rotation by $180$ degrees around $(0,0)$ (a single 
component may be asymmetric).

Since $\widehat{HF}(S^3)$ is $1$-dimensional, the gluing theorem says there is only one intersection point (even with multiplicity from local systems) between the peg-board diagram of $\widehat{HF}(M)$ and the vertical line $\{1/2\}\times \real$. By symmetry around $(0,0)$, the intersection point must be $\{1/2\}\times \{0\}$. Without loss of generality, we assume $\ga_0$ in (\ref{eq: multicurve inv}) is the component contributing to the intersection point, which must have a trivial local system by the lack of the dimension. We call $\ga_0$ the \emph{distinguished component} and other $\ga_i$ the \emph{acyclic} components. Note that acyclic components all lie in a small neighborhood of $\{0\}\times \real$. 

\begin{figure}[ht]
\centering
\begin{overpic}[width=1\textwidth]{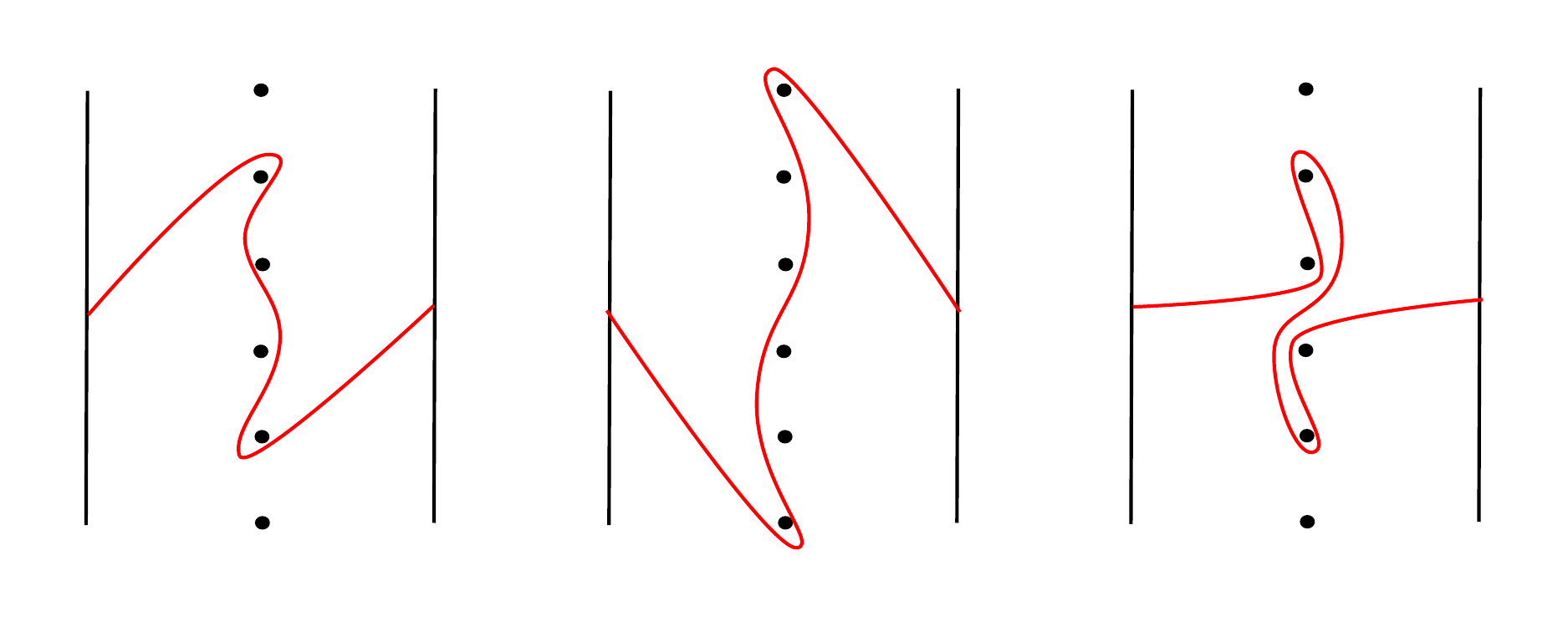}
	\put(11,0){$\tau = 2$ $\epsilon = 1$}
	\put(43,0){$\tau = -3$ $\epsilon = -1$}
	\put(78,0){$\tau = 0$ $\epsilon = -1$}
\end{overpic}
\vspace{-0.2in}
\caption{Examples of distinguished components}\label{fig: tau and ep}
\end{figure}

The name ``acyclic" comes from Hom's result about concordant knots \cite[Theorem 1]{Hom2017concordance}. From Hanselman-Watson's reinterpretation \cite[Proposition 2]{HW2023cable}, the distinguished component $\ga_0$ is a concordance invariant, and the $\tau$ invariant and the $\epsilon$ invariant in Heegaard Floer theory can be read from $\ga_0$ easily. Explicitly, suppose $\ga_0$ is oriented from left to right and starts at some point on $\{-1/2\}\times \real$. Then the integer $\tau(K)$ is the height of the first vertical arc in $\{0\}\times \real$ where $\ga_0$ meets, and $\epsilon(K)\in\{-1,1,0\}$ corresponds to the situation after the intersection with the first vertical arc. When $\ga_0$ turns downwards, upwards, or continues straight to wrap around the cylinder, the invariant becomes $-1,1,0$, respectively. Note that $\epsilon =0$ only if $\tau=0$, where the curve is a horizontal line. Some examples of curves and corresponding invariants (which are not necessarily from $\widehat{HF}(M)$) can be found in Figure \ref{fig: tau and ep}.

\subsection{First differentials on local maxima}\label{subsec: d1 curve}
In this subsection, we use the setup in \S \ref{subsec: curve inv} to study the first differential $\Psi$ and $\Phi$ for the dual knot $\widetilde{K}_{p/q}$ of a knot $K\subset S^3$. The main reason to focus on $S^3$ is that $\widehat{HF}(S^3)$ is $1$-dimensional, so many results can be generalized to knots in other closed $3$-manifolds, with extra care.

We consider the infinite cylinder $[-1/2,1/2]\times \real$ with basepoints $\{0\}\times (\intg+1/2)$. To be clear, we lift it periodically to the plane with punctures (\ref{eq: plane with points}) but still consider intersections and bigons are counted in the infinite cylinder. A \emph{compact curve} is an immersed curve with no endpoints on punctures. A \emph{non-compact curve} is an immersed curve with endpoints on punctures. We always assume the curves are in the peg-board diagram mentioned in \S \ref{subsec: d1 curve}, {\it i.e.}, there are no trivial bigons ({\it i.e.}, bigon containing no basepoints) between two curves. Recall that curves in $\widehat{HF}(M)$ in (\ref{eq: multicurve inv}) are all compact and the arc $L_{p/q,h}$ of slope $p/q$ and height $h$ in (\ref{eq: arc}) is a non-compact curve. Though many results hold for general non-compact curves, we only focus on arcs of shape $L_{p/q,h}$ for simplicity. We specify the name ``arc" only for some $L_{p/q,h}$.

We start with the following obvious lemma for curves.
\blem\label{lem: must have intersection}
Suppose $L_1,L_2,L_3$ are three arcs that form a triangle in the plane (not in the infinite cylinder) with basepoints, {\it i.e.}, any two of them share a common endpoint. Suppose the interior of the triangle does not contain any basepoints. Suppose $\ga$ is a compact curve in the plane with any orientation. If $\ga$ intersects $L_1$ at a point $a$ with the orientation pointing into the triangle, then the next intersection point $b$ with the triangle $L_1\cup L_2\cup L_3$ along the orientation must exist and lie on $L_2\cup L_3$.
\elem
\bpf
The algebraic intersection number of $\ga$ and $L_1\cup L_2\cup L_3$ is zero, so the number of intersection points is even. Hence there must be another intersection point other than $a$. For the next intersection point $b$, if $b$ is also on $L_1$, then there is a trivial bigon between $a$ and $b$, which contradicts the assumption of the peg-board diagram.
\epf
\brem\label{rem: two lifts}
Note that if we project the triangle onto the infinite cylinder, then the interior of the triangle is not well-defined. Indeed, a triangle on the infinite cylinder can be lifted into two different triangles in the plane modulo the horizontal translation. This phenomenon corresponds to the two different bypass triangles for three fixed sutured instanton homologies illustrated in \cite[Remarks 4.15 and 4.38]{LY2020}. Note that in the reference, the authors consider the manifolds with opposite orientations so all arrows between arcs are clockwise with respect to basepoints, while in this paper there is no orientation reversal so all arrows are counterclockwise, which is compatible with the usual orientations about bigons in Heegaard Floer theory. Though not used, we always write the arcs in the order such that the arrows between them are counterclockwise. Hence one may think of an exact triangle between the corresponding knot Floer homologies. The correspondence between arrows and bypass maps in Heegaard Floer theory may be proved by translating the computations of bordered sutured Floer homology by Etnyre-Vela-Vick-Zarev \cite[\S 6]{etnyre2017sutured} to the language of curve invariants. However, we do not need such correspondence in this paper, so we leave the verification to the readers.
\erem

Recall that an immersed curve has a well-defined tangent line at any point. We define a \emph{local maximum} of a compact curve to be a point where the tangent line is horizontal ({\it i.e.}, constant on the second coordinate) and above the curve locally ({\it i.e.}, in a neighborhood of the point). Similarly, a \emph{local minimum} of a compact curve is a point where the tangent line is horizontal and below the curve locally. We set that a horizontal line (corresponding to $\epsilon=\tau=0$) has no local maximum and minimum.

Note that we can isotope the curve to create or cancel pairs of a local maximum and a local minimum, but no cancellation can happen in a peg-board diagram otherwise there are two trivial bigons between the curve and the horizontal tangent line.

In the peg-board diagram, a neighborhood of local maximum or minimum must contain an intersection point with a vertical arc $L_{1/0,h}$. We also regard the height $h\in\intg$ as the height of the local maximum or minimum. A local maximum is a \emph{global maximum} if its height is not less than the height of any other local maximum. Note that the global maximum of any component of curve must exist and may not be unique. We also define a \emph{global minimum} similarly. Given a compact curve $\ga$, we write $n_+(\ga,h)$ for the number of local maxima of height $h$ and write $n_-(\ga,h)$ for the number of local minima of height $h$. The definitions of $n_\pm(\ga,h)$ extend to a collection of compact curves with local systems by counting the points with multiplicities from the ranks of the local systems.

\begin{figure}[ht]
\centering
\begin{overpic}[width=0.6\textwidth]{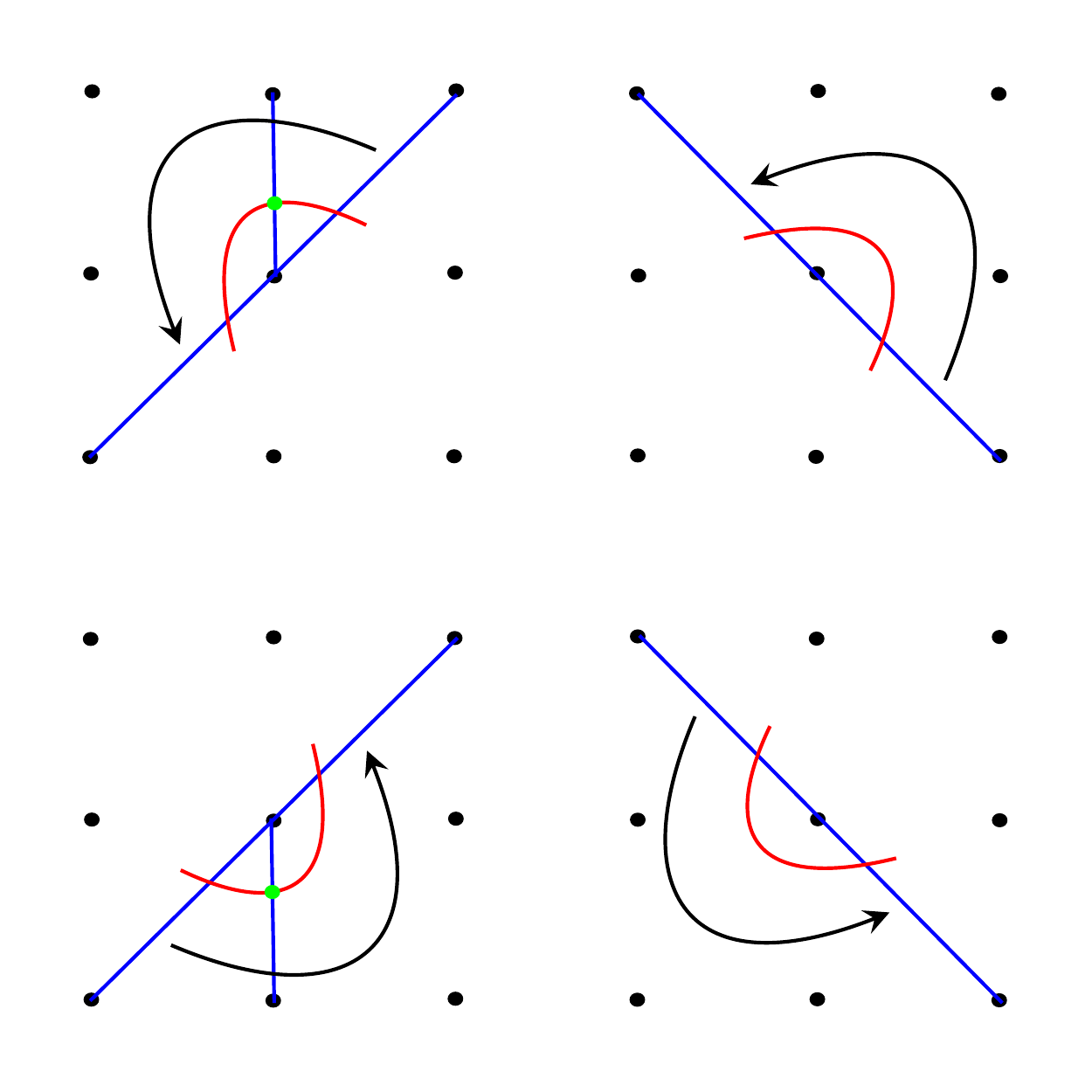}
	\put(23.5,3){(2)}
	\put(23.5,53){(1)}
	\put(73,3){(4)}
	\put(73,53){(3)}
	
	\put(10,83){$\Phi$}
	\put(26.5,82.5){\color{lygreen}$a$}
	\put(22,75){$z$}
	\put(26,72){$w$}
	\put(10,83){$\Phi$}
	
	\put(90,83){$\Psi$}
	\put(76,76){$w$}
	\put(72,72){$z$}
	
	\put(30,13){$\Psi$}
	\put(22,16){\color{lygreen}$a$}
	\put(76,26){$w$}
	\put(72,22.5){$z$}
	
	\put(65,17){$\Phi$}
	\put(22,25){$z$}
	\put(26,22){$w$}
\end{overpic}
\vspace{-0.2in}
\caption{local maxima and local minima}\label{fig: maxima}
\end{figure}
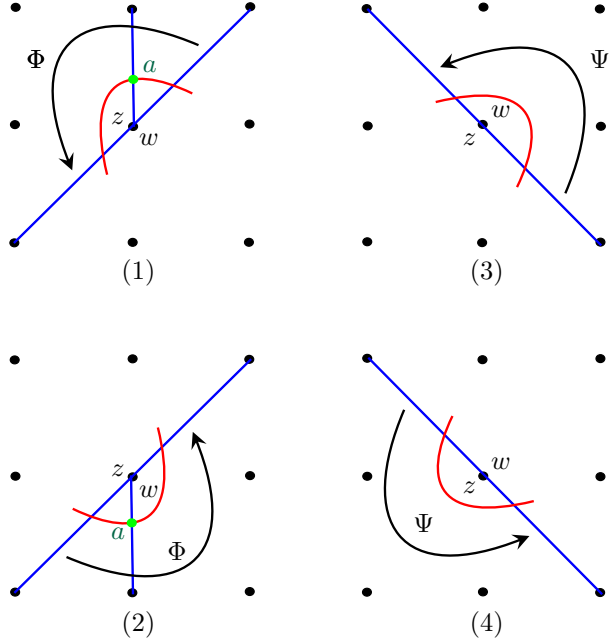

Now we state the main results of this subsection. Suppose $M=S^3\backslash   N(K)$ is a knot complement and $\widehat{HF}(M)$ is the curve invariant in (\ref{eq: multicurve inv}) drawn in a peg-board diagram of the infinite cylinder with basepoints. Suppose $\ga$ is a component of the compact curves in $\widehat{HF}(M)$. For the slope $p/q\in \mathbb{Q}$ with $p>0$ and the height $h\in \intg+(p-1)/2$, suppose $\widehat{HFK}(\widetilde{K}_{p/q},h)$ is the knot Floer homology of the dual knot $\widetilde{K}_{p/q}\subset S^3_{p/q}(K)$ in Alexander grading $h$ and suppose \begin{equation}\label{eq: phi psi}
\begin{aligned}
     \Psi_h:\widehat{HFK}(\widetilde{K}_{p/q},h)\to \widehat{HFK}(\widetilde{K}_{p/q},h+p)\\
     \Phi_h:\widehat{HFK}(\widetilde{K}_{p/q},h)\to \widehat{HFK}(\widetilde{K}_{p/q},h-p)
\end{aligned}
\end{equation} are differentials constructed by Sarkar \cite[\S 4]{sarkar15moving} corresponding to basepoints $w$ and $z$, respectively. We will omit the subscript $h$ if it is clear.
\bprop\label{prop: local maximum 1}
Suppose $p/q\in \mathbb{Q}$ with $p>0$ and $h\in \intg$. Suppose $\ga$ is an acyclic component ({\it i.e.}, $\ga\neq \ga_0$). Then we have the following cases related to Figure \ref{fig: maxima}.
 \benu 
    \item If $q>0$, then any local maximum $a$ of height $h$ corresponds to an element \begin{equation}\label{eq: xa 1}
     x_a\in \widehat{HFK}(\widetilde{K}_{p/q},h+\frac{-1+p}{2})\text{ such that } \Phi(x_a)\neq 0.
 \end{equation}The elements $x_{a_i}$ for all local maxima $a_i$ are linearly independent, so are their images under $\Phi$. In particular, we have $\rk \Phi_{h+\frac{-1+p}{2}}\ge n_+(\ga,h)$.
    \item If $q>0$, then any local minimum $a$ of height $h$ corresponds to an element \begin{equation}\label{eq: xa 2}
     x_a\in \widehat{HFK}(\widetilde{K}_{p/q},h+\frac{1-p}{2})\text{ such that } \Psi(x_b)\neq 0.
 \end{equation}There are similar linear independence results and we have $\rk \Psi_{h+\frac{1-p}{2}}\ge n_-(\ga,h).$
    \item If $q<0$, then any local maximum $a$ of height $h$ corresponds to an element \begin{equation}\label{eq: xa 3}
     x_a\in \widehat{HFK}(\widetilde{K}_{p/q},h+\frac{-1-p}{2})\text{ such that } \Psi(x_a)\neq 0.
 \end{equation}There are similar linear independence results and we have $\rk \Psi_{h+\frac{-1-p}{2}}\ge n_+(\ga,h).$
    \item If $q<0$, then any local minimum $a$ of height $h$ corresponds to an element \begin{equation}\label{eq: xa 4}
     x_a\in \widehat{HFK}(\widetilde{K}_{p/q},h+\frac{1+p}{2})\text{ such that } \Phi(x_a)\neq 0.
 \end{equation}There are similar linear independence results and we have $\rk \Psi_{h+\frac{1+p}{2}}\ge n_-(\ga,h).$
 \eenu 
\eprop
\bprop\label{prop: local maximum 2}
 Suppose $p/q\in \mathbb{Q}$ with $p>0$ and $h\in \intg$. Suppose $\ga$ is the distinguished component $\ga_0$. Then similar results as in Proposition \ref{prop: local maximum 1} hold with the following exception.
 \benu 
    \item If $\tau>0$, $\epsilon=1$ ({\it i.e.}, the part of $\ga$ disjoint from the vertical arc has a positive slope and it turns downwards after the first intersection point with the vertical arc from the left to the right) and $p/q>2\tau-1$, then the local maximum $a$ corresponds to $\tau$ does not correspond to an element in (\ref{eq: xa 1}).
    \item If $\tau>0$, $\epsilon=1$, and $p/q>2\tau-1$, then the local minimum $a$ corresponds to $-\tau$ does not correspond to an element in (\ref{eq: xa 2}).
    \item If $\tau<0$, $\epsilon=-1$, and $p/q<2\tau+1$, then the local maximum $a$ corresponds to $-\tau$ does not correspond to an element in (\ref{eq: xa 3}).
    \item If $\tau<0$, $\epsilon=-1$, and $p/q<2\tau+1$, then the local minimum $a$ corresponds to $\tau$ does not correspond to an element in (\ref{eq: xa 4}).
 \eenu 
 Here the local maximum and local minimum corresponding to $\tau$ or $-\tau$ is the one of height $\tau$ or $-\tau$ near the first or the last intersection point with the vertical arc from the left to the right, respectively.
\eprop
\bpf[Proof of Proposition \ref{prop: local maximum 1}]
The proofs of the four cases are similar, so we only prove the first case. Without loss of generality, we assume $a$ is an intersection point of $\ga$ and the vertical arc $L_{1/0,h}$. We consider the orientation of $\ga$ such that it goes from the left to the right of $L_{1/0,h}$ in the neighborhood of $a$. 

\begin{figure}[ht]
\centering
\begin{overpic}[width=0.85\textwidth]{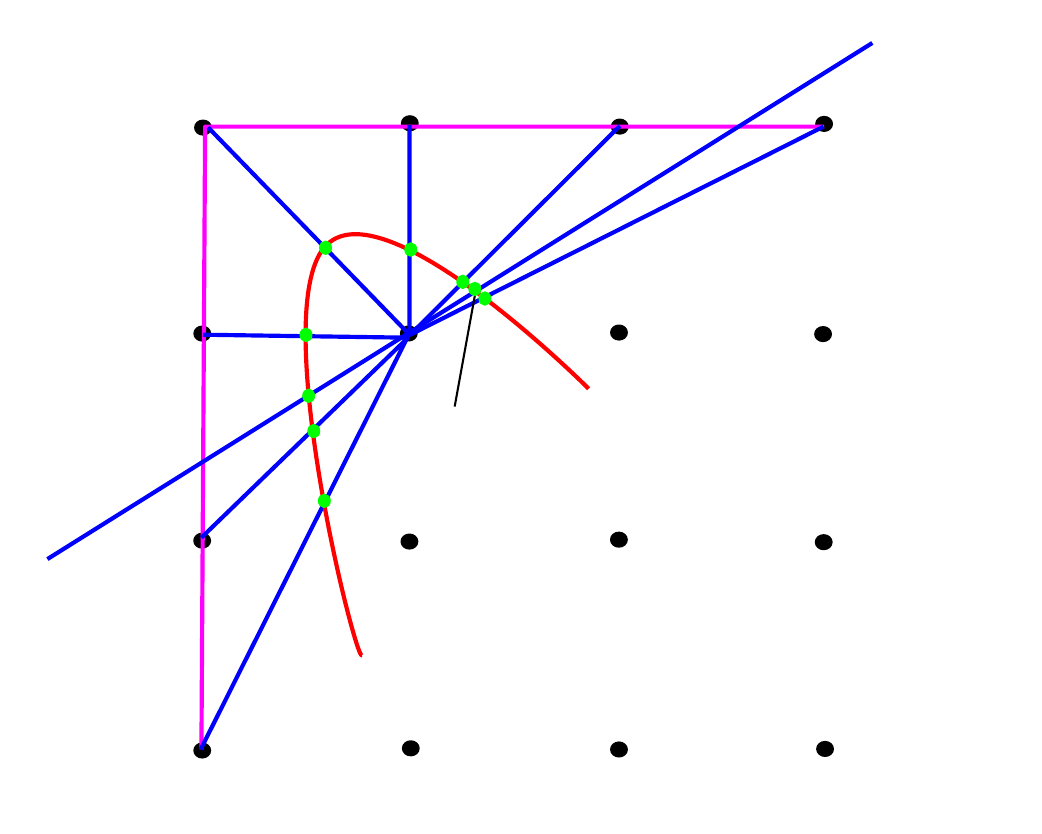}
	\put(45,70){\color{purple}$\widetilde{L}'_{0/1}$}
	\put(65,70){\color{purple}$\widetilde{L}'_{0/1}$}
	\put(15,57){\color{purple}$\widetilde{L}'_{1/0}$}
	
	\put(5,32){\color{blue}$\widetilde{L}_{p/q}$}
	\put(75,75){\color{blue}$\widetilde{L}_{p/q}$}
	
	\put(75,63){\color{blue}$\widetilde{L}_{1/2}$}
	\put(49,63){\color{blue}$\widetilde{L}_{1/1}$}
	\put(34,63){\color{blue}$\widetilde{L}_{1/0}$}
	\put(26,63){\color{blue}$\widetilde{L}_{1/(-1)}$}
	\put(21,49){\color{blue}$\widetilde{L}_{0/1}$}
	\put(23,29){\color{blue}$\widetilde{L}_{1/1}$}
	\put(23,11){\color{blue}$\widetilde{L}_{2/1}$}
	
	\put(32,31){\color{lygreen}$c_{2/1}$}
	\put(23.5,37){\color{lygreen}$c_{1/1}$}
	\put(25,43){\color{lygreen}$c_{p/q}$}
	\put(30,49){\color{lygreen}$c_{1/0}$}

	\put(23,56){\color{lygreen}$b_{1/(-1)}$}
	\put(37.5,57.5){\color{lygreen}$a$}
	\put(42,55){\color{lygreen}$b_{1/1}$}
	\put(48.5,50.5){\color{lygreen}$b_{1/2}$}
	\put(42,39){\color{lygreen}$b_{p/q}$}
\end{overpic}
\vspace{-0.2in}
\caption{Sequence of intersection points, where we omit the heights of the arcs}\label{fig: intersection}
\end{figure}
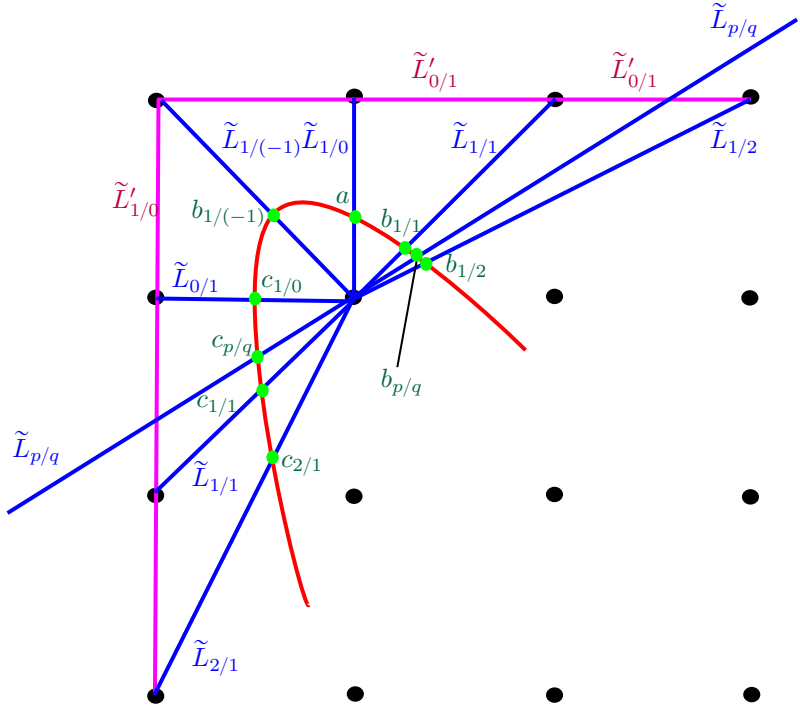

Then we apply Lemma \ref{lem: must have intersection} to the point $a$ and the lift of arcs $L_{1/0,h},L_{0/1,h+1/2},L_{1/1,h}$ and the lift $\widetilde{\ga}$ in Figure \ref{fig: intersection}, which we denote by $\widetilde{L}_{1/0,h},\widetilde{L}_{0/1,h+1/2},\widetilde{L}_{1/1,h}$, respectively. Note that the choices of lifts are important by Remark \ref{rem: two lifts} and indeed there are two different choices.

Since $a$ is a local maximum, the next intersection point $b$ from Lemma \ref{lem: must have intersection} cannot be on the lift $\widetilde{L}_{0/1,h+1/2}$. So it must be on the lift $\widetilde{L}_{1/1,h}$. We write this point as $b_{1/1}$. Then we apply Lemma \ref{lem: must have intersection} to the point $b_{1/1}$ and the lifts of arcs $L_{1/0,h},L_{0/1,h+1/2},L_{1/2,h}$ in Figure \ref{fig: intersection}. Note that the lift $\widetilde{L}^\p_{0/1,h+1/2}$ of $L_{0/1,h+1/2}$ in this triangle is different from the previous lift $\widetilde{L}_{0/1,h+1/2}$. Again from the fact that $a$ is a local maximum, the next intersection point cannot be on the lift $\widetilde{L}^\p_{0/1,h+1/2}$ and must be on the lift $\widetilde{L}_{1/2,h}$, which we write as $b_{1/2}$. By an obvious induction argument (using different lifts of $L_{0/1,h+1/2}$), we can find intersection points $b_{1/n}$ on the lifts $\widetilde{L}_{1/n,h}$ for all numbers $n\in \mathbb{N}$, where $b_{1/0}=a$. 

Since the curve $\widetilde{\ga}$ connects $b_{1/n}$ to $b_{1/{n+1}}$ and a lift $\widetilde{L}_{{p/q,h+(-1+p)/2}}$ of $L_{p/q,h+(-1+p)/2}$ for $q>0$ lies in between $\widetilde{L}_{1/n_0,h}$ and $\widetilde{L}_{1/{n_0+1},h}$ for some $n_0\in \mathbb{N}$, the curve also intersects $\widetilde{L}_{p/q,h+(-1+p)/2}$ at one point $b_{p/q}$ for $q>0$. The construction of points $b_{p/q}$ can also be generalized to the case $q<0$ by reserving the orientation of $\ga$.

From (\ref{eq: alexander}), for $q>0$, the point $b_{p/q}$ induces an element $x_a$ of the desired grading in (\ref{eq: xa 1}). Later we prove $\Phi(x_a)\neq 0$ using the curve interpretation of $\Phi$ after (\ref{eq: xa 1}). It is also worth mentioning that for $q<0$, the point $b_{p/q}$ is the element $\Psi(x_a)$ in the third case.

Now we reverse the orientation of $\ga$ such that it goes from the right to the left of $L_{1/0,h}$ in the neighborhood of $a$. Previously we already obtained the intersection point $b_{1/(-1)}$ between $\widetilde{\ga}$ and $\widetilde{L}_{1/(-1),h}$. Then we apply Lemma \ref{lem: must have intersection} to $b_{1/(-1)}$ and the lifts of the arcs $L_{1/(-1),h},L_{0/1,h-1/2},L_{1/0,h}$ in Figure \ref{fig: intersection}. Note that the lift $\widetilde{L}^\p_{1/0,h}$ in this triangle is different from the previous lift $\widetilde{L}_{1/0,h}$ of $L_{1/0,h}$ where $a$ lies, though the projections are the same. 

Since $\ga$ is acyclic, it lies in the neighborhood of the vertical arcs. Hence its lift $\widetilde{\ga}$ cannot intersect the lift $\widetilde{L}^\p_{1/0,h}$ in the above triangle. Hence, the next intersection point must be on the lift $\widetilde{L}_{0/1,h-1/2}$, which we denote by $c_{0/1}$. Similarly, by induction, we obtain intersection points $c_{n/1}$ on the lifts $\widetilde{L}_{n/1,h+(-1-n)/2}$ in Figure \ref{fig: intersection} for all number $n\ge \mathbb{N}$. The main fact we used is that $\widetilde{\ga}$ lies in the neighborhood of the vertical arcs and hence disjoint from the lifts $\widetilde{L}^\p_{1/0,h-n}$ in the corresponding triangles. Also, we can find intersection points $c_{p/q}$ on a lift $\widetilde{L}_{p/q,h+(-1-p)/2}$ for $q>0$ by inserting it in between two lifts $\widetilde{L}_{n/1,h+(-1-n)/2}$ for some $n=n_1$ and $n_1+1$.

In summary, the lift of the curve $\widetilde{\ga}$ with the second orientation intersects the lifts of arcs in a sequence of points \begin{equation}\label{eq: sequence}
    b_{p/q},b_{1/n_0},b_{1/(n_0-1)},\dots,b_{1/1},a=b_{0/1},b_{1/(-1)},c_{0/1},c_{1/1},\dots, c_{n_1/1},c_{p/q}.
\end{equation}Hence $c_{p/q}$ corresponds to the element $\Phi(x_a)\neq 0$ by the interpretation after (\ref{eq: xa 1}).

Finally, for different local maxima, the sequences of points in (\ref{eq: sequence}) lie in a disjoint part of the lift $\widetilde{\ga}$, so the corresponding elements are linearly independent.
\epf
\bpf[Proof of Proposition \ref{prop: local maximum 2}]
The proof of Proposition \ref{prop: local maximum 1} applies without change. Note that the existence of $c_{p/q}$ now depends on the fact that $p/q\le 2\tau -1$.
\epf

Now an analog of Corollary \ref{cor: dual injective} is obvious.
\bcor
Suppose $K \subset S^3$ is a nontrivial knot of genus $g$ and let $\widetilde{K}_1$ be the dual knot in $S_1^3(K)$. Then the maps $\Phi$ on $\widehat{HFK}(\widetilde{K}_1,g)$ and the map $\Psi$ on $\widehat{HFK}(\widetilde{K}_1,-g)$, both from (\ref{eq: phi psi}), are injective.
\ecor
\bpf
We simply apply Proposition \ref{prop: local maximum 1} and Proposition \ref{prop: local maximum 2} to global maxima and minima (possibly with multiplicities from local systems). Note that for nontrivial knot $K$, we have $g\ge 1$. If the first exception in Proposition \ref{prop: local maximum 2} happens, then $\tau=g$ and $p/q>2g-1$, which is impossible for $p/q=1/1$. The other exceptions do not happen for similar reasons.
\epf

Indeed we obtain stronger corollaries based on similar reasons.
\bcor\label{cor: main results}
Suppose $K\subset S^3$ is a nontrivial knot of genus $g$ and let $\widetilde{K}_{p/q}$ be the dual knot in $S_{p/q}^3(K)$ with $p/q\in \mathbb{Q}$ and $p>0$. Suppose $M=S^3\backslash N(K)$. Then we have the following.
\benu
    \item The map $\Phi$ in (\ref{eq: phi psi}) on \[\begin{cases}
        \widehat{HFK}(\widetilde{K}_{p/q},g+\frac{-1+p}{2})&\mathrm{for}~q>0\\
        \widehat{HFK}(\widetilde{K}_{p/q},-g+\frac{1+p}{2})&\mathrm{for}~q<0
    \end{cases}\] and the map $\Psi$ in (\ref{eq: phi psi}) on \[\begin{cases}
        \widehat{HFK}(\widetilde{K}_{p/q},-g+\frac{1-p}{2})&\mathrm{for}~q>0\\ \widehat{HFK}(\widetilde{K}_{p/q},g+\frac{-1-p}{2})&\mathrm{for}~q<0
    \end{cases}\] have at most $1$-dimensional kernel. If the kernel is nontrivial, then $\tau\in \{g,-g\}$ and the exceptions in Proposition \ref{prop: local maximum 2} happen.
    \item The total ranks of $\Phi$ and $\Psi$ on $\widehat{HFK}(\widetilde{K}_{p/q})$ are at least the number $n$ of acyclic components in the multicurve invariant $\widehat{HF}(M)$ in (\ref{eq: multicurve inv}). If the exceptions in Proposition \ref{prop: local maximum 2} do not happen ({\it e.g.}, $p/q=1 /(\pm1)$), then the ranks are at least the number $n+1$ of all components.
    \item If the total rank of either $\Phi$ or $\Psi$ on $\widehat{HFK}(\widetilde{K}_{p/q})$ is at most one, then there are no acyclic components in $\widehat{HF}(M)$. In such case, we know $\widehat{HFK}(K)$ has the minimal rank among all knots concordant to $K$. Moreover, if either $\Phi$ or $\Psi$ vanishes, then there is one local maximum or one local minimum in the distinguished component, and the exceptions in Proposition \ref{prop: local maximum 2} happen with $\tau=\pm g$. In such a case, the knot $K$ is an L-space knot.
\eenu
\ecor
\bpf
The first two arguments are from Proposition \ref{prop: local maximum 1} and Proposition \ref{prop: local maximum 2} by focusing on global maxima and global minima. Note that each component has at least one global maximum and one global minimum since it is compact. For the third argument, the acyclic components are ruled out by the second argument. The concordance argument follows from the fact that the distinguished component $\ga_0$ is a concordance invariant \cite[Proposition 2]{HW2023cable}.  Similarly, if either $\Phi$ or $\Psi$ vanishes, the global maximum or the global minimum must satisfy the exceptions in Proposition \ref{prop: local maximum 2}, which should also be the unique local maximum or local minimum. In such a case, if we pair the curve either with a compact line of slope $p^\p/q^\p$, or the arc $L_{p^\p/q^\p}$ with \[p^\p/q^\p\begin{cases}
    >2\tau-1&\tau=g>0;\\
    \neq 0 &\tau=g=0;\\
    <2\tau+1&\tau=-g<0,
\end{cases}\]we will obtain $|p^\p|$ intersection points, which implies that $K$ is an L-space knot.
\epf

Now we prove the main result of this section.

\quad

\noindent
{\bf Theorem \ref{thm: Heegaard dual simple}.} {\it Suppose $K\subset S^3$. If there exists $p/q\in \mathbb{Q}\backslash\{0\}$ satisfying \begin{equation*}\label{eq: HF simple 2}
    \dim \widehat{HFK}(S^3_{p/q}(K),\widetilde{K}_{p/q};\ft)=\dim \widehat{HF}(S^3_{p/q}(K);\ft),
\end{equation*}then both dimensions must equal to $|H_1(S^3_{p/q}(K);\intg)|=|p|$ and hence $K$ is a (Heegaard Floer) L-space knot. Consequently, we know $K$ is fibered and strongly quasi-positive \cite{Ozsvath2005,ni2007knot,hedden2010positivity}, and $|p/q|>2g(K)-1$ \cite[Corollary 3.6]{Rasmussen2017}.
}
\bpf
Recall that $\Phi$ or $\Psi$ is the differential in the first page of the spectral sequence from $\widehat{HFK}(Y,K;\ft)$ to $\widehat{HF}(Y;\ft)$ for rationally null-homologous knot $K\subset Y$. Then the theorem follows from the third argument in Corollary \ref{cor: main results} where $\Phi$ or $\Psi$ vanishes. The fact that $S^3_{p/q}(K)$ is an L-space is from \cite[Corollary 3.6]{Rasmussen2017}.
\epf

\appendix
\section{Closed \texorpdfstring{$3$}{3}-manifolds without \texorpdfstring{$2$}{2}-torsion}\label{sec: examples for no 2-torsion}In this appendix, we provide a few families of closed $3$-manifolds $Y$ such that $I^{\sharp}(Y;\intg)$ has no $2$-torsion. We make this collection for the completeness of the current paper, but we believe all of these examples are known to experts and we claim no originality for the results presented here.

Scaduto \cite[Theorem 1.1]{scaduto2015instanton} constructed a spectral sequence from the reduced odd Khovanov homology of a knot $K$ to the framed instanton Floer homology $\ish{\Sigma_2(\widebar{K})}$ of the double branched cover of the mirror $\widebar{K}$. As a corollary, when $K$ is a quasi-alternating link, \cite[Corollary 1.2]{scaduto2015instanton} says that $\ish{\Sigma_2(\widebar{K});\intg}$ is free over $\intg$ and has minimal possible rank, and hence has no $2$-torsion. In particular, the double-branched cover of a $2$-bridge knot is a lens space $L(p,q)$, so $\ish{L(p,q);\intg}$ has no $2$-torsion.

Additionally, Agol \cite{agol2023chainmail} constructed examples of L-spaces coming from surgeries along chainmail links. Since his argument only used the results for lens spaces and surgery exact triangles ({\it cf.} Lemma \ref{lem: surgery exact triangle, BS}), it applies verbatim to $\ft$ coefficients and hence we have the following proposition.
\bprop
Suppose $Y$ is an $L$-space obtained from surgery with large enough slopes on negative alternating chainmail links and partially augmented negative alternating chainmail links that arise from the proofs of \cite[Theorem 3.2 and Theorem 4.2]{agol2023chainmail}, then
\[
\dim\ish{Y;\ft}=\dim \ish{Y;\comp}= |H_1(Y;\intg)|
,\]and $\ish{Y;\intg}$ has no $2$-torsion.
\eprop
\brem
Agol showed that large enough surgery slopes on these links yield L-spaces, for the technical reason of allowing induction to work. The above lemma works exactly for those L-spaces. There might be the possibility that small surgery slopes of these links also generate an L-space but we know nothing about the $I^{\sharp}$ of such $3$-manifolds.
\erem

Moreover, we can adapt the proof of \cite[Proposition 4.12]{baldwin2018stein} over $\ft$ since it only uses the surgery exact triangle. We obtain the following result.

\bprop\label{prop: slope}
Suppose $K\subset S^3$ is a knot and $n$ is a positive integer. If\[\dim\ish{S_n^3(K);\ft}=\dim \ish{S_n^3(K);\comp}= |H_1(S_n^3(K);\intg)|=n,\]then for all rational number $p/q\ge n$, we have \[\dim\ish{S_{p/q}^3(K);\ft}=\dim \ish{S_{p/q}^3(K);\comp}= |H_1(S_{p/q}^3(K);\intg)|=|p|,\]and $\ish{S_{p/q}^3(K);\intg}$ has no $2$-torsion.
\eprop

\brem
Other results in \cite[\S 4.2]{baldwin2018stein} do not apply over $\ft$ directly since they are based on the adjunction inequality \cite[Remark 4.14]{baldwin2018stein}, which is only available over $\comp$ as mentioned in the introduction. A natural question from Proposition \ref{prop: slope} is to find the minimal positive real number $r_0$ such that for all slope $r>r_0$, we have $I^\sharp$ of the corresponding manifolds do not have $2$-torsion. From Lemma \ref{lem: false statement for Poincare sphere} and Proposition \ref{prop: wrong for trefoil}, we know $r_0\in (1,5]$ for the trefoil knot but cannot decide it by current techniques, which already becomes different from the case over $\comp$.
\erem

\section{\texorpdfstring{$2$}{2}-torsion for alternating knots}\label{sec: alternating knot}
It is well-known to experts that $\ish{S^3,K}$ of an alternating knot $K$ has $2$-torsion. However, the authors have not found a complete computation over $\intg$ in the literature so we do this in this section. We do not claim originality for these results.

\blem\label{lem: rank inequal}
For any link $L$ in a closed $3$-manifold $Y$, we have $\dim \ish{Y,K;\comp}\ge \dim \ina{Y,K;\comp}$.
\elem
\bpf
From an unoriented skein relation as in \cite[\S 3]{xie2021earring}, we have the following exact triangle
\begin{equation*}
	\xymatrix{
	\ish{Y,L;\comp}\ar[rr]^{}&&\ish{Y,L;\comp}\ar[dl]\\
	&\ish{Y,L^\natural;\comp}\ar[lu]&
	}
\end{equation*}
From \cite[Proposition 5.2]{xie2021earring} and the connected sum formula for sutured Floer homology (for example, an instanton version of \cite[Corollary 4.13]{li2018contact}), we obtain \[\ish{Y,L^\natural;\comp}\cong I(Y,L^\natural;\comp)\ot \comp^2=\ina{Y,L,\comp}\ot\comp^2.\]Then the inequality follows from the dimension inequality from the triangle.
\epf

\bthm\label{thm: thin knot has 2 torsion}
Suppose $K$ is a quasi-alternating knot (in particular, a non-split alternating knot) \cite{MO2008khthin}. Then the spectral sequence from Khovanov homology to instanton homology collapses, and we have \[Kh(\widebar{K};\intg)\cong \ish{S^3,K;\intg}\cong \intg^{\delta+1}\op (\intg/2)^{(\de-1)/2},\]
\[\aand \widebar{Kh}(\widebar{K};\intg)\cong \ina{S^3,K;\intg}\cong \intg^{\de},\]where $\de=|\det(K)|=|\Delta_K(-1)|$, $Kh$ and $\widebar{Kh}$ are the unreduced and the reduced Khovanov homologies. For the reduced version, the results also hold for quasi-alternating links.
\ethm

\bpf
Note that quasi-alternating links are thin for $Kh$ and $\widebar{Kh}$ over any coefficients. By the universal coefficient theorem, it is also $H$-slim in the sense of \cite[1.A]{Shu21torsion}. By Lee's and Shumakovitch's work \cite{LEE02,LEE05,Shu14torsion,Shu21torsion} (in particular, \cite[1.E, 1.G, and 1.I]{Shu21torsion}), for an alternating link $K$, we have \[\dim Kh(K;\comp)=\de+1\aand \dim  Kh(K;\ft)=2\de.\]Moreover, over $\intg$ the only possible torsion summand is $\intg/2$. Hence \[Kh(K;\intg)\cong \intg^{\delta+1}\op (\intg/2)^{(\de-1)/2}.\]Also from \cite{LEE02}, we know $\widebar{Kh}(K)\cong \intg^\de$. Note that $\Delta_K(t)$ does not depend on the mirror, and neither does $\de$. Kronheimer-Mrowka \cite[Proposition 1.2 and Theorem 8.2]{kronheimer2011khovanov}  showed that there are spectral sequences from $Kh(\widebar{K})$ and $\widebar{Kh}(\widebar{K})$ to $\ish{S^3,K}$ and $\ina{S^3,K}$ over $\intg$, respectively, which provide upper bounds on $I^\sharp$ and $I^\natural$. 

On the other hand, Kronheimer-Mrowka's results $\chi(KHI(K))=-\Delta_K(t)$ \cite{kronheimer2010instanton}, together with the second equation in (\ref{eq: 2dim}), provide a lower bound on $I^\natural$ by the coefficients of $\Delta_K(t)$. Since $\Delta_K(t)$ for a quasi-alternating link $K$ has alternating coefficients, the lower bound matches with the upper bound $\de$ (over $\comp$ coefficients). Since $\widebar{Kh}$ is free, the spectral sequence about $I^\natural$ collapses and we conclude that the associated graded homology of $I^\natural$ is $\intg^{\delta}$. Since it is free, we also have $I^\natural\cong \intg^{\delta}$. 

Now we assume $K$ is a quasi-alternating knot. From the first equation in (\ref{eq: 2dim}) and the universal coefficient theorem, we know $\dim \ish{S^3,K;\comp}$ is even. From Lemma \ref{lem: rank inequal} and the fact that $\delta$ is odd for knots, we know \[\dim \ish{S^3,K;\comp}\ge \de+1.\]Hence the spectral sequence from $Kh(\widebar{K})$ to $\ish{S^3,K}$ collapses over $\comp$. Hence the differential in the spectral sequence can only have a torsion image. On the other hand, the spectral sequence collapses over $\ft$ by (\ref{eq: 2dim}) and the result for $I^{\natural}$. Since $Kh$ has only $\intg/2$ torsion summand, the differential vanishes over $\intg$, and the associated graded homology of $I^\sharp$ is $\intg^{\delta+1}\op (\intg/2)^{(\de-1)/2}$. Potentially $I^\sharp$ and its associated graded homology could be different, but that cannot happen since the spectral sequence collapses both over $\comp$ and $\ft$. Hence, we conclude the result for $I^\sharp$.
\epf

\bibliographystyle{alpha}

\end{document}